\documentclass[12pt]{article}
\usepackage{amssymb}
\usepackage{amsmath}
\usepackage{array}
\usepackage[russian]{babel}
\usepackage[pdftex]{graphics}
\usepackage[X2,T2A]{fontenc}
\usepackage[cp1251]{inputenc}
\usepackage{mathrsfs}
\usepackage{longtable}
\usepackage{mathrsfs}
\usepackage{textcase}
\usepackage{setspace}
\usepackage{dsfont}
\usepackage{bbm}
\usepackage{pgf,pgfarrows,pgfnodes,pgfautomata,pgfheaps,pgfshade}
\usepackage{pdflscape}
\usepackage{fancyhdr}
\DeclareMathOperator{\mmod}{mod}

\multlinegap=0em \DeclareFontFamily{T1}{msb}{}
\DeclareFontShape{T1}{msb}{m}{ol}{<5> <6> <7> <8> <9> gen * msbm
<10> <10.95> <12> <14.4> <17.28> <20.74> <24.88> msbm10}{}
\DeclareSymbolFont{AMSb}{T1}{msb}{m}{ol} \multlinegap=0em

\renewcommand{\S}{\mathhexbox278}

\begin{document}

\begin{center}
{\rmfamily\bfseries\Large A distribution related to Farey sequences - I}
\end{center}

\begin{center}
{\normalsize M.A.~Korolev\footnote{Steklov Mathematical Institute of Russian Academy of Sciences. E-mail: \texttt{korolevma@mi-ras.ru}, \texttt{hardy\_ramanujan@mail.ru}}}
\end{center}
\vspace{0.5cm}

\fontsize{11}{13pt}\selectfont

We study some arithmetical properties of Farey sequences by the method introduced by F.~Boca, C.~Cobeli and A.~Zaharescu (2001).
Let $\Phi_{Q}$ be the classical Farey sequence of order $Q$. Having the fixed integers $D\geqslant 2$ and $0\leqslant c\leqslant D-1$, we colour to the red the fractions in $\Phi_{Q}$ with denominators $\equiv c \; (\mmod D)$. Consider the gaps in $\Phi_{Q}$ with coloured endpoints, that do not contain the fractions $a/q$ with $q\equiv c\;(\mmod D)$ inside. The question is to find the limit proportions $\nu(r;D,c)$ (as $Q\to +\infty$) of such gaps with precisely $r$ fractions inside in the whole set of the gaps under considering ($r = 0,1,2,3,\ldots$).

In fact, the expression for this proportion can be derived from the general result obtained by C.~Cobeli, M.~V\^{a}j\^{a}itu and A.~Zaharescu (2014). However, such formula expresses $\nu(r;D,c)$ in the terms of areas of some polygons related to a special geometrical transform. In the present paper, we obtain an explicit formulas for $\nu(r;D,c)$ for the cases $D = 2, 3$ and $c=0$.

This paper is an extended version of the authors' paper <<A distribution related to Farey series>> published in <<Chebyshevski\v{i} Sbornik>> (2023).
\vspace{0.3cm}
\begin{center}
\textcolor{blue}{\textsc{Сontents}}
\end{center}

\textcolor{blue}{
\contentsline{subsection} {\numberline {\S\,1} Introduction}{\textcolor{blue}{1}}{}
\contentsline{subsection} {\numberline {\S\,2} Farey sequences and lattice points in polygons}{\textcolor{blue}{6}}{}
\contentsline{subsection} {\numberline {\S\,3} Some properties of the polygons $\mathcal{T}(k_{1},\ldots,k_{r})$}{\textcolor{blue}{17}}{}
\contentsline{subsection} {\numberline {\S\,4} Expression for $N(Q;r,D,c_{0})$ in terms of areas of polygons}{\textcolor{blue}{24}}{}
\contentsline{subsection} {\numberline {\S\,5} Precise formulas for the continuants and the proportion $\nu(Q;r,D,c_{0})$}{\textcolor{blue}{29}}{}
\contentsline{subsection} {\numberline {\S\,6} Theorem 1: case $1\leqslant r\leqslant 7$}{\textcolor{blue}{34}}{}
\contentsline{subsection} {\numberline {\S\,7} Explicit form of some polygons $\mathcal{T}(\mathbf{k})$}{\textcolor{blue}{41}}{}
\contentsline{subsection} {\numberline {\S\,8} Theorem 1: case $r\geqslant 8$}{\textcolor{blue}{67}}{}
\contentsline{subsection} {\numberline {\S\,9} The formula of three authors}{\textcolor{blue}{70}}{}
\contentsline{subsection} {\numberline {\S\,10} Theorem 2: the case $D = 2, c_{0}=0$}{\textcolor{blue}{74}}{}
\contentsline{subsection} {\numberline {\S\,10} Final remarks}{\textcolor{blue}{75}}{}
\contentsline{subsection} {\numberline {} Appendix I. Precise form of the polygons $\mathcal{T}(\mathbf{k}_{r})$}{\textcolor{blue}{77}}{}
\contentsline{subsection} {\numberline {} Appendix II. Table of the proportions $\nu(r;3,0), r\geqslant 8$}{\textcolor{blue}{85}}{}
\contentsline{subsection} {\numberline {} Appendix III. The sets $\mathcal{A}_{r}^{\circ}(3;0,1)$}{\textcolor{blue}{87}}{}
\contentsline{subsection} {\numberline {} References}{\textcolor{blue}{88}}{}
}
\fontsize{12}{15pt}\selectfont

\section{Introduction}

Suppose that $Q\geqslant 1$. Farey sequence $\Phi_{Q}$ of order $Q$ is the set of all non-negative irreducible subunitary fractions whose denominators do not exceed $Q$, arranged in the ascending order. That is,
\[
\Phi_{Q} = \biggl\{\frac{a}{q}\,:\,0\leqslant a\leqslant q\leqslant Q,\;\text{GCD}(a,q)=1\biggr\}.
\]
Farey sequences appeared in a wide class of number-theoretic problems, and thus have become the object of comprehensive research. The bibliography of Farey sequences is too wide; that's why we refer the readers to the survey \cite{Cobeli_Zaharescu_2003} (2003). The recent studies are partly mentioned in the reference lists of the papers \cite{Haynes_2003}-\cite{Boca_Siskaki_2022}.

The statistical laws connected with this sequence are of a great interest. As example, one can mention a class of problems concerning the distribution of the fractions of $\Phi_{Q}$ with denominators in a given arithmetic progression modulo $D\geqslant 2$ (see \cite{Haynes_2003}-\cite{Boca_Cobeli_Zaharescu_2001}). A typical problem of such type is the following. Let
$r\geqslant 1$  be any fixed integer and let $\overline{\mathbf{c}} = (c_{1},\ldots,c_{r})$ be the tuple of integers with the conditions $0\leqslant c_{1},\ldots, c_{r}\leqslant D-1$. Denote by
$N(Q)$ the cardinality of $\Phi_{Q}$ and by $N(Q;r,D,\overline{\mathbf{c}})$ the number of tuples of $r$ consecutive fractions of the sequence $\Phi_{Q}$ of the type
\[
\frac{a_{1}}{q_{1}} < \frac{a_{2}}{q_{2}} < \ldots < \frac{a_{r}}{q_{r}},
\]
where $q_{i}\equiv c_{i} \pmod{D}$, $i = 1,2,\ldots, r$. In 2012, C.~Cobeli, M.~V\^{a}j\^{a}itu and A.~Zaharescu \cite{Cobeli_Vajaitu_Zaharescu_2012} proved that the limit
\[
\varrho(r,D,\overline{\mathbf{c}}) = \lim_{Q\to +\infty}\frac{N(Q;r,D,\overline{\mathbf{c}})}{N(Q)}
\]
exists for any fixed $r$, $D$ and $\overline{\mathbf{c}}$. They also give a general formula for such limit\footnote{For the first time, this result was published in 2005; see \cite{Cobeli_Zaharescu_2005}. Here we use another notation: $\Phi_{Q}$ instead of $\mathfrak{F}_{Q}$,  $N(Q)$ instead of $\# \mathfrak{F}_{Q}$, $D$ instead of $\mathfrak{d}$, $N(Q;r,D,\overline{\mathbf{c}})$ instead of $\mathbf{N}_{Q}^{r}(\mathbf{c},\mathfrak{d})$ and $\varrho(r,D,\overline{\mathbf{c}})$ instead of $\rho^{r}(\mathbf{c},\mathfrak{d})$ (see \S 9).}.

The main tool for deriving such formula is the transformation $T$ of the triangle
\[
\mathcal{T} = \{(x,y): 0<x,y\leqslant 1,\,x+y>1\}
\]
defined as follows:
\[
T(x,y) = \biggl(y,\biggl[\frac{1+x}{y}\biggr]y-x\biggr).
\]
The triangle $\mathcal{T}$ is called sometimes as \textit{Farey triangle} and the map $T$ as \textit{BCZ-transform} (named after the initials of the authors of \cite{Boca_Cobeli_Zaharescu_2001} who introduced this map).

Given an integer $k\geqslant 1$, we define $\mathcal{T}(k)$ as the set of points $(x,y)$ of Farey triangle such that $[(1+x)/y] = k$ (see \cite{Boca_Cobeli_Zaharescu_2001}). Further, given integers $k_{1}, k_{2}, \ldots, k_{r}\geqslant 1$ we define the region
\[
\mathcal{T}(k_{1},k_{2},\ldots, k_{r}) = \mathcal{T}(k_{1})\cap T^{-1}\mathcal{T}(k_{2})\cap\ldots\cap T^{-(r-1)}\mathcal{T}(k_{r}),
\]
It turned out that such region is either convex polygon or empty set.

The aforementioned formula for $\varrho(r,D,\overline{\mathbf{c}})$ contains the sum (in general case, infinite) of the areas of the polygons $\mathcal{T}(k_{1},k_{2},\ldots, k_{r})$  where $\mathbf{k} = (k_{1},k_{2},\ldots, k_{r})$ runs through a special set. As a rule, for given $r, D$ and $\mathbf{c}$, it is not easy to derive a closed expression for the proportion $\varrho(r,D,\overline{\mathbf{c}})$ from such a general formula. However, some examples of such calculations are given in \cite{Cobeli_Vajaitu_Zaharescu_2012} (see also \S 10 of the present paper).

Denote by $N(Q;r,D,c_{0})$ the number of  $(r+2)$-tuples of consecutive fractions of the sequence $\Phi_{Q}$ of the type
\begin{equation}\label{lab_01}
\frac{a_{0}}{q_{0}} < \frac{a_{1}}{q_{1}} < \ldots < \frac{a_{r}}{q_{r}} < \frac{a_{r+1}}{q_{r+1}},
\end{equation}
where $q_{0}, q_{r+1}\equiv c_{0} \pmod{D}$ and $q_{i}\not\equiv c_{0}\pmod{D}$, $i=1,2,\ldots,r$.

In what follows, we extend the sequence $\Phi_{Q}$ by requiring that $a_{i+N}/q_{i+N} = a_{i}/q_{i}+1$ for any $i$, where $N = N(Q)$.

Hence, $N(Q;r,D,c_{0})$ takes into account the tuples (\ref{lab_01}) with the condition $a_{0}/q_{0}>0$ and the tuples containing $a_{i}/q_{i} = 1/1$.

Further, denote by $N(Q;D,c_{0})$ the number of elements in $\Phi_{Q}$ whose denominators are congruent to $c$ modulo $D$ and consider the proportion
\[
\nu(Q;r,D,c_{0}) = \frac{N(Q;r,D,c_{0})}{N(Q;D,c_{0})}.
\]
In view of the results of \cite{Cobeli_Vajaitu_Zaharescu_2012}, the quantity $\nu(Q;r,D,c_{0})$ tends to a limit $\nu(r,D,c_{0})$ as $Q$ grows. This limit coincides (up to some multiplier depending only on $D$) with the sum of values $\varrho(r+2,D,\overline{\mathbf{c}})$ over all the tuples $\overline{\mathbf{c}} = (c_{0},\ldots,c_{r+1})$ satisfying the conditions  $c_{i}\not\equiv c_{0}\pmod{D}$, $i=1,2,\ldots, r$, $c_{r+1} = c_{0}$, $\text{GCD}(c_{0},c_{1},D)=1$ (see \S 9).

In the case $D=3$, $c_{0}=0$ and $Q = 3\cdot 10^5$, when  $N(Q;3,0) = 68\,395\,970$, the calculations lead us to the following approximate values of the proportions $\nu(Q;r,3,0)$ for small $r$:

\begin{longtable}
{|c|c|c|c|c|c|}
\hline $r$ & $N(Q;r,3,0)$   & $\nu(Q;r,3,0)$       & $r$   & $N(Q;r,3,0)$ & $\nu(Q;r,3,0)$ \\
\hline 1   & 11\,982\,989   & 0.17520\,02201       & 8     & 390\,976     & 0.00571\,63602 \\
\hline 2   & 25\,653\,970   & 0.37506\,84726       & 9     & 347\,028     & 0.00507\,38077 \\
\hline 3   & 14\,259\,027   & 0.20847\,75901       & 10    & 236\,232     & 0.00345\,38877 \\
\hline 4   &  6\,293\,540   & 0.09201\,62401       & 11    & 171\,556     & 0.00250\,82764 \\
\hline 5   &  4\,843\,722   & 0.07081\,88216       & 12    & 127\,068     & 0.00185\,78288 \\
\hline 6   &  2\,482\,708   & 0.03629\,90393       & 13    &  82\,264     & 0.00120\,27609 \\
\hline 7   &     953\,418   & 0.01393\,96810       & 14    &  83\,750     & 0.00122\,44874 \\
\hline
\end{longtable}
The aim of the present paper is to find an explicit expressions for $\nu(r) = \nu(r;3,0)$, $r = 1,2,3,\ldots$. Our main assertions is the following:
\vspace{0.2cm}

\textsc{Theorem 1.} \emph{Let $Q\to +\infty$. Then, for $r\geqslant 1$, $r = o((Q/\ln{Q})^{1/3})$, the following asymptotic formula holds:}
\[
\nu(Q;r,3,0) = \nu(r,3,0) + O\biggl(\frac{\ln{Q}}{Q}\biggr),
\]
\emph{where the implied constant is absolute. For $1\leqslant r\leqslant 7$, the values $\nu(r) = \nu(r;3,0)$ are listed below:}
\begin{align*}
& \nu(1) = 6-2\biggl(\frac{\pi}{\sqrt{3}}+\ln{3}\biggr) = 0.17517\,66941\ldots, \\[6pt]
& \nu(2) = 4\biggl(\frac{\pi}{\sqrt{3}}-\ln{3}\biggr)-\frac{87}{35} = 0.37503\,40165\ldots,\\
& \nu(3) = 12\ln{3}-\frac{53\,132}{4095} = 0.20850\,00891\ldots, \\
& \nu(4) = \frac{528\,904}{45\,045}-4\biggl(\frac{\pi}{\sqrt{3}}+\ln{3}\biggr) = 0.09203\,39300\ldots, \\
& \nu(5) = 2\biggl(\frac{\pi}{\sqrt{3}}-\ln{3}\biggr)-\frac{4\,164\,383}{3\,063\,060} = 0.07082\,42239\ldots, \\
& \nu(6) = \frac{3\,089}{85\,085} = 0.0363048715989892\ldots ,\\
& \nu(7) = \frac{54\,097}{3\,879\,876}= 0.01394\,29713\ldots\;.
\end{align*}
\emph{If $r\geqslant 8$, then, setting $r = 5m+i$ where $m\geqslant 1$, $0\leqslant i\leqslant 4$, one has}
\[
\nu(r) = \sum\limits_{j=1}^{n}P_{ij}(m),\quad \mbox{\emph{where}}\quad n = n(i) =
\begin{cases}
  3, & \mbox{\emph{for} } i=3, \\
  2, & \mbox{\emph{otherwise}},
\end{cases}
\]
\emph{and $P_{ij}(m)$ are rational functions defined as follows:}
\begin{align*}
& P_{01}(m)  = \frac{6(8m-1)}{(3m-1)(6m-1)(12m-1)(12m+1)}, \quad P_{02}(m)  = \frac{2}{(6m-1)(6m+1)(12m-1)},\\
& P_{11}(m)  = \frac{6(8m+1)}{(3m+1)(6m+1)(12m-1)(12m+1)}, \quad P_{12}(m)  = \frac{2}{(6m-1)(6m+1)(12m+1)}, \\
& P_{21}(m)  = \frac{6(4m+1)}{(3m+1)(6m+1)(12m+1)(12m+5)}, \\
& \qquad \qquad \qquad P_{22}(m)  = \frac{2(9m+4)}{3(2m+1)(3m+1)(6m+1)(12m+7)}, \\
& P_{31}(m)  = \frac{6(2m+1)}{(3m+1)(3m+2)(12m+5)(12m+7)}, \quad P_{32}(m)  = \frac{2}{3(2m+1)(6m+1)(12m+5)},\\
& \qquad \qquad \qquad P_{33}(m)  = \frac{2}{3(2m+1)(6m+5)(12m+7)}, \\
& P_{41}(m)  = \frac{6(4m+3)}{(3m+2)(6m+5)(12m+7)(12m+11)}, \\
& \qquad \qquad \qquad P_{42}(m)  = \frac{2(9m+5)}{3(2m+1)(3m+2)(6m+5)(12m+5)}.
\end{align*}

\textsc{Remark 1.} Theorem 1 implies that all the proportions $\nu(r)$ are rational for $r\geqslant 6$. At the same time, the numbers $\nu(r)$, $1\leqslant r\leqslant 5$ are transcendental.
Namely, if $r = 1,2,4,5$ then $\nu(r)$ are linear combinations with algebraic coefficients of the numbers $1$, $\pi i$ and $\ln{3}$, and its transcendence follows from the classical results of A.~Baker (see \cite[theorems 2.1, 2.2]{Baker_1975}). The transcendence of $\nu(3)$ is a direct corollary of Lindemann's theorem (see \cite{Lindeman_1882}).

\textsc{Remark 2.} One can easily see that $\nu(2)>\nu(3)>\nu(1)$. Next, we also have that
\begin{align*}
& \nu(14) = \tfrac{3\,931}{3\,209\,430} = 0.00122\,48280\ldots > \nu(13) = \tfrac{100\,349}{83\,445\,180} = 0.00120\,25739\ldots, \\[6pt]
& \nu(19) = \tfrac{206\,834}{440\,240\,493} = 0.00046\,98204\,\ldots > \nu(18) = \tfrac{387\,197}{889\,847\,805} = 0.00043\,51272\,\ldots, \\[6pt]
& \nu(29) = \tfrac{154\,342}{1\,214\,248\,035} = 0.00012\,71091\,\ldots > \nu(30) = \tfrac{12\,916}{114\,103\,745} = 0.00011\,31952\,\ldots > \\[6pt]
& > \nu(28) = \tfrac{473\,497}{4\,241\,317\,080} = 0.00011\,16391\,\ldots
\end{align*}
and so on (see Appendix II). Thus, it appears that the proportion $\nu(r)$ does not decrease monotonically when $r$ grows. At the same time, the monotonicity of $\nu(r)$ takes place inside any progression $r\pmod{5}$.
\vspace{0.3cm}

The last section of \cite{Cobeli_Vajaitu_Zaharescu_2012} contains the examples of the precise calculation of the values $\varrho(r,D,\overline{\mathbf{c}})$ for some vectors $\overline{\mathbf{c}}$ in the case $D = 2$. In particular, the authors of \cite{Cobeli_Vajaitu_Zaharescu_2012} give the precise values of $\varrho(r,2,\overline{\mathbf{c}}_{r})$ for $r = 3,4,5$ and for $\overline{\mathbf{c}}_{r} = (0,\underbrace{1,\ldots,1}_{r-2},0)$, namely
\[
\varrho(3,2,\overline{\mathbf{c}}_{3}) = 2-\frac{8}{3}\ln{2},\quad
\varrho(4,2,\overline{\mathbf{c}}_{4}) = \frac{16}{3}\ln{2} - \frac{1132}{315}, \quad \varrho(5,2,\overline{\mathbf{c}}_{5}) = \frac{599}{315}-\frac{8}{3}\ln{2}.
\]
The multiplication to the factor
\[
D\prod_{p|D}\biggl(1+\frac{1}{p}\biggr) = 3
\]
(for details, see \S 9) yields to the values of the proportions $\nu(r;2,0)$:
\begin{align*}
& \nu(1;2,0) = 6-8\ln{2} = 0.45482\,25555\ldots,\\[6pt]
& \nu(2;2,0) = 16\ln{2}-\frac{1132}{105} = 0.30940\,25080\ldots,\\[6pt]
& \nu(3;2,0) = \frac{599}{105}-8\ln{2} = 0.15958\,44603\ldots\;.
\end{align*}
The above arguments allow us to continue these calculations. In particular, the following assertion holds true.
\vspace{0.3cm}

\textsc{Theorem 2.} \emph{Let $Q\to +\infty$. Then, for $r\geqslant 1$, $r = o((Q/\ln{Q})^{1/3})$, the following asymptotic formula holds:}
\[
\nu(Q;r,2,0) = \nu(r,2,0) + O\biggl(\frac{\ln{Q}}{Q}\biggr),
\]
\emph{where the implied constant is absolute and}
\[
\nu(4;2,0) = \frac{2}{45}, \quad \nu(r;2,0) = \frac{8}{(2r-3)(2r-1)(2r+1)}, \quad\textit{for }\; r\geqslant 5.
\]
\emph{In particular,}
\begin{align*}
& \nu(5;2,0) = \frac{8}{693} = 0.01154\,40115\ldots && \nu(11;2,0) = \frac{8}{9\,177} = 0.00087\,17445\ldots \\[6pt]
& \nu(6;2,0) = \frac{8}{1\,287} = 0.00621\,60062\ldots && \nu(12;2,0) = \frac{8}{12\,075} = 0.00066\,25258\ldots \\[6pt]
& \nu(7;2,0) = \frac{8}{2\,145} = 0.00372\,96037\ldots && \nu(13;2,0) = \frac{8}{15\,525} = 0.00051\,52979\ldots \\[6pt]
& \nu(8;2,0) = \frac{8}{3\,315} = 0.00241\,32730\ldots && \nu(14;2,0) = \frac{8}{19\,575} = 0.00040\,86845\ldots \\[6pt]
& \nu(9;2,0) = \frac{8}{4\,845} = 0.00165\,11868\ldots && \nu(15;2,0) = \frac{8}{24\,273} = 0.00032\,95843\ldots \\[6pt]
& \nu(10;2,0) = \frac{8}{6\,783} = 0.00117\,94191\ldots && \nu(16;2,0) = \frac{8}{29\,667} = 0.00026\,96598\ldots \\[6pt]
\end{align*}

\textsc{Remark 2.} The above theorems appear for the first time in authors' paper \cite{Korolev_2023}. In this paper, the most part of auxilliary lemmas were proved for the particular case $D=3$, $c_{0}=0$. In this paper, we reprove such assertions for the general case. This will allow us to use these results in a forthcoming paper where we treat the case $D=3$, $c_{0}\ne 0$.
\vspace{0.3cm}

\textsc{Acknowledgements.} The author is grateful to Kirill S. Arzhanykh for carefull reading of the text and for many valuable remarks.
\vspace{0.3cm}

\section{Farey sequences and lattice points in polygons}

For the convenience of the reader, we give a sketch of the proof of general formula for $\nu(r;D,c)$ that follows to the method of \cite{Cobeli_Vajaitu_Zaharescu_2012}.
This also allows us to formulate all auxilliary assertions from \cite{Cobeli_Vajaitu_Zaharescu_2012} and its corollaries in the form most appropriated for the explicit calculation of $\nu(r;D,c)$.

Let $Q\geqslant 2$ and suppose that $a/q<a'/q'<a''/q''$ are consecutive fractions of the sequence $\Phi_{Q}$. Then the substraction of $a''q'-a'q''=1$ from the equality $a'q-aq'=1$ yields $(a+a'')q'=a'(q+q'')$. Since $a'$ and $q'$ are coprime, we get $q+q'' = kq'$ for some integer $k$. Since $kq'-q = q''\leqslant Q$, then $k\leqslant (q+Q)/q'$. At the same time, for neighbouring denominators we have $q'+q''>Q$ or, that is the same, $(k+1)q'-q>Q$. Hence we get $k+1>(q+Q)/q'$. Thus we obtain
\begin{equation}\label{lab_02}
k = \biggl[\frac{q+Q}{q'}\biggr] = \biggl[\frac{1+q/Q}{q'/Q}\biggr].
\end{equation}
The application of $BCZ$-transform to the point with coordinates $(x,y) = (q/Q,q'/Q)$ yields:
\[
T\biggl(\frac{q}{Q},\frac{q'}{Q}\biggr) = \biggl(\frac{q'}{Q},\biggl[\frac{1+q/Q}{q'/Q}\biggr]\frac{q'}{Q}-\frac{q}{Q}\biggr) = \biggl(\frac{q'}{Q},\frac{kq'-q}{Q}\biggr) = \biggl(\frac{q'}{Q},\frac{q''}{Q}\biggr).
\]
Therefore, if
\begin{equation}\label{lab_03}
\frac{a_{0}}{q_{0}} < \frac{a_{1}}{q_{1}} < \ldots < \frac{a_{r}}{q_{r}} < \frac{a_{r+1}}{q_{r+1}}
\end{equation}
is any tuple of consecutive elements of $\Phi_{Q}$ then
\begin{equation}\label{lab_04}
T^{2}\biggl(\frac{q_{0}}{Q},\frac{q_{1}}{Q}\biggr) = \biggl(\frac{q_{2}}{Q},\frac{q_{3}}{Q}\biggr),\quad T^{\,3}\biggl(\frac{q_{0}}{Q},\frac{q_{1}}{Q}\biggr) = \biggl(\frac{q_{3}}{Q},\frac{q_{4}}{Q}\biggr),\ldots, T^{\,r}\biggl(\frac{q_{0}}{Q},\frac{q_{1}}{Q}\biggr) = \biggl(\frac{q_{r}}{Q},\frac{q_{r+1}}{Q}\biggr).
\end{equation}
Setting $k_{i} = [(q_{i-1}+Q)/q_{i}]$ we get the chain of equations
\begin{equation}\label{lab_05}
\begin{cases}
q_{0} + q_{2} = k_{1}q_{1},\\
q_{1} + q_{3} = k_{2}q_{2},\\
\cdots\\
q_{r-1}+q_{r+1} = k_{r}q_{r}
\end{cases}
\quad\text{or, that is the same,}\quad
\begin{cases}
k_{1}q_{1}-q_{2} = q_{0},\\
-q_{1} + k_{2}q_{2} - q_{3} = 0,\\
\cdots\\
-q_{r-1}+k_{r}q_{r} = q_{r+1}.
\end{cases}
\end{equation}
One can treat it as the system of linear equations in variables $q_{1},\ldots, q_{r}$ with the matrix
\begin{equation*}
\begin{pmatrix}
k_{1} & -1 & 0 & 0 & \cdots & 0 & 0 \\
-1    & k_{2} & -1 & 0 & \cdots & 0 & 0 \\
0     & -1 & k_{3} & -1 &  \cdots & 0 & 0 \\
\cdots \\
0     & 0  & 0     & 0 & \cdots & -1 & k_{r}
\end{pmatrix}
\end{equation*}
In what follows, we denote its determinant by $\mathbb{K}_{r}(k_{1},\ldots, k_{r})$ or, briefly, by $\mathbb{K}(k_{1},\ldots, k_{r})$.
One can check an obvious properties of such determinants:
\begin{align}
& \mathbb{K}_{1}(k_{1}) = k_{1}, \quad \mathbb{K}_{2}(k_{1},k_{2})=k_{1}k_{2}-1; \label{lab_06} \\[6pt]
& \mathbb{K}_{r}(k_{1},\ldots,k_{r}) = k_{r}\mathbb{K}_{r-1}(k_{1},\ldots,k_{r-1})-\mathbb{K}_{r-2}(k_{1},\ldots,k_{r-2}), \quad r\geqslant 3; \label{lab_07} \\[6pt]
& \mathbb{K}_{r}(k_{1},\ldots,k_{r}+1) =  \mathbb{K}_{r}(k_{1},\ldots,k_{r}) +  \mathbb{K}_{r-1}(k_{1},\ldots,k_{r-1}),\quad r\geqslant 1 \label{lab_08}
\end{align}
(by definition, we set $\mathbb{K}_{-1}=0$, $\mathbb{K}_{0}=1$). The expressions $\mathbb{K}_{r}(k_{1},\ldots, k_{r})$ turns out to be so-called \textit{modified continuants} (see, for example, \cite[p.~133]{Ustinov_2009}). In particular, they satisfy the identities
\begin{align}
& \mathbb{K}_{r}(k_{r},\ldots,k_{1}) =  \mathbb{K}_{r}(k_{1},\ldots,k_{r}); \label{lab_09} \\[6pt]
& \mathbb{K}_{r+s}(k_{1},\ldots,k_{r+s}) = \notag \\[6pt]
& =\mathbb{K}_{r}(k_{1},\ldots,k_{r})\mathbb{K}_{s}(k_{r+1},\ldots,k_{r+s}) - \mathbb{K}_{r-1}(k_{1},\ldots,k_{r-1})\mathbb{K}_{s-1}(k_{r+2},\ldots,k_{r+s}),\label{lab_10}
\end{align}
which are similar to the relations (6.131), (6.133) from \cite[Section 6.7]{Graham_Knuth_Patashnik_1998}. Also, they satisfy the relation
\begin{equation}\label{lab_11}
\mathbb{K}_{r}(k_{1},\ldots,k_{r})\mathbb{K}_{r}(k_{2},\ldots,k_{r+1}) - \mathbb{K}_{r-1}(k_{2},\ldots,k_{r})\mathbb{K}_{r+1}(k_{1},\ldots,k_{r+1}) = 1.
\end{equation}
Using (\ref{lab_06}) and (\ref{lab_07}), one can prove (\ref{lab_11}) by induction with respect to $r$ (see, for example, \cite[\S\S 2.2, 2.3]{Smirnov_2022}).
For brevity, in what follows, if the tuple contains the same components then we write $a^{n}$ instead of $\underbrace{a,a,\ldots,a}_{n}$. In this case, $a$ means both the single term and the sequence of terms: $(1,2,2,2,2) = (1,2^{4})$, $(2,1,4,1,4,3) = (2,(1,4)^{2},3)$ etc. According to this, we write $\mathbb{K}_{5}(1,2^{4})$, $\mathbb{K}_{5}(2,(1,4)^{2},3)$ instead of
$\mathbb{K}_{5}(1,2,2,2,2)$, $\mathbb{K}_{5}(2,1,4,1,4,3)$ and so on.

The formulas (\ref{lab_06}), (\ref{lab_07}) together with the system (\ref{lab_05}) imply that
\begin{equation}\label{lab_12}
q_{i+1} = q_{1}\mathbb{K}_{i}(k_{1},\ldots,k_{i}) - q_{0}\mathbb{K}_{i-1}(k_{2},\ldots,k_{i}).
\end{equation}
Given integers $c_{0}\not\equiv c_{1}\pmod{D}$, we define $\mathcal{A}_{r}^{\circ}(D,c_{0},c_{1})$ as the set of tuples $\mathbf{k} = (k_{1},\ldots,k_{r})$ with integer components $k_{1},\ldots,k_{r}\geqslant 1$ satisfying the conditions
\begin{equation}\label{lab_13}
\begin{cases}
c_{1}\,\mathbb{K}_{i}(k_{1},\ldots,k_{i})-c_{0}\,\mathbb{K}_{i-1}(k_{2},\ldots,k_{i})\not\equiv c_{0}\pmod{D},\quad i = 1,2,\ldots, r-1,\\
c_{1}\,\mathbb{K}_{r}(k_{1},\ldots,k_{r})-c_{0}\,\mathbb{K}_{r-1}(k_{2},\ldots,k_{r})\equiv c_{0}\pmod{D}.
\end{cases}
\end{equation}
If $r=1$ then (\ref{lab_13}) takes a form
\[
c_{1}\,\mathbb{K}_{1}(k_{1})-c_{0}\,\mathbb{K}_{0}\equiv c_{0}\pmod{D},\quad\text{that is,}\quad k_{1}c_{1}\equiv 2c_{0}\pmod{D}.
\]
Let $r\geqslant 2$. Fix the tuple $\mathbf{k}\in \mathcal{A}_{r}^{\circ}(D,c_{0},c_{1})$ and integers $q_{0},q_{1}$ satisfying the conditions
\begin{equation*}
q_{0}\equiv c_{0}\pmod{D}, \quad q_{1}\equiv c_{1}\pmod{D}.
\end{equation*}
Then, determining $q_{2}, q_{3}, \ldots, q_{r+1}$ by (\ref{lab_12}) we find that
\begin{equation}\label{lab_14}
q_{i}\not\equiv c_{0}\pmod{D}, \quad i = 2,3,\ldots,r, \quad q_{r+1}\equiv c_{0}\pmod{D}.
\end{equation}
Conversely, if the fractions (\ref{lab_03}) correspond to the tuple $\mathbf{k}$ and satisfy the conditions (\ref{lab_14}), then the relations (\ref{lab_13}) hold true, or, that is the same, $\mathbf{k}\in \mathcal{A}_{r}^{\circ}(D,c_{0},c_{1})$.

In fact, (\ref{lab_04}) means that the order $Q$ of Farey sequence and the pair of consecutive denominators $q_{0}$, $q_{1}$ generate the chains of denominators $q_{2}, q_{3}, \ldots$  and natural numbers $k_{1}, k_{2}, \ldots$. Now we formulate the condition for the quantities $q_{0}$ и $q_{1}$ that guarantee that the corresponding chain of denominators $q_{0}, q_{1}, q_{2}, \ldots, q_{r+1}$ will satisfy (\ref{lab_14}).
\vspace{0.5cm}

\textsc{Lemma 1.} \emph{Consider the following conditions:} \\[3pt]

(a) \emph{there exists the set of $(r+2)$ consecutive fractions of $\Phi_{Q}$ corresponding to the tuple} $\mathbf{k} = (k_{1},\ldots,k_{r})\in \mathcal{A}_{r}^{\circ}(D,c_{0},c_{1})$;\\[3pt]

(b) \emph{there exists the pair of coprime integers $q_{0}, q_{1}$ satisfying the relations}
\[
1\leqslant q_{0}, q_{1}\leqslant Q,\quad q_{0}+q_{1}>Q,\quad q_{0}\equiv c_{0}\pmod{D},\quad q_{1}\equiv c_{1}\pmod{D}
\]
\emph{and such that the inequalities}
\begin{equation}\label{lab_15}
\frac{Q+q_{0}\mathbb{K}_{i-1}(k_{2},\ldots,k_{i-1},k_{i}+1)}{\mathbb{K}_{i}(k_{1},\ldots,k_{i-1},k_{i}+1)}<q_{1}\leqslant
\frac{Q+q_{0}\mathbb{K}_{i-1}(k_{2},\ldots,k_{i-1},k_{i})}{\mathbb{K}_{i}(k_{1},\ldots,k_{i-1},k_{i})}
\end{equation}
\emph{hold for any} $i$, $1\leqslant i\leqslant r$. \\ %[3pt]

\emph{Then} (a) \emph{and} (b) \emph{are equivalent. Thus, there exists a one-to-one correspondence between the above $(r+2)$-tuples of fractions and the set of integer pairs $(q_{0},q_{1})$ with the aforementioned properties.}

\vspace{0.3cm}

\textsc{Prof.} Necessity. Suppose that the fractions (\ref{lab_03}) correspond to the tuple $\mathbf{k}$. Then $q_{0}, q_{1},\ldots,q_{r},q_{r+1}$ satisfy to the system (\ref{lab_05}) and, therefore, to the relations (\ref{lab_12}). Thus, the condition $q_{i+1}\leqslant Q$ leads to the upper bound for $q_{1}$. Similarly, both the relation (\ref{lab_08}) and (\ref{lab_12}) imply the identity
\begin{equation}\label{lab_16}
q_{i+1}+q_{i} = q_{1}\mathbb{K}_{i}(k_{1},\ldots,k_{i}+1)-q_{0}\mathbb{K}_{i-1}(k_{2},\ldots,k_{i}+1).
\end{equation}
Then the desired lower bound for $q_{1}$ follows from (\ref{lab_16}) and from the inequality $q_{i+1}+q_{i}>Q$.

Sufficiency. Let $q_{0}, q_{1}$ be the integers with the required properties, and let $q_{2},\ldots,q_{r+1}$ satisfy to the system (\ref{lab_05}) (or defined by (\ref{lab_12})).
Then the conditions of the lemma together with the identity (\ref{lab_08}) imply the inequalities
\[
q_{i}\leqslant Q,\quad  q_{i}+q_{i+1}>Q
\]
for any $i$, $1\leqslant i\leqslant r$.

Further, let us find the numbers $a_{0}$ and $a_{1}$ such that the fractions $a_{0}/q_{0}$ and $a_{1}/q_{1}$ are consecutive in $\Phi_{Q}$ and the inequalities $0<a_{0}/q_{0}<1$ hold true.
In view of the condition $q_{0}\equiv 0\pmod{D}$, we have $q_{0} \geqslant D >1$. Then the equation $a_{1}q_{0} - a_{0}q_{1}=1$ imply that $a_{0}\equiv -q_{1}^{*}\pmod{q_{0}}$.
The last congruence together with the inequalities $1\leqslant a_{0}\leqslant q_{0}-1$ determine $a_{0}$ uniquely. So, it remains to define $a_{1} = (a_{0}q_{1}+1)/q_{0}$. Given the pair $a_{0}, a_{1}$ и and the tuple $\mathbf{k}$, consider the system
\begin{equation*}
\begin{cases}
a_{0} + a_{2} = k_{1}a_{1},\\
a_{1} + a_{3} = k_{2}a_{2},\\
\cdots\\
a_{r-1}+a_{r+1} = k_{r}a_{r}
\end{cases}
\end{equation*}
in variables $a_{2}, a_{3},\ldots, a_{r+1}$. Similarly to (\ref{lab_12}), we find:
\begin{equation}\label{lab_17}
a_{i+1} = a_{1}\mathbb{K}_{i}(k_{1},\ldots,k_{i}) - a_{0}\mathbb{K}_{i-1}(k_{2},\ldots,k_{i}).
\end{equation}
Using (\ref{lab_11}), (\ref{lab_12}) and (\ref{lab_16}), we conclude that
\begin{equation}\label{lab_18}
a_{i+1}q_{i} - a_{i}q_{i+1} = 1.
\end{equation}
Now let us consider the fractions $a_{i}/q_{i}$, $0\leqslant i\leqslant r$. In view of (\ref{lab_18}), they are irreducible and $a_{i}/q_{i}<a_{i+1}/q_{i+1}$ for any $i$. Since $q_{i}\leqslant Q$, all these fractions lies in $\Phi_{Q}$ (modulo $1$). Finally, the inequality $q_{i}+q_{i+1}>Q$ means that there are no fractions in $\Phi_{Q}$ between $a_{i}/q_{i}$ and $a_{i+1}/q_{i+1}$. Hence, the above fractions form the  $(r+2)$-tuple of consecutive elements of $\Phi_{Q}$ with the desired properties. Lemma is proved. $\square$
\vspace{0.3cm}

\textsc{Lemma 2.} \emph{Suppose that $r\geqslant 1$, $k_{1},\ldots,k_{r}\geqslant 1$ are integers and the set $\mathcal{T}(k_{1},\ldots,k_{r})$ is non-empty. Then}\\
(a) $\mathbb{K}_{r}(k_{1},\ldots,k_{r})\geqslant 1$;\\
(b) $\mathcal{T}(k_{1},\ldots,k_{r})$ \emph{is defined by the system of the following inequalities}
\[
0<x,y\leqslant 1,\quad x+y>1,\quad f_{i}(x;k_{1},\ldots,k_{i})<y\leqslant  g_{i}(x;k_{1},\ldots,k_{i}),\quad i = 1,2,\ldots, r,
\]
\emph{where}
\[
f_{i}(x;k_{1},\ldots,k_{i}) = \frac{1+x\,\mathbb{K}_{i-1}(k_{2},\ldots,k_{i}+1)}{\mathbb{K}_{i}(k_{1},\ldots,k_{i}+1)},\quad
g_{i}(x;k_{1},\ldots,k_{i}) = \frac{1+x\,\mathbb{K}_{i-1}(k_{2},\ldots,k_{i})}{\mathbb{K}_{i}(k_{1},\ldots,k_{i})}.
\]
\textsc{Proof.} We prove (a), (b) by induction with respect to $r$. In the case $r=1$, we obviously have $\mathbb{K}_{1}(k_{1}) = k_{1}\geqslant 1$ and the set $\mathcal{T}(k_{1})$ is defined by the relation
\[
\biggl[\frac{1+x}{y}\biggr] = k_{1},\quad\text{that is, by the inequalities}\quad k_{1}\leqslant \frac{1+x}{y} < k_{1}+1.
\]
After easy calculations, we find that
\[
f_{1}(x;k_{1}) = \frac{1+x}{k_{1}+1} < y \leqslant \frac{1+x}{k_{1}} =  g_{1}(x;k_{1}).
\]

Suppose now that (a) is true for any $r$, $1\leqslant r\leqslant m-1$, and let $\omega = \mathcal{T}(k_{1},\ldots,k_{m})$ be non-empty set such that $\mathbb{K}_{m}(k_{1},\ldots,k_{m})\leqslant 0$.
Taking any point $(x,y) = (x_{1},y_{1})\in \omega$ we define
\[
(x_{s},y_{s}) = T(x_{s-1},y_{s-1}) = T^{s-1}(x_{1},y_{1}).
\]
Since
\[
(x_{s},y_{s})\in T^{s-1}\omega \subseteq T^{s-1}\circ T^{-(s-1)}\mathcal{T}(k_{s}) = \mathcal{T}(k_{s})
\]
then
\begin{multline}\label{lab_18a}
\biggl[\frac{1+x_{s}}{y_{s}}\biggr] = k_{s}\quad\text{for any}\quad s = 1,2,\ldots,m,\\
(x_{s+1},y_{s+1}) = T(x_{s},y_{s}) = (y_{s},k_{s}y_{s}-x_{s}),
\end{multline}
and therefore
\begin{equation}\label{lab_18b}
\begin{cases}
x_{s+1} = y_{s},\\
y_{s+1} = k_{s}y_{s} - x_{s},
\end{cases}
\quad s = 1,2,\ldots, m-1.
\end{equation}
Using (\ref{lab_18b}) together with (\ref{lab_07}), we easily find that
\begin{equation}\label{lab_18c}
\begin{cases}
y_{s+1} = y_{1}\mathbb{K}_{s}(k_{1},\ldots,k_{s})-x_{1}\mathbb{K}_{s-1}(k_{2},\ldots,k_{s}),\quad s = 1,2,\ldots,m-1,\\
x_{s+1} = y_{1}\mathbb{K}_{s-1}(k_{1},\ldots,k_{s-1})-x_{1}\mathbb{K}_{s-2}(k_{2},\ldots,k_{s-1}),\quad s = 2,3,\ldots,m-1.\\
\end{cases}
\end{equation}
Since $0<x_{1},y_{1}\leqslant 1$, taking $s=m$, we get
\[
y_{m+1} = y_{1}\mathbb{K}_{m}(k_{1},\ldots,k_{m})-x_{1}\mathbb{K}_{m-1}(k_{2},\ldots,k_{m})\leqslant -x_{1}\mathbb{K}_{m-1}(k_{2},\ldots,k_{m}).
\]
Next, since $\omega = \mathcal{T}(k_{1})\bigcap T^{-1}\mathcal{T}(k_{2},\ldots,k_{m})\ne \varnothing$, then the set $\mathcal{T}(k_{2},\ldots,k_{m})$ is also non-empty. By induction, we have
$\mathbb{K}_{m-1}(k_{2},\ldots,k_{m})\geqslant 1$. Therefore, $y_{m+1}\leqslant -x_{1}<0$. This yields to contradiction. Hence, (a) holds true for any $r\geqslant 1$.

Further, (\ref{lab_18a}) yields for $s=m$:
\begin{equation}\label{lab_18d}
k_{m}\leqslant \frac{1+x_{m}}{y_{m}}<k_{m}+1,\quad\text{and hence}\quad
\begin{cases}
k_{m}y_{m}\leqslant 1+x_{m},\\
1+x_{m}<(k_{m}+1)y_{m}.
\end{cases}
\end{equation}
Using (\ref{lab_18c}) we rewrite the first inequality from (\ref{lab_18d}) as follows:
\begin{multline*}
k_{m}\bigl(y_{1}\mathbb{K}_{m-1}(k_{1},\ldots,k_{m-1})-x_{1}\mathbb{K}_{m-2}(k_{2},\ldots,k_{m-1})\bigr)\leqslant \\
\leqslant 1 + y_{1}\mathbb{K}_{m-2}(k_{1},\ldots,k_{m-2})-x_{1}\mathbb{K}_{m-3}(k_{1},\ldots,k_{m-2}),
\end{multline*}
or, that is the same,
\begin{multline*}
y_{1}\bigl(k_{m}\mathbb{K}_{m-1}(k_{1},\ldots,k_{m-1})-\mathbb{K}_{m-2}(k_{1},\ldots,k_{m-2})\bigr) \leqslant \\
\leqslant 1 + x_{1}\bigl(k_{m}\mathbb{K}_{m-2}(k_{2},\ldots,k_{m-1}) - \mathbb{K}_{m-3}(k_{1},\ldots,k_{m-2})\bigr).
\end{multline*}
In view of (\ref{lab_07}), we get
\[
y_{1}\mathbb{K}_{m}(k_{1},\ldots,k_{m}) \leqslant 1 + x_{1}\mathbb{K}_{m-1}(k_{2},\ldots,k_{m}).
\]
By the assertion (a) of this lemma, $\mathbb{K}_{m}(k_{1},\ldots,k_{m})\geqslant 1 >0$. Hence,
\[
y_{1}\leqslant \frac{1 + x_{1}\,\mathbb{K}_{m-1}(k_{2},\ldots,k_{m})}{\mathbb{K}_{m}(k_{1},\ldots,k_{m})} = g_{m}(x_{1};k_{1},\ldots,k_{m}).
\]
Transforming the second inequality from (\ref{lab_18d}) in the same way, we get
\begin{multline*}
1 + y_{1}\mathbb{K}_{m-2}(k_{1},\ldots,k_{m-2})-x_{1}\mathbb{K}_{m-3}(k_{2},\ldots,k_{m-2})\bigr) < \\
< (k_{m}+1)\bigl(y_{1}\mathbb{K}_{m-1}(k_{1},\ldots,k_{m-1}) - x_{1}\mathbb{K}_{m-2}(k_{2},\ldots,k_{m-1})\bigr),\\
1 + x_{1}\bigl((k_{m}+1)\mathbb{K}_{m-2}(k_{2},\ldots,k_{m-1}) - \mathbb{K}_{m-3}(k_{2},\ldots,k_{m-2})\bigr) < \\
< y_{1}\bigl((k_{m}+1)\mathbb{K}_{m-1}(k_{1},\ldots,k_{m-1}) - \mathbb{K}_{m-2}(k_{1},\ldots,k_{m-2})\bigr)
\end{multline*}
and, finally,
\[
1+x_{1}\mathbb{K}_{m-1}(k_{2},\ldots,k_{m-1},k_{m}+1) < y_{1}\mathbb{K}_{m}(k_{1},\ldots,k_{m}+1).
\]
In view of (\ref{lab_08}) and the assertion (a) of this lemma, we have
\[
\mathbb{K}_{m}(k_{1},\ldots,k_{m}+1) = \mathbb{K}_{m}(k_{1},\ldots,k_{m}) + \mathbb{K}_{m-1}(k_{1},\ldots,k_{m-1})\geqslant 1+1>0.
\]
Therefore,
\[
y_{1} > \frac{1+x_{1}\,\mathbb{K}_{m-1}(k_{2},\ldots,k_{m-1},k_{m}+1)}{\mathbb{K}_{m}(k_{1},\ldots,k_{m}+1)} = f_{m}(x_{1};k_{1},\ldots,k_{m}).
\]
Lemma is proved. $\square$ \\

\textsc{Remark 1.} From this lemma, it follows that any set $\mathcal{T}(k_{1},\ldots,k_{r})$ is either empty or convex polygon.

\textsc{Remark 2.} The assertion (a) of the lemma gives only the necessary condition for the non-emptiness of the region $\omega = \mathcal{T}(k_{1},\ldots,k_{r})$. If $\omega$ is empty then
the value of $\mathbb{K}_{r}(k_{1},\ldots,k_{r})$ can be both positive and negative. Thus, one can check that
\[
\mathcal{T}(3,2^{4},1,6),\quad\mathcal{T}(3,2^{5},1,6) = \varnothing,
\]
but
\[
\mathbb{K}_{7}(3,2^{4},1,6) = 1 \quad\text{and}\quad \mathbb{K}_{8}(3,2^{5},1,6) = -1.
\]
For $\omega\ne \varnothing$, the lower bound $\mathbb{K}_{r}(k_{1},\ldots,k_{r})\geqslant 1$ is best possible. For example, $\mathcal{T}(1,2^{r-1})\ne \varnothing$ for any $r\geqslant 2$, and $\mathbb{K}_{r}(1,2^{r-1}) = 1$ (see formulae (\ref{lab_53}) of Lemma 20 below). \\

Lemmas 1 and 2 allow one to reduce the calculation of $N(Q;r,D,c_{0})$ to the counting of lattice points in polygons.
\vspace{0.3cm}

\textsc{Lemma 3.} \emph{Let $r,D,c_{0}, c_{1}$ be integers such that $r\geqslant 1$, $D\geqslant 2$, $0\leqslant c_{0}, c_{1}\leqslant D-1$, $c_{0}\ne c_{1}$}, $\text{GCD}(D,c_{0},c_{1})=1$ \emph{and let $Q\geqslant r$. Suppose also that $(k_{1},\ldots,k_{r})\in \mathcal{A}_{r}^{\circ}(D,c_{0},c_{1})$. Then the number of appearances of the tuple}
(\ref{lab_01}) \emph{of consecutive Farey fractions in $\Phi_{Q}$ satisfying the conditions}
\[
q_{0},q_{r+1}\equiv c_{0}\pmod{D},\quad q_{1}\equiv c_{1}\pmod{D},\quad q_{i}\not\equiv c_{0}\pmod{D},\quad i=2,\ldots, r,
\]
\emph{and corresponding to the tuple $(k_{1},\ldots,k_{r})$ equals to the number of primitive points $(q_{0},q_{1})$ in $Q\mathcal{T}(k_{1},\ldots,k_{r})$ such that}
\[
q_{0}\equiv c_{0}\pmod{D},\quad q_{1}\equiv c_{1}\pmod{D}.
\]
%\vspace{0.3cm}

\textsc{Proof.} By Lemma 1, the desired number of tuples (\ref{lab_01}) coincides with the number of coprime numbers $(q_{0},q_{1})$ satisfying the conditions
\[
0<q_{0},q_{1}\leqslant Q, \quad q_{0}+q_{1}>Q, \quad q_{0}\equiv c_{0}\pmod{D}, \quad q_{1}\equiv c_{1}\pmod{D}
\]
and such that
\begin{equation}\label{lab_22}
\frac{Q+q_{0}\mathbb{K}_{i-1}(k_{2},\ldots,k_{i-1},k_{i}+1)}{\mathbb{K}_{i}(k_{1},\ldots,k_{i-1},k_{i}+1)}<q_{1}\leqslant\frac{Q+q_{0}\mathbb{K}_{i}(k_{2},\ldots,k_{i})}{\mathbb{K}_{i}(k_{1},\ldots,k_{i})}
\end{equation}
for any $i$, $1\leqslant i\leqslant r$. After the change of variables $x = q_{0}/Q$, $y = q_{1}/Q$, the inequalities (\ref{lab_22}) take the form (in the notations of Lemma 2)
\[
f_{i}(x;k_{1},\ldots,k_{i})<y\leqslant g_{i}(x;k_{1},\ldots,k_{i}),\quad i = 1,2,\ldots,r
\]
(obviously, one has $0<x,y\leqslant 1$, $x+y>1$). By Lemma 2, the point $(x,y)$ lies in $\mathcal{T}(k_{1},\ldots,k_{r})$. Hence, $(q_{0},q_{1}) = Q\cdot (x,y)$ is a primitive point of $Q\mathcal{T}(k_{1},\ldots,k_{r})$. Lemma is proved. $\square$
\vspace{0.3cm}

\textsc{Corollary.} \emph{Let $r, D, c_{0}$ be integers such that $r\geqslant 1$, $D\geqslant 2$, $0\leqslant c_{0}\leqslant D-1$, and suppose that $Q\geqslant r$. Then}
\begin{equation}\label{lab_23}
N(Q;r,D,c_{0}) = \sum\limits_{\substack{c_{1} = 1, \; c_{1}\ne c_{0} \\ \text{GCD}(D,c_{0},c_{1})=1}}^{D-1}\sum\limits_{\mathbf{k}\in \mathcal{A}_{r}^{\circ}(D,c_{0},c_{1})}N_{c_{0},c_{1}}\bigl(Q\mathcal{T}(\mathbf{k})\bigr),
\end{equation}
\emph{where the symbol $N_{\sigma_{0},\sigma_{1}}(\Omega)$ stands for the number of primitive lattice points $(q_{0},q_{1})$ lying in the region $\Omega$ and satisfying the conditions}
$q_{0}\equiv \sigma_{0}\pmod{D}$, $q_{1}\equiv \sigma_{1}\pmod{D}$.
\vspace{0.3cm}

The map $T$ is discontinuous but it preserves measure (see \cite[Lemma 3]{Boca_Cobeli_Zaharescu_2001}). Obviously, the transformation $(x,y)\mapsto Q\cdot T(x,y)$ maps the lattice points onto the lattice points. It appears also that such transformation and its iterations $QT^{2}$, $QT^{3}, \ldots$ establish a one-to-one correspondence between the primitive points of the initial region and the primitive points of its images (see Lemmas 4 and 5 below).

Using this fact, we act to the region $Q\cdot \mathcal{T}(\mathbf{k})$ by some iteration of the map $T$ to move it to some special region in order to simplify the calculation of the number of its primitive points (see Lemma 7).

\textsc{Lemma 4.} \emph{Let $r\geqslant 2$, $k_{1},\ldots,k_{r}\geqslant 1$ be integers such that the region $\omega = \mathcal{T}(k_{1},\ldots,k_{r})$ in non-empty. Further, let $(m,n)$ be a lattice point of $Q\cdot\omega$. Then, for any $s$, $1\leqslant s\leqslant r-1$, the point}
\[
(u,v) = QT^{s}\biggl(\frac{m}{Q},\frac{n}{Q}\biggr)
\]
\emph{is a lattice point of the region $QT^{s}(\omega)$ and}
\begin{multline}\label{lab_24}
u = n\,\mathbb{K}_{s-1}(k_{1},\ldots,k_{s-1}) - m\,\mathbb{K}_{s-2}(k_{2},\ldots,k_{s-1}),\\
v = n\,\mathbb{K}_{s}(k_{1},\ldots,k_{s}) - m\,\mathbb{K}_{s-1}(k_{2},\ldots,k_{s}).
\end{multline}
\emph{Moreover, if $(m,n)$ is primitive then $(u,v)$ is also primitive.}
\vspace{0.5cm}

\textsc{Proof.} Induction on $s$. Since
\[
\biggl(\frac{m}{Q},\frac{n}{Q}\biggr) \in \omega \subseteq \mathcal{T}(k_{1})\quad\text{then}\quad
T\biggl(\frac{m}{Q},\frac{n}{Q}\biggr) = \biggl(\frac{n}{Q},k_{1}\frac{n}{Q} - \frac{m}{Q}\biggr) = \biggl(\frac{u}{Q},\frac{v}{Q}\biggr)
\]
and therefore
\[
u = n = n\mathbb{K}_{0}-m\mathbb{K}_{-1},\quad v = k_{1}n-m = n\mathbb{K}_{1}(k_{1})-m\mathbb{K}_{0}.
\]
In particular, $\text{GCD}(u,v) = \text{GCD}(n,k_{1}n-m) = \text{GCD}(n,m)$. So, the assertion is true for $s=1$.

It remains to check its truth after passing from $s$ to $s+1$ (we assume here that $r\geqslant 3$ and $1\leqslant s\leqslant r-2$). By induction, $(u_{s},v_{s}) = QT^{s}(m/Q,n/Q)$ is a lattice point of the region $QT^{s}(\omega)$, and its coordinates satisfy (\ref{lab_24}). Moreover, if $\text{GCD}(m,n)=1$ then $\text{GCD}(u_{s},v_{s})=1$.

Since
\begin{multline*}
T^{s}(\omega) = T^{s}\Bigl(\mathcal{T}(k_{1})\bigcap T^{-1}\mathcal{T}(k_{2})\bigcap \ldots \bigcap T^{-s}\mathcal{T}(k_{s+1})\bigcap \ldots \bigcap T^{-(r-1)}\mathcal{T}(k_{r})\Bigr) = \\
= T^{s}\mathcal{T}(k_{1})\bigcap T^{s-1}\mathcal{T}(k_{2})\bigcap \ldots \bigcap \mathcal{T}(k_{s+1})\bigcap \ldots \bigcap T^{-(r-s-1)}\mathcal{T}(k_{r})\subseteq \mathcal{T}(k_{s+1}),
\end{multline*}
then the point $\bigl(u_{s}/Q,v_{s}/Q\bigr)$ lies in $\mathcal{T}(k_{s+1})$. Hence,
\[
T^{s+1}\biggl(\frac{m}{Q},\frac{n}{Q} \biggr) = \biggl(\frac{u_{s+1}}{Q},\frac{v_{s+1}}{Q} \biggr) = T\biggl(\frac{u_{s}}{Q},\frac{v_{s}}{Q}\biggr) =
\biggl(\frac{v_{s}}{Q},k_{s+1}\frac{v_{s}}{Q}-\frac{u_{s}}{Q}\biggr),
\]
so one has $u_{s+1} = v_{s}$, $v_{s+1} = k_{s+1}v_{s}-u_{s}$. After the replacement $u = u_{s}$ and $v = v_{s}$ by the expressions (\ref{lab_24}) we get
\begin{multline*}
u_{s+1} = n\,\mathbb{K}_{s}(k_{1},\ldots,k_{s}) - m\,\mathbb{K}_{s-1}(k_{2},\ldots,k_{s}) = \\
= n\,\mathbb{K}_{(s+1)-1}(k_{1},\ldots,k_{(s+1)-1}) - m\,\mathbb{K}_{(s+1)-2}(k_{2},\ldots,k_{(s+1)-1}),\\
v_{s+1} = k_{s+1}\bigl(n\,\mathbb{K}_{s}(k_{1},\ldots,k_{s}) - m\,\mathbb{K}_{s-1}(k_{2},\ldots,k_{s})\bigr)-\\
-\bigl(n\,\mathbb{K}_{s-1}(k_{1},\ldots,k_{s-1}) - m\,\mathbb{K}_{s-2}(k_{2},\ldots,k_{s-1})\bigr) =\\
= n\,\mathbb{K}_{s+1}(k_{1},\ldots,k_{s+1}) - m\,\mathbb{K}_{s}(k_{2},\ldots,k_{s+1}).
\end{multline*}
To finish the proof, it remains to note that
\[
\text{GCD}(u_{s+1},v_{s+1}) = \text{GCD}(v_{s},k_{s+1}v_{s}-u_{s}) = \text{GCD}(u_{s},v_{s}) = \text{GCD}(m,n).\quad \square
\]

\textsc{Lemma 5.} \emph{Let $r\geqslant 2$, $k_{1},\ldots,k_{r}\geqslant 1$ be integers such that the region $\omega = \mathcal{T}(k_{1},\ldots,k_{r})$ in non-empty. Further, let $1\leqslant s\leqslant r-1$, $\Omega = T^{s}(\omega)$ and let $(u,v)$ be a point of $Q\cdot\Omega$. Then}
\[
(m,n) = QT^{-s}\biggl(\frac{u}{Q},\frac{v}{Q}\biggr)
\]
\emph{is a point of $Q\cdot \omega$ and}
\begin{multline}\label{lab_25}
m = u\,\mathbb{K}_{s}(k_{s},\ldots,k_{1}) - v\,\mathbb{K}_{s-1}(k_{s-1},\ldots,k_{1}),\\
n = u\,\mathbb{K}_{s-1}(k_{s},\ldots,k_{2}) - v\,\mathbb{K}_{s-2}(k_{s-1},\ldots,k_{2}).
\end{multline}
\emph{Moreover, if $(u,v)$ is an integer lattice point then $(m,n)$ is also integer lattice point.}
\vspace{0.5cm}

\textsc{Proof.} Prove that, for any $j$, $1\leqslant j\leqslant s$, the numbers $m_{j}, n_{j}$ defined by the relation
\[
\biggl(\frac{m_{j}}{Q},\frac{n_{j}}{Q} \biggr) = T^{-j}\biggl(\frac{u}{Q},\frac{v}{Q} \biggr),
\]
satisfy to the equations
\begin{multline}\label{lab_26}
m_{j} = u\,\mathbb{K}_{j}(k_{s},\ldots,k_{s-j+1}) - v\,\mathbb{K}_{j-1}(k_{s-1},\ldots,k_{s-j+1}),\\
n_{j} = u\,\mathbb{K}_{j-1}(k_{s},\ldots,k_{s-j+2}) - v\,\mathbb{K}_{j-2}(k_{s-1},\ldots,k_{s-j+2}).
\end{multline}
Since
\begin{multline*}
\biggl(\frac{m_{1}}{Q},\frac{n_{1}}{Q} \biggr) = T^{-1}\biggl(\frac{u}{Q},\frac{v}{Q} \biggr)\in T^{-1}\circ T^{s}(\omega) = T^{s-1}(\omega) = \\
= T^{s-1}\Bigl(\mathcal{T}(k_{1})\bigcap \ldots \bigcap T^{-(s-1)}\mathcal{T}(k_{s})\bigcap \ldots \bigcap T^{-(r-1)}\mathcal{T}(k_{r})\Bigr)\subseteq \mathcal{T}(k_{s}),
\end{multline*}
then
\[
\biggl(\frac{u}{Q},\frac{v}{Q} \biggr) = T\biggl(\frac{m_{1}}{Q},\frac{n_{1}}{Q} \biggr) = \biggl(\frac{n_{1}}{Q},k_{s}\frac{n_{1}}{Q}-\frac{m_{1}}{Q}\biggr),\quad u=n_{1},\quad v=k_{s}n_{1}-m_{1}
\]
and hence
\begin{equation}\label{lab_27}
m_{1} = k_{s}u-v,\quad n_{1} = v.
\end{equation}
So, the formulas (\ref{lab_26}) hold true for $j=1$.

Suppose now that $s\geqslant 2$, $1\leqslant j\leqslant s-1$ and check that (\ref{lab_26}) holds when passing from $j$ to $j+1$. Indeed, we have
\begin{equation}\label{lab_28}
T\biggl(\frac{m_{j+1}}{Q},\frac{n_{j+1}}{Q}\biggr) = T\circ T^{-(j+1)}\biggl(\frac{u}{Q},\frac{v}{Q} \biggr) = T^{-j}\biggl(\frac{u}{Q},\frac{v}{Q} \biggr) =
\biggl(\frac{m_{j}}{Q},\frac{n_{j}}{Q}\biggr).
\end{equation}
Moreover, since the point $\bigl(u/Q,v/Q\bigr)$ lies in $T^{s}(\omega)$ then  $\bigl(m_{j+1}/Q,n_{j+1}/Q\bigr)$ lies in $T^{-(j+1)}\circ T^{s}(\omega) = T^{s-j-1}(\omega)$.
Obviously, one has
\[
T^{s-j-1}(\omega) = T^{s-j-1}\Bigl(T(k_{1})\bigcap \ldots\bigcap T^{-(s-j-1)\mathcal{T}(k_{s-j})}\bigcap \ldots \bigcap T^{-(r-1)}\mathcal{T}(k_{r})\Bigr)\subseteq \mathcal{T}(k_{s-j}),
\]
so $\bigl(m_{j+1}/Q,n_{j+1}/Q\bigr)$ is a point of $\mathcal{T}(k_{s-j})$. Thus (\ref{lab_28}) implies that
\[
\biggl(\frac{m_{j}}{Q},\frac{n_{j}}{Q}\biggr) = T\biggl(\frac{m_{j+1}}{Q},\frac{n_{j+1}}{Q}\biggr) =\biggl(\frac{n_{j+1}}{Q},k_{s-j}\frac{n_{j+1}}{Q}-\frac{m_{j+1}}{Q}\biggr)
\]
and therefore
\begin{equation}\label{lab_29}
m_{j+1} = k_{s-j}m_{j}-n_{j},\quad n_{j+1} = m_{j}.
\end{equation}
By induction, we replace $m_{j}$ and $n_{j}$ in (\ref{lab_29}) by the expressions (\ref{lab_26}). Thus we find
\begin{multline*}
m_{j+1} = k_{s-j}\bigl(u\,\mathbb{K}_{j}(k_{s},\ldots,k_{s-j+1}) - v\,\mathbb{K}_{j-1}(k_{s-1},\ldots,k_{s-j+1})\bigr) - \\
-\bigl(u\,\mathbb{K}_{j-1}(k_{s},\ldots,k_{s-j+2}) - v\,\mathbb{K}_{j-2}(k_{s-1},\ldots,k_{s-j+2})\bigr) =\\
= u\bigl(k_{s-j}\,\mathbb{K}_{j}(k_{s},\ldots,k_{s-j+1})-\mathbb{K}_{j-1}(k_{s},\ldots,k_{s-j+2})\bigr) - \\
- v\bigl(k_{s-j}\,\mathbb{K}_{j-1}(k_{s-1},\ldots,k_{s-j+1})-\mathbb{K}_{j-2}(k_{s-1},\ldots,k_{s-j+2})\bigr) =\\
= u\,\mathbb{K}_{j+1}(k_{s},\ldots,k_{s-j}) - v\,\mathbb{K}_{j}(k_{s-1},\ldots,k_{s-j})
\end{multline*}
and, similarly,
\[
n_{j+1} = u\,\mathbb{K}_{j}(k_{s},\ldots,k_{s-j}) - v\,\mathbb{K}_{j}(k_{s-1},\ldots,k_{s-j+1}).
\]
Thus, the formulas (\ref{lab_26}) hold true for any $j$. Taking $j=s$, we get (\ref{lab_25}).

Now it remains to note that $m_{j+1}, n_{j+1}$ are integers if and only if $m_{j}, n_{j}$ are integers. Moreover, (\ref{lab_29}) implies that $\text{GCD}(m_{j+1},n_{j+1}) = \text{GCD}(m_{j},n_{j})$ (for integers $m_{j}, n_{j}$), and (\ref{lab_27}) implies that $\text{GCD}(m_{j},n_{j}) = \text{GCD}(u,v)$. Lemma is proved. $\square$
\vspace{0.5cm}

\textsc{Lemma 6.} \emph{Let $c_{0}, c_{1}, D$ be integers such that} $D\geqslant 2$, $0\leqslant c_{0},c_{1}\leqslant D-1$, $\text{GCD}(D,c_{0},c_{1})=1$. \emph{Further, let $r\geqslant 2$, $k_{1},\ldots,k_{r}\geqslant 1$ be integers such that $\omega = \mathcal{T}(k_{1},\ldots,k_{r})$ is non-empty. Finally, let $j\geqslant 1$ and suppose that $0\leqslant \sigma_{0}, \sigma_{1}\leqslant D-1$ are defined by}
\begin{equation}\label{lab_30}
\begin{cases}
\sigma_{0}\equiv c_{1}\,\mathbb{K}_{j-1}(k_{1},\ldots,k_{j-1}) - c_{0}\,\mathbb{K}_{j-2}(k_{2},\ldots,k_{j-1})\pmod{D},\\
\sigma_{1}\equiv c_{1}\,\mathbb{K}_{j}(k_{1},\ldots,k_{j}) - c_{0}\,\mathbb{K}_{j-1}(k_{2},\ldots,k_{j})\pmod{D}.
\end{cases}
\end{equation}
\emph{Then }
\[
N_{c_{0},c_{1}}(Q\cdot \omega) = N_{\sigma_{0},\sigma_{1}}\bigl(Q\cdot T^{j}(\omega)\bigr).
\]
\textsc{Proof.} Let $(q_{0},q_{1})$ be a primitive point of $Q\cdot \omega$ such that
\begin{equation}\label{lab_31}
q_{0}\equiv c_{0}\pmod{D},\quad q_{1}\equiv c_{1}\pmod{D}.
\end{equation}
By Lemma 4, $(u,v) = QT^{j}(q_{0}/Q,q_{1}/Q)$ is a primitive point of $Q\cdot T^{j}(\omega)$ and
\begin{multline}\label{lab_32}
u = q_{0}\,\mathbb{K}_{j-1}(k_{1},\ldots,k_{j-1})-q_{1}\,\mathbb{K}_{j-2}(k_{2},\ldots,k_{j-1})\bigr),\\
v = q_{0}\,\mathbb{K}_{j}(k_{1},\ldots,k_{j})-q_{1}\,\mathbb{K}_{j-1}(k_{2},\ldots,k_{j})\bigr).
\end{multline}
From (\ref{lab_31}), it follows that
\begin{equation}\label{lab_33}
u\equiv \sigma_{0}\pmod{D},\quad v\equiv \sigma_{1}\pmod{D},
\end{equation}
so we have $N_{c_{0},c_{1}}(Q\cdot \omega) \leqslant N_{\sigma_{0},\sigma_{1}}\bigl(Q\cdot T^{j}(\omega)\bigr)$.

Conversely, let $(u,v)$ be a primitive point of $Q\cdot T^{j}(\omega)$ and its coordinates satisfy to (\ref{lab_33}). By Lemma 5,
\[
(q_{0},q_{1}) = Q\cdot T^{-j}\biggl(\frac{u}{Q},\frac{v}{Q}\biggr)
\]
is a primitive point of $Q\cdot\omega$ and
\begin{multline}\label{lab_34}
q_{0} = u\,\mathbb{K}_{j}(k_{j},\ldots,k_{1})-v\,\mathbb{K}_{j-1}(k_{j-1},\ldots,k_{1}),\\
q_{1} = u\,\mathbb{K}_{j-1}(k_{j},\ldots,k_{2})-v\,\mathbb{K}_{j-2}(k_{j-1},\ldots,k_{2}).
\end{multline}
Let us check that $q_{0}$ and $q_{1}$ satisfy (\ref{lab_31}). Indeed, passing to the congruences modulo $D$ and using (\ref{lab_30}) we get
\begin{multline*}
q_{0} \equiv \sigma_{0}\,\mathbb{K}_{j}(k_{j},\ldots,k_{1})-\sigma_{1}\,\mathbb{K}_{j-1}(k_{j-1},\ldots,k_{1}) \pmod{D}\equiv\\
\equiv \bigl\{c_{1}\,\mathbb{K}_{j-1}(k_{1},\ldots,k_{j-1})-c_{0}\,\mathbb{K}_{j-2}(k_{2},\ldots,k_{j-1})\bigr\}\mathbb{K}_{j}(k_{j},\ldots,k_{1}) -\\[6pt]
- \bigl\{c_{1}\,\mathbb{K}_{j}(k_{1},\ldots,k_{j})-c_{0}\,\mathbb{K}_{j-1}(k_{2},\ldots,k_{j})\bigr\}\mathbb{K}_{j-1}(k_{j-1},\ldots,k_{1})\pmod{D}\equiv \\[6pt]
\equiv c_{1}\bigl\{\mathbb{K}_{j}(k_{j},\ldots,k_{1})\mathbb{K}_{j-1}(k_{1},\ldots,k_{j-1}) -\mathbb{K}_{j}(k_{1},\ldots,k_{j})\mathbb{K}_{j-1}(k_{j-1},\ldots,k_{1})\bigr\}\,+ \\[6pt]
+ c_{0}\bigl\{\mathbb{K}_{j-1}(k_{j-1},\ldots,k_{1})\mathbb{K}_{j-1}(k_{2},\ldots,k_{j}) -\mathbb{K}_{j}(k_{j},\ldots,k_{1})\mathbb{K}_{j-2}(k_{2},\ldots,k_{j-1})\bigr\}\pmod{D}.
\end{multline*}
Now both (\ref{lab_09}) and (\ref{lab_11}) imply: $q_{0}\equiv c_{0}\pmod{D}$. The congruence $q_{1}\equiv c_{1}\pmod{D}$ is checked in a similar way.
Therefore, $N_{\sigma_{0},\sigma_{1}}\bigl(Q\cdot T^{j}(\omega)\bigr)\leqslant N_{c_{0},c_{1}}(Q\cdot \omega)$.

Lemma is proved. $\square$

\section{Some properties of the polygons $\boldsymbol{\mathcal{T}(k_{1},\ldots,k_{r})}$}

For the following, we need some properties of the convex polygons $\mathcal{T}(k_{1},\ldots,k_{r})$. We organize these properties as the series of lemmas.
\vspace{0.5cm}

\textsc{Lemma 7.} \emph{If $k=1$ then the closure of $\mathcal{T}(k)$ is a triangle with vertices $\bigl(\tfrac{1}{3},\tfrac{2}{3}\bigr)$, $(0,1)$, $(1,0)$ and area $\tfrac{1}{6}$. If $k\geqslant 2$ then the closure\footnote{In what follows, when speaking about polygons $\mathcal{T}(\mathbf{k})$, we omit the word <<closure>>. Since the closure of $\mathcal{T}(\mathbf{k})$ differs from the initial region by some of its edges only, this will not lead us to misunderstanding.} of $\mathcal{T}(k)$ is a quadrangle with vertices}
\begin{equation*}
A=\biggl(\frac{k}{k+2},\frac{2}{k+2}\biggr),\quad B=\biggl(1,\frac{2}{k}\biggr),\quad C = \biggl(1,\frac{2}{k+1}\biggr),\quad D=\biggl(\frac{k-1}{k+1},\frac{2}{k+1}\biggr)
\end{equation*}
\emph{and area}
\begin{equation*}
\frac{4}{k(k+1)(k+2)}.
\end{equation*}

In what follows, by $|\Omega|$ we denote the area of the measurable compact region $\Omega \subset \mathbb{R}^{2}$. So, Lemma 7 claims that
\begin{equation}\label{lab_35}
\setstretch{1.1}
|\mathcal{T}(k)| =
\begin{cases}
1/6, & \text{for } k=1, \\
\displaystyle \frac{4}{k(k+1)(k+2)}, & \text{for } k\geqslant 2.
\end{cases}
\end{equation}
%\vspace{0.5cm}

\textsc{Lemma 8.}  \emph{For any $k\geqslant 1$, the regions $\mathcal{T}(k)$ and $T(\mathcal{T}(k))$ are symmetrical with respect to the line} $y = x$.
\vspace{0.5cm}

These assertions are well-known (see, for example, \cite{Boca_Cobeli_Zaharescu_2003}, \cite{Boca_Cobeli_Zaharescu_2001}). Lemma 7 follows the definition of the region $\mathcal{T}(k)$, and Lemma 8 can be proved by the direct calculation (see fig. 1). Denoting by $S$ the symmetry with respect to the line $y=x$ we get $T(\mathcal{T}(k)) = S(\mathcal{T}(k))$, but $T\ne S$ (as two transformations of Farey triangle).

\begin{center}
\includegraphics{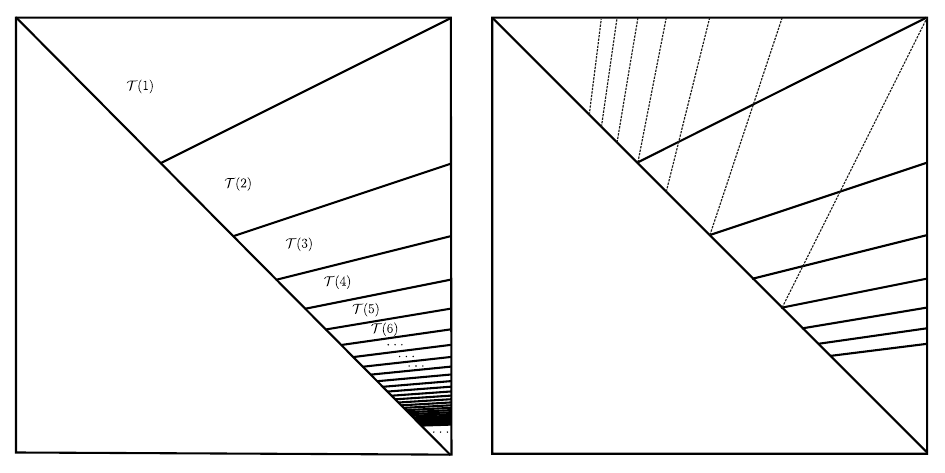}

\fontsize{10}{12pt}\selectfont
\emph{Fig.~1.} The regions $\mathcal{T}(k)$, $k = 1,2,3,\ldots$ (to the left); \\ relative position of the regions $\mathcal{T}(k)$ и $T(\mathcal{T}(k)), 1\leqslant k\leqslant 7$ (to the right).
\fontsize{12}{15pt}\selectfont
\end{center}

\textsc{Lemma 9.} \emph{One has $T^{\,-1}(\mathcal{T}(k))\subset \mathcal{T}(1)$, $T(\mathcal{T}(k))\subset \mathcal{T}(1)$ for any} $k\geqslant 5$.
\vspace{0.5cm}

\textsc{Lemma 10.} \emph{The region $\mathcal{T}(k,m)$ is empty in  the following cases:}\\ [6pt]
\begin{tabular}{l p{0.5\linewidth}}
(a) $m=1$, $k=1$;  & (b) $m=2$, $k\geqslant 5$;  \\ [6pt]
(c) $m=3,4$, $k\geqslant 3$; & (d) $m\geqslant 5$, $k\geqslant 2$.
\end{tabular}

\vspace{0.5cm}

Lemma 9, 10 easily follow from Lemma 8 (see fig. 1).

\vspace{0.5cm}

\textsc{Lemma 11.} \emph{Suppose that $\mathcal{T}(k_{1},\ldots,k_{r})$ is non-empty. Then the linear transformation defined by}
\begin{equation}\label{lab_36}
\begin{cases}
X = y\,\mathbb{K}_{r}(k_{1},\ldots,k_{r})-x\,\mathbb{K}_{r-1}(k_{2},\ldots,k_{r}),\\
Y = y\,\mathbb{K}_{r-1}(k_{1},\ldots,k_{r-1})-x\,\mathbb{K}_{r-2}(k_{2},\ldots,k_{r-1})
\end{cases}
\end{equation}
\emph{maps $\mathcal{T}(k_{1},\ldots,k_{r})$ onto $\mathcal{T}(k_{r},\ldots,k_{1})$.}

\textsc{Proof.} Let $(x,y) = (x_{0},x_{1})$ be the point of $\omega = \mathcal{T}(k_{1},\ldots,k_{r})$. Define the quantities $x_{s}$ as follows: $(x_{s},x_{s+1}) = T^{s}(x_{0},x_{1})$, $s = 1,2,\ldots, r$. Then, noting that $(x_{s},x_{s+1})$ lies in $\mathcal{T}(k_{s+1})$ for $s = 0,1,\ldots, r-1$, we have
\[
(x_{s},x_{s+1}) = T\circ T^{s-1}(x_{0},x_{1}) = T(x_{s-1},x_{s}) = (x_{s},k_{s}x_{s}-x_{s-1}),
\]
and hence
\begin{equation}\label{lab_37}
x_{s+1} = k_{s}x_{s}-x_{s-1},\quad s = 1,2,\ldots, r.
\end{equation}
Further, since the point $(x_{s-1},x_{s})$ belongs to $\mathcal{T}(k_{s})$, we have
\[
k_{s}\leqslant \frac{1+x_{s-1}}{x_{s}}< k_{s}+1,\quad k_{s}x_{s}-x_{s-1} \leqslant 1,\quad k_{s}x_{s}-x_{s-1} > 1-x_{s}.
\]
In view of (\ref{lab_37}), we write it as follows:
\[
1-x_{s} < x_{s+1}\leqslant 1,\quad s = 0,1,\ldots, r.
\]
Thus, one can define the region $\omega$ as the set of all points $(x_{0},x_{1})$ with the condition $0<x_{0}\leqslant 1$, such that the quantities $x_{s}$ defined by (\ref{lab_37}) satisfy the system of inequalities
\begin{equation}\label{lab_38}
\begin{cases}
1-x_{0}<x_{1}\leqslant 1,\\
1-x_{1}<x_{2}\leqslant 1,\\
1-x_{2}<x_{3}\leqslant 1,\\
\cdots \\
1-x_{r}<x_{r+1}\leqslant 1.
\end{cases}
\end{equation}
From (\ref{lab_38}), we conclude that
\[
x_{s+1} = x_{1}\,\mathbb{K}_{s}(k_{1},\ldots,k_{s})-x_{0}\,\mathbb{K}_{s-1}(k_{2},\ldots,k_{s}),\quad s = 1,\ldots, r.
\]
Now we make the change of variables
\begin{equation*}
\begin{cases}
y_{0} = x_{r+1} = x_{1}\,\mathbb{K}_{r}(k_{1},\ldots,k_{r})-x_{0}\,\mathbb{K}_{r-1}(k_{2},\ldots,k_{r}),\\
y_{1} = x_{r} = x_{1}\,\mathbb{K}_{r-1}(k_{1},\ldots,k_{r-1})-x_{0}\,\mathbb{K}_{r-2}(k_{2},\ldots,k_{r-1}).
\end{cases}
\end{equation*}
In view of (\ref{lab_10}), its Jacobian equals to $(-1)$. Hence, the corresponding linear map preserves area.

Now we introduce the auxilliary notations
\[
y_{2} = x_{r-1},\quad y_{3} = x_{r-2},\quad \ldots, \quad y_{r} = x_{1},\quad y_{r+1} = x_{0}.
\]
Then the system (\ref{lab_38}) takes the form
\begin{equation}\label{lab_39}
\begin{cases}
1-y_{0}<y_{1}\leqslant 1,\\
1-y_{1}<y_{2}\leqslant 1,\\
1-y_{2}<y_{3}\leqslant 1,\\
\cdots \\
1-y_{r}<y_{r+1}\leqslant 1,
\end{cases}
\end{equation}
and the relations (\ref{lab_37}) take the form
\begin{equation}\label{lab_40}
y_{s+1} = k_{r-s+1}y_{s}-y_{s-1},\quad s = 1,2,\ldots, r.
\end{equation}
The above arguments mean that both the relations (\ref{lab_39}), (\ref{lab_40}) and the obvious inequalities $0<y_{0}\leqslant 1$ define the region $\mathcal{T}(k_{r},\ldots,k_{1})$.
\vspace{0.5cm}

\textsc{Corollary 1.} \emph{For any integers $k_{1},\ldots,k_{r}$, one has $|\mathcal{T}(k_{r},\ldots,k_{1})| = |\mathcal{T}(k_{1},\ldots,k_{r})|$.}
\vspace{0.3cm}

\textsc{Corollary 2.} \emph{Suppose that $\mathcal{T}(k_{1},\ldots,k_{r})$ is non-empty. Then the linear transformation defined by}
\begin{equation}\label{lab_41}
\begin{cases}
x = Y\,\mathbb{K}_{r}(k_{1},\ldots,k_{r})-X\,\mathbb{K}_{r-1}(k_{1},\ldots,k_{r-1}),\\
y = Y\,\mathbb{K}_{r-1}(k_{2},\ldots,k_{r})-X\,\mathbb{K}_{r-2}(k_{2},\ldots,k_{r-1})
\end{cases}
\end{equation}
\emph{maps $\mathcal{T}(k_{r},\ldots,k_{1})$ onto $\mathcal{T}(k_{1},\ldots,k_{r})$.}
\vspace{0.3cm}

For the proof of Corollary 1, it is sufficient to note that Jacobian of the above transformation is equal to
\[
J = \mathbb{K}_{r}(k_{1},\ldots,k_{r})\mathbb{K}_{r-2}(k_{2},\ldots,k_{r-1}) - \mathbb{K}_{r-1}(k_{1},\ldots,k_{r-1})\mathbb{K}_{r-1}(k_{2},\ldots,k_{r}) = -1.
\]
The Corollary 2 follows directly from (\ref{lab_36}) and the equality $J=-1$.
\vspace{0.5cm}

\textsc{Lemma 12.} \emph{The following relations hold true:} \\[6pt]
\begin{tabular}{l p{0.5\linewidth}}
(a) $|\mathcal{T}(k,1)| = |\mathcal{T}(k)|$ \emph{for any} $k\geqslant 5$;  & (b) $|\mathcal{T}(k,1,2)| = |\mathcal{T}(k)|$ \emph{for any} $k\geqslant 9$; \\ [6pt]
(c) $|\mathcal{T}(2,1,k)| = |\mathcal{T}(k)|$ \emph{for any} $k\geqslant 9$; & (d) $|\mathcal{T}(2,1,k,1)| = |\mathcal{T}(k)|$ \emph{for any} $k\geqslant 9$; \\[6pt]
(e) $|\mathcal{T}(2,1,k,1,2)| = |\mathcal{T}(k)|$ \emph{for any} $k\geqslant 9$; & (f) $|\mathcal{T}(1,k,1)| = |\mathcal{T}(k)|$ \emph{for any} $k\geqslant 6$.\\[6pt]
\end{tabular}

\textsc{Proof.} (a) By definition, $\mathcal{T}(k)$ is defined by the system
\begin{equation}\label{lab_42}
0<x\leqslant 1,\quad \frac{1+x}{k+1}<y\leqslant \frac{1+x}{k}.
\end{equation}
Passing to $\mathcal{T}(k,1)$, we should add to the system (\ref{lab_42}) an extra condition $f(x)<y\leqslant g(x)$, where
\[
f(x) =f(x;k,1) = \frac{1+2x}{2k-1},\quad g(x) = g(x;k,1) = \frac{1+x}{k-1}.
\]
Obviously,
\begin{equation}\label{lab_43}
f(-1/2) = 0,\quad \quad g(-1) = 0.
\end{equation}
Using the notation of Lemma 7 and denoting by $x_{A}, y_{A}, x_{B}, y_{B}$ etc. the coordinates of vertices of $\mathcal{T}(k)$, we see that
\begin{align*}
& f(1) = \frac{3}{2k-1}\leqslant \frac{2}{k+1} = y_{C}\quad\text{for}\quad k\geqslant 5, \\
& g(1) = \frac{2}{k-1} > \frac{2}{k} = y_{B}\quad\text{for}\quad k\geqslant 2.
\end{align*}
The last inequalities together with (\ref{lab_43}) mean that $\mathcal{T}(k)$ lies over the line $y = f(x)$ and lies under the line $y=g(x)$. Therefore, one has $\mathcal{T}(k,1) = \mathcal{T}(k)$ and in particular $|\mathcal{T}(k,1)| = |\mathcal{T}(k)|$.

(b) In fact, if $k\geqslant 5$ then we prove in (a) that $\mathcal{T}(k,1)$ is defined by the same system as $\mathcal{T}(k)$, that is, by (\ref{lab_42}). Passing to $\mathcal{T}(k,1,2)$, we should add to (\ref{lab_42}) the conditions $f(x)<y\leqslant g(x)$, where
\[
f(x) =f(x;k,1,2) = \frac{1+2x}{2k-3},\quad g(x) = g(x;k,1,2) = \frac{1+x}{k-2}.
\]
Obviously, the relations (\ref{lab_43}) still hold. Since
\[
f(1) = \frac{3}{2k-3} = y_{C}\quad \text{for} \quad k\geqslant 9,\quad g(1) = \frac{2}{k-2} > \frac{2}{k} = y_{B}\quad\text{for}\quad k\geqslant 3
\]
then $\mathcal{T}(k,1)$ lies over the line $y = f(x)$ and lies under the line $y = g(x)$. Hence, $\mathcal{T}(k,1,2) = \mathcal{T}(k)$ for all $k\geqslant 9$.

Assertion (c) easily follows from (b) and from the Corollary 1 of Lemma 11.

(d) By Lemma 10, the transformation
\begin{equation*}
\begin{cases}
x = Y\,\mathbb{K}_{3}(k,1,2)-X\,\mathbb{K}_{2}(k,1) = (k-2)Y-(k-1)X,\\[6pt]
y = Y\,\mathbb{K}_{2}(1,2)-X\,\mathbb{K}_{1}(1) = Y-X
\end{cases}
\end{equation*}
maps the region $\mathcal{T}(2,1,k)$ to $\mathcal{T}(k,1,2)$. Correspondingly, the lines
\[
y=1-x,\quad y = \frac{1+x}{k},\quad x=1,\quad y = \frac{1+x}{k+1}
\]
forming the boundary of $\mathcal{T}(k,1,2)=\mathcal{T}(k)$ have the following images
\[
Y = \frac{1+kX}{k-1},\quad Y = \frac{1+X}{2},\quad Y = \frac{1+(k-1)X}{k-2},\quad Y = \frac{1+2X}{3}.
\]
These lines form the boundary of $\mathcal{T}(2,1,k)$. Using this, one can check that $\mathcal{T}(2,1,k)$ is defined by the system
\begin{equation}\label{lab_44}
\setstretch{1.5}
\begin{cases}
\displaystyle \frac{k-5}{k+1}<X\leqslant \frac{k-4}{k+2},\quad \frac{1+2X}{3} < Y\leqslant \frac{1+(k-1)X}{k-2},\\[6pt]
\displaystyle \frac{k-4}{k+2}<X\leqslant \frac{k-4}{k},\quad \frac{1+kX}{k-1} < Y\leqslant \frac{1+(k-1)X}{k-2},\\[6pt]
\displaystyle \frac{k-4}{k}<X\leqslant \frac{k-3}{k+1},\quad \frac{1+kX}{k-1} < Y\leqslant \frac{1+X}{2}
\end{cases}
\end{equation}
and appears to be a quadrangle with vertices
\begin{multline*}
P=\biggl(\frac{k-5}{k+1},\frac{k-3}{k+1}\biggr),\quad Q=\biggl(\frac{k-4}{k},\frac{k-2}{k}\biggr),\quad R=\biggl(\frac{k-3}{k+1},\frac{k-1}{k+1}\biggr),\\ S=\biggl(\frac{k-4}{k+2},\frac{k-2}{k+2}\biggr).
\end{multline*}
Passing to $\mathcal{T}(2,1,k,1)$, we should add to (\ref{lab_44}) the inequalities $f(X)<Y\leqslant g(X)$, where
\[
f(X) = f(X;2,1,k,1) = \frac{1+X(2k-3)}{2k-5},\quad g(X) = g(X;2,1,k,1) = \frac{1+X(k-2)}{k-3}.
\]
A direct calculation shows that
\begin{equation}\label{lab_45}
f(X)\leqslant \frac{1+2X}{3}\quad\text{for}\quad X\leqslant \frac{2k-8}{2k+1}\quad\text{and}\quad f(X)\leqslant \frac{1+kX}{k-1}\quad\text{for}\quad X\leqslant \frac{k-4}{3}
\end{equation}
Since
\[
\frac{2k-8}{2k+1}\geqslant \frac{k-4}{k+2} = x_{S},\quad \frac{k-4}{3}\geqslant \frac{k-3}{k+1} = x_{R}\quad\text{for all}\quad k\geqslant 5,
\]
then the inequalities (\ref{lab_45}) hold on the segments $x_{P}\leqslant X\leqslant x_{S}$ and $x_{S}\leqslant X\leqslant x_{R}$ respectively
 (the sense of notations $x_{P}, x_{S}, x_{R}$ is clear). Hence, the region $\mathcal{T}(2,1,k)$ lies entirely over the line $Y = f(X)$.
 Similarly, one can check that
\[
\frac{1+(k-1)X}{k-2}\leqslant g(X)\quad\text{for}\quad X\geqslant -1\quad\text{and}\quad \frac{1+X}{2}\leqslant g(X)\quad\text{for}\quad X\geqslant \frac{k-5}{k-1}.
\]
Since $(k-5)/(k-1)\leqslant x_{Q}$ then $\mathcal{T}(2,1,k)$ lies entirely under the line $Y = g(X)$. Hence, $\mathcal{T}(2,1,k,1) = \mathcal{T}(2,1,k)$ for all $k\geqslant 9$.

(e) By (d), for any $k\geqslant 9$, the region $\mathcal{T}(2,1,k,1)$ is defined by the system (\ref{lab_44}). Passing to $\mathcal{T}(2,1,k,1,2)$, we have to add an additional inequalities $f(X)<Y\leqslant g(X)$ where
\begin{multline*}
f(X) = f(X;2,1,k,1,2) = \frac{1+X(2k-5)}{2k-7},\quad g(X) = g(X;2,1,k,1,2) = \frac{1+X(k-3)}{k-4}.
\end{multline*}
The direct calculation shows that
\begin{equation}\label{lab_46}
f(X)\leqslant \frac{1+2X}{3}\quad\text{for}\quad X\leqslant \frac{2k-10}{2k-1}\quad\text{and}\quad f(X)\leqslant \frac{1+(k-1)X}{k-2}\quad\text{for}\quad X\leqslant \frac{k-5}{3}.
\end{equation}
Since
\[
\frac{2k-10}{2k-1}\geqslant \frac{k-4}{k+2} = x_{S},\quad \frac{k-5}{3}\geqslant \frac{k-1}{k+1} = x_{R}\quad\text{for all}\quad k\geqslant 8,
\]
then the inequalities (\ref{lab_46}) hold on the segments $x_{P}\leqslant X\leqslant x_{S}$ and $x_{S}\leqslant X\leqslant x_{R}$ correspondingly. This means that $\mathcal{T}(2,1,k,1)$ lies entirely under the line $Y = f(X)$. Similarly,
\[
\frac{1+(k-1)X}{k-2}\leqslant g(X)\quad\text{for}\quad X\geqslant -1\quad\text{and}\quad \frac{1+X}{2}\leqslant g(X)\quad\text{for}\quad X\geqslant \frac{k-6}{k-2}.
\]
Next, since $(k-6)/(k-2)\leqslant x_{Q}$ then $\mathcal{T}(2,1,k,1)$ lies entirely over the line $Y = g(X)$. Therefore, $\mathcal{T}(2,1,k,1,2) = \mathcal{T}(2,1,k,1) = \mathcal{T}(2,1,k)$ for $k\geqslant 9$.

(f) Since $\mathcal{T}(1,k,1) = \mathcal{T}(1)\bigcap T^{-1}\bigl(\mathcal{T}(k,1)\bigr)$ then
\[
T\mathcal{T}(1,k,1) = T \mathcal{T}(1)\bigcap \mathcal{T}(k,1)\subseteq \mathcal{T}(k,1).
\]
Since $T$ preserves an area, we have
\begin{equation}\label{lab_47}
|\mathcal{T}(1,k,1)| = |T\bigl(\mathcal{T}(1,k,1)\bigr)|\leqslant |\mathcal{T}(k,1)|.
\end{equation}
 Next, since $\mathcal{T}(1,k,1,2)\subseteq \mathcal{T}(1,k,1)$ then
\begin{equation}\label{lab_48}
|\mathcal{T}(1,k,1)|\geqslant |\mathcal{T}(1,k,1,2)|.
\end{equation}
In view of (a),(d) we have $|\mathcal{T}(1,k,1,2)| = |\mathcal{T}(k,1)| = |\mathcal{T}(k)|$ for $k\geqslant 9$. Together with (\ref{lab_47}), (\ref{lab_48}), this means that $|\mathcal{T}(1,k,1)| = \mathcal{T}(k)$ for all $k\geqslant 9$. For $6\leqslant k\leqslant 8$, this equality is established by the direct calculation. Lemma is proved.
\vspace{0.3cm}

\textsc{Lemma 13.} \emph{Suppose that} $r\geqslant 1$, $k\geqslant 4r+2$. \emph{Then} $T^{\,-j}(\mathcal{T}(k))\subset \mathcal{T}(2)$ and $T^{\,j}(\mathcal{T}(k))\subset \mathcal{T}(2)$ \emph{for any} $j$, $2\leqslant j\leqslant r$.
\vspace{0.5cm}

\textsc{Lemma 14.} \emph{Suppose that} $r\geqslant 1$, $\max{\{k_{1},\ldots,k_{r}\}} = k_{j}\geqslant 4r+2$, \emph{and let} $\mathcal{T}(k_{1},\ldots,k_{r}) \ne \varnothing$. \emph{Then the tuple} $(k_{1},\ldots,k_{r})$ \emph{satisfies the following conditions:}
\begin{equation}\label{lab_49}
k_{i} = 1\quad\textit{for} \quad |i-j|=1,\quad k_{i} = 2\quad\textit{for} \quad |i-j|\geqslant 2.
\end{equation}

Lemmas 13, 14 follow from the assertions (\textit{i}), (\textit{ii}) of \cite[Lemma 3.4]{Boca_Cobeli_Zaharescu_2003} (see also \cite[p.~570]{Boca_Cobeli_Zaharescu_2003}).
\vspace{0.5cm}

\textsc{Lemma 15.} \emph{The region $\mathcal{T}(k,m,n)$ is empty iff\;\footnote{To make this cumbersome assertion a little more clear, we do not try to make the intersection of the lists (1) and (2a),\ldots, (6) empty.}} \\ [6pt]
\begin{tabular}{l p{0.5\linewidth}}
(1) \emph{if $(k,m)$ or $(m,n)$ belongs to the list of Lemma 10}; & \\ [6pt]
\end{tabular}
\begin{tabular}{l p{0.5\linewidth}}
(2a) $m=1$, $k=2$, $1\leqslant n\leqslant 5$;  & (2b) $m=1$, $k=3$, $n\ne 4,5,6,7,8$;\\ [6pt]
(2c) $m=1$, $k=4$, $n\ne 3,4,5$;  & (2d) $m=1$, $k=5$, $n\ne 3,4$; \\ [6pt]
(2e) $m=1$, $6\leqslant k\leqslant 8$, $n\geqslant 4$; & (2f) $m=1$, $k\geqslant 9$, $n\geqslant 3$; \\ [6pt]
(3a) $m=2$, $k=1$, $n\ne 2,3,4$; & (3b) $m=2$, $k=2$, $n\geqslant 4$; \\ [6pt]
(3c) $m=2$, $k=3$, $n\geqslant 3$; & (3d) $m=2$, $k=4$, $n\geqslant 2$; \\ [6pt]
(4a) $m=3$, $k=1$, $n\geqslant 3$; & (4b) $m=3$, $k=2$, $n\geqslant 3$; \\ [6pt]
(5a) $m=4$, $k=1$, $n\geqslant 3$; & (5b) $m=4$, $k=2$, $n\geqslant 2$; \\ [6pt]
(6)  $m\geqslant 5$, $k=1$, $n\geqslant 2$. & \\ [6pt]
\end{tabular}

\textsc{Proof.} We start from the case $m=2$. By Lemma 10, $\mathcal{T}(2,n)$ and $\mathcal{T}(k,2)$ are empty for $k,n\geqslant 5$. Hence, it is sufficient to check the regions $\mathcal{T}(k,2,n)$ with $1\leqslant k,n\leqslant 4$. The direct calculation shows that $\mathcal{T}(k,2,n)$ is empty for the tuples
\[
(1,2,1),\quad (2,2,4),\quad (3,2,3),\quad (3,2,4),\quad (4,2,2),\quad (4,2,3),\quad (4,2,4).
\]
Further, let $3\leqslant m\leqslant 4$. Since $\mathcal{T}(3,n)$ and $\mathcal{T}(4,n)$ are empty for $n\geqslant 3$, it is sufficient to check the cases $1\leqslant k,n\leqslant 2$. The calculation shows that all such regions are non-empty.

Let $m\geqslant 5$. Since $\mathcal{T}(m,n)$ is empty for $n\ne 1$, then $\mathcal{T}(k,m,n)$ is empty in all cases when $(k,n)\ne (1,1)$. By Lemma 12 (f), $|\mathcal{T}(1,m,1)| = |\mathcal{T}(m)|\ne 0$ for any $m\geqslant 6$. Finally, the direct calculation shows that $|\mathcal{T}(1,5,1)|\ne 0$.

Collecting together all the above assertions, we get pp. (3)-(6) of the lemma.

Suppose now that $m=1$. Since both $\mathcal{T}(k,m,n)$ and $\mathcal{T}(n,m,k)$ are empty or non-empty, it is sufficient to check the pairs $(k,n)$ with $2\leqslant n\leqslant k$.

If $k\geqslant 4\cdot 3+2=14$ then Lemma 14 implies that $n=2$. In this case, Lemma 12 (b) implies that $\mathcal{T}(k,1,2)\neq \varnothing$. For $2\leqslant k\leqslant 13$, the direct calculation shows that $\mathcal{T}(k,1,n)$ is empty iff the triple $(k,1,m)$ (or the triple symmetric to it) belongs to the list \\

\begin{tabular}{l p{0.5\linewidth}}
\hspace{-8mm} $(2,1,2),\quad (3,1,2),\quad (3,1,3),\quad (4,1,2), \quad (5,1,2),\quad (5,1,5)$; & \\ [3pt]
\end{tabular}
\begin{tabular}{l p{0.5\linewidth}}
$(6,1,n)$, $4\leqslant n \leqslant 6$; & $(10,1,n)$, $3\leqslant n\leqslant 10$; \\ [3pt]
$(7,1,n)$, $4\leqslant n \leqslant 7$; & $(11,1,n)$, $3\leqslant n\leqslant 11$; \\ [3pt]
$(8,1,n)$, $4\leqslant n \leqslant 8$; & $(12,1,n)$, $3\leqslant n\leqslant 12$; \\ [3pt]
$(9,1,n)$, $3\leqslant n \leqslant 9$; & $(13,1,n)$, $3\leqslant n\leqslant 13$. \\ [6pt]
\end{tabular}

Lemma is proved. $\square$

\vspace{0.5cm}

\textsc{Lemma 16.} \emph{Suppose that} $r\geqslant 1$, $Q>2(r+1)$, \emph{and let} $\max{\{k_{1},\ldots,k_{r}\}} = k_{j}\geqslant 4r+2$. \emph{Suppose also that} $\mathcal{T}(k_{1},\ldots,k_{r})$ \emph{is non-empty. Then } $k_{j}\leqslant 2Q$.
\vspace{0.5cm}

For the proof of this assertion, see \cite[p. 384; Remark 2.3]{Boca_Gologan_Zaharescu_2002}.

\section{Expression for $\boldsymbol{N(Q;r,D,c_{0})}$ in terms of areas of polygons}

The asymptotic expression for $N_{e_{0},e_{1}}(Q\omega)$, $\omega = \mathcal{T}(k_{1},\ldots, k_{r})$ is based on the following version of Jarnik inequality (see \cite[Lemma 1]{Ustinov_2009}): suppose that $\Omega$ is simply-connected planar region with rectifiable boundary and denote by $|\Omega|$ and $|\partial \Omega|$ its area and perimeter. Then the number $n(\Omega)$ of lattice points lying in the interior of $\Omega$ satisfies the inequality
\[
\bigl|n(\Omega) - |\Omega|\bigr|<4\bigl(|\partial \Omega|+1\bigr).
\]
In particular, if $\Omega = Q\omega$ then
\begin{equation}\label{lab_50}
n(\Omega) = Q^{2}|\omega| + O\bigl(Q|\partial \omega|+1\bigr)
\end{equation}
where the implied constant is absolute.
\vspace{0.3cm}

In what follows, we denote
\[
\kappa_{D} = \sum\limits_{\text{GCD}(d,D)=1}\frac{\mu(d)}{d^{2}} = \frac{6}{\pi^{2}}\,\prod\limits_{p|D}\biggl(1-\frac{1}{p^{2}}\biggr)^{-1}.
\]

\textsc{Lemma 17.} \emph{Suppose that $c_{0}, c_{1}, D$ are integers such that} $D\geqslant 2$, $0\leqslant c_{0}, c_{1}\leqslant D-1$, $\text{GCD}(D,c_{0},c_{1})=1$. \emph{Further, let $r\geqslant 1$, $k_{1},\ldots, k_{r}\geqslant 1$ be integers such that $\omega = \mathcal{T}(k_{1},\ldots, k_{r})$ is non-empty. Finally, let $\max{\{k_{1},\ldots, k_{r}\}} = k_{j}$ and $Q>2k_{j}$. Then the number $N_{c_{0},c_{1}}(Q\omega)$ of primitive lattice points lying in $Q\omega$ and satisfying to the conditions $q_{0}\equiv c_{0}\pmod{D}$, $q_{1}\equiv c_{1}\pmod{D}$ is given by the following expression:}
\[
N_{c_{0},c_{1}}(Q\omega) = \frac{\kappa_{D}}{D^{2}}\,Q^{2}\,|\omega|\,+\,O\bigl(Q(\rho_{1}+\rho_{2}+\rho_{3})\bigr),
\]
\emph{where}
\[
\rho_{1} = \frac{Q}{k_{j}},\quad \rho_{2} = \frac{Q(\ln{Q})}{D}\mathbbm{1}_{\{k_{j}<4r+2\}},\quad \rho_{3} = \frac{Q(\ln{Q})}{Dk_{j}}\mathbbm{1}_{\{k_{j}\geqslant 4r+2\}}.
\]
\emph{and the implied constant is absolute\footnote{Here and in what follows, we assume that $\mathbbm{1}_{A}=1$ if the condition $A$ holds true and $\mathbbm{1}_{A}=0$ otherwise.}.}
\vspace{0.3cm}

\textsc{Proof.} Suppose that $2\leqslant j\leqslant r$ and denote the region $T^{\,j-1}(\omega)$ by $\omega_{j}$. In view of Lemma 6, $N_{c_{0},c_{1}}(Q\omega) = N_{\sigma_{0},\sigma_{1}}(Q\omega_{j})$, where
\begin{equation*}
\begin{cases}
\sigma_{0}\equiv c_{1}\,\mathbb{K}_{j-2}(k_{1},\ldots,k_{j-2}) - c_{0}\,\mathbb{K}_{j-3}(k_{2},\ldots,k_{j-2})\pmod{D},\\
\sigma_{1}\equiv c_{1}\,\mathbb{K}_{j-1}(k_{1},\ldots,k_{j-1}) - c_{0}\,\mathbb{K}_{j-2}(k_{2},\ldots,k_{j-1})\pmod{D}.
\end{cases}
\end{equation*}
Further, let $\Omega_{j}$ be the set of all points
\[
\biggl(\frac{1}{D}\bigl(Qx-\sigma_{0}\bigr),\frac{1}{D}\bigl(Qy-\sigma_{1}\bigr)\biggr),\quad\text{where}\quad (x,y)\in \omega_{j}.
\]
Then $N_{\sigma_{0},\sigma_{1}}(Q\omega_{j})$ equals to the number of lattice points $(m,n)$ of the region $Q\Omega_{j}$  satisfying the condition
\begin{equation}\label{lab_51}
\text{GCD}(Dm+\sigma_{0},Dn+\sigma_{1}) = 1.
\end{equation}
Hence,
\[
N_{\sigma_{0},\sigma_{1}}(Q\omega_{j}) = \sum\limits_{\substack{(m,n)\in Q\Omega_{j} \\ \text{GCD}(Dm + \sigma_{0},Dn + \sigma_{1})=1}}1 =
\sum\limits_{(m,n)\in Q\Omega_{j}}\;\;\sum\limits_{\delta | \text{GCD}(Dm + \sigma_{0},Dn + \sigma_{1})}\mu(\delta).
\]
In order to estimate the maximal value of $\delta$, we note that
\[
\omega = \mathcal{T}(k_{1})\cap T^{-1}\mathcal{T}(k_{2})\cap \cdots \cap T^{\,-(r-1)}\mathcal{T}(k_{r})\subseteq T^{\,-(j-1)}\mathcal{T}(k_{j}),
\]
so we have $\omega_{j} = T^{j-1}(\omega) \subseteq\mathcal{T}(k_{j})$. By Lemma 10, there are no two neighbouring units in the tuple $(k_{1},\ldots,k_{r})$. Hence, $k_{j}\geqslant 2$. By Lemma 5, $\mathcal{T}(k_{j})$ is a quadrangle with vertices $A,B,C$ and $D$. Therefore, $Q\cdot \omega_{j}$ is entirely contained in the rectangle
\[
Q\cdot\frac{k_{j}-1}{k_{j}+1}\leqslant x\leqslant Q,\quad 0<y\leqslant \frac{2Q}{k_{j}}.
\]
Thus for any lattice point $(u,v)\in Q\cdot\omega_{j}$ we have $1\leqslant v\leqslant 2Q/k_{j}$.
Therefore, for any lattice point $(m,n)\in \Omega_{j}$ we get $1\leqslant Dn+\sigma_{1}\leqslant 2Q/k_{j}$. This means that $\delta\leqslant 2Q/k_{j}$ and
\[
N_{\sigma_{0},\sigma_{1}}(Q\omega_{j}) = \sum\limits_{1\leqslant\delta\leqslant 2Q/k_{j}}\mu(\delta)\;\;\sum\limits_{\substack{(m,n)\in Q\Omega_{j} \\ Dm + \sigma_{0}\equiv  0\;(\mmod{\delta}) \\ Dn + \sigma_{1}\equiv  0\;(\mmod{\delta})}}1 = \sum\limits_{\substack{1\leqslant\delta\leqslant 2Q/k_{j} \\ \text{GCD}(\delta,D)=1}}\mu(\delta)\;\;
\sum\limits_{\substack{(m,n)\in Q\Omega_{j} \\ m \equiv  \sigma_{0}\;(\mmod{\delta}) \\ n \equiv  \sigma_{1}\;(\mmod{\delta})}}1,
\]
where $\sigma_{j}\equiv -\varepsilon_{j}D^{*}\pmod{\delta}$, $j = 0,1$. It is easy to see that the inner sum is equal to the number of lattice point lying in the region
$Q\Omega_{j}^{*}$, which can be obtained from $\Omega_{j}$ by the shift to the vector $(-\sigma_{0}/Q,-\sigma_{1}/Q)$ and by the horizontal and vertical dilations with coefficient $\delta$. Hence, in view of (\ref{lab_50}), this sum equals to
\[
|\Omega_{j}^{*}| + O\bigl(|\partial(\Omega_{j}^{*})|+1\bigr) = \frac{|\Omega_{j}|}{\delta^{2}}\,+\,O\biggl(\frac{|\partial\Omega_{j}|}{\delta}+1\biggr) =
\frac{Q^{2}|\omega_{j}|}{(\delta D)^{2}}\,+\,O\biggl(\frac{Q}{\delta D}|\partial\omega_{j}|+1\biggr).
\]
Since $T$ preserves an area, then $|\omega_{j}| = |\omega|$. Further, if $k_{j}\geqslant 4r+2$, then, in view of Lemma 14, the tuple $(k_{1},\ldots,k_{r})$ satisfies to (\ref{lab_49}).
Next, let us show that $\omega$ coincides with $T^{-(j-1)}\bigl(\mathcal{T}(k_{j})\bigr)$. Indeed, if $|i-j|=1$ then $k_{i}=1$. Consider the region
\begin{multline*}
\omega^{(1)} = T^{-\,(j-2)}\mathcal{T}(k_{j-1})\bigcap T^{-\,(j-1)}\mathcal{T}(k_{j}) = T^{-\,(j-2)}\mathcal{T}(1)\bigcap T^{-\,(j-1)}\mathcal{T}(k_{j}) = \\
= T^{-(j-2)}\bigl(\mathcal{T}(1)\bigcap T^{-1}\mathcal{T}(k_{j})\bigr).
\end{multline*}
By Lemma 9, $T^{-1}\mathcal{T}(k_{j})$ is contained entirely in $\mathcal{T}(1)$. Hence,
\[
\omega^{(1)} = T^{-\,(j-2)}\circ T^{-1}\mathcal{T}(k_{j}) = T^{-\,(j-1)}\mathcal{T}(k_{j}).
\]
Similarly, if we denote by $\omega^{(2)}$ the intersection
\[
T^{-\,(j-1)}\mathcal{T}(k_{j})\bigcap T^{-j}\mathcal{T}(k_{j+1})
\]
then
\[
\omega^{(2)} = T^{-(j-1)}\mathcal{T}(k_{j})\bigcap T^{-j}\mathcal{T}(1) = T^{-j}\bigl(T\mathcal{T}(k_{j})\bigcap \mathcal{T}(1)\bigr) = T^{-j}\circ T\mathcal{T}(k_{j}) = T^{-(j-1)}\mathcal{T}(k_{j}).
\]
Further, let $j\geqslant 3$ and $1\leqslant i\leqslant j-2$. Setting
\[
\omega^{(3)} = T^{-(i-1)}\mathcal{T}(k_{i})\bigcap T^{-(j-1)}\mathcal{T}(k_{j})
\]
we get
\[
\omega^{(3)} = T^{-(i-1)}\mathcal{T}(2)\bigcap T^{-(j-1)}\mathcal{T}(k_{j}) = T^{-(i-1)}\bigl(\mathcal{T}(2)\bigcap T^{-(j-i)}\mathcal{T}(k_{j})\bigr).
\]
Since $j-i\geqslant 2$, then Lemma 13 implies that $T^{-(j-i)}\mathcal{T}(k_{j})\subseteq \mathcal{T}(2)$. Hence,
\[
\omega^{(3)} = T^{-(i-1)}\circ T^{-(j-i)}\mathcal{T}(k_{j}) = T^{-(j-1)}\mathcal{T}(k_{j}).
\]
Finally, if $j\leqslant r-2$, $j+2\leqslant i\leqslant r$ and
\[
\omega^{(4)} = T^{-(j-1)}\mathcal{T}(k_{j})\bigcap T^{-(i-1)}\mathcal{T}(k_{i})
\]
then one can check that $\omega^{(4)} = T^{-(j-1)}\mathcal{T}(k_{j})$.

Therefore,
\[
\omega = T^{-(j-1)}\mathcal{T}(k_{j}),\quad \omega_{j} = T^{j-1}(\omega) = \mathcal{T}(k_{j}).
\]
Hence,
\[
|\partial \omega_{j}| = |\partial \mathcal{T}(k_{j})|\ll \frac{1}{k_{j}}.
\]
If $k_{j}<4r+2$ then
\[
|\partial \omega_{j}| \ll 1.
\]
Thus, for any $k_{j}$ we have
\begin{equation*}
|\partial \omega_{j}| \ll \mathbbm{1}_{\{k_{j}<4r+2\}} + \frac{\mathbbm{1}_{\{k_{j}\geqslant 4r+2\}}}{k_{j}}.
\end{equation*}
This means that
\[
|\Omega_{j}^{*}| = \frac{Q^{2}|\omega|}{(\delta D)^{2}}\, + O\biggl(\frac{Q}{\delta D}\biggl(\mathbbm{1}_{\{k_{j}<4r+2\}} + \frac{\mathbbm{1}_{\{k_{j}\geqslant 4r+2\}}}{k_{j}}\biggr)+1\biggr).
\]
Substituting this relation into the formula for $N_{\sigma_{0},\sigma_{1}}(Q\omega_{j})$ and using the equality
\[
\sum\limits_{\substack{\delta\leqslant 2Q/k_{j} \\ \text{GCD}(\delta,D)=1}}\frac{\mu(\delta)}{\delta^{2}} = \kappa_{D}+O\biggl(\frac{k_{j}}{Q}\biggr)
\]
together with the bound $|\omega|\leqslant |\mathcal{T}(k_{j})|\ll k_{j}^{-3}$, we get to the desired formula for $N_{c_{0},c_{1}}(Q\omega)$. $\square$
\vspace{0.3cm}

The following assertion states the general formula for $N(Q;r,D,c_{0})$.
\vspace{0.3cm}

\textsc{Lemma 18.} \emph{Uniformly in $Q$, $1\leqslant r\leqslant Q/3$, $2\leqslant D\leqslant Q/3$, $0\leqslant c_{0}\leqslant D-1$ one has}
\[
N(Q;r,D,c_{0}) = \mathfrak{c}_{r}\,\frac{\kappa_{D}}{D^{2}}\,Q^{2} + O\biggl(Q\Bigl(1+\frac{\ln{Q}}{D}\Bigr)\bigl(\Sigma_{0}+\Sigma_{1}\bigr)+\frac{Q^{2}}{D^{2}}\Sigma_{2}\biggr),
\]
\emph{where}
\begin{multline*}
\mathfrak{c}_{r} = \sum\limits_{c_{1}}\sum\limits_{\mathbf{k}\in\, \mathcal{A}_{r}^{\circ}(D,c_{0},c_{1})}|\mathcal{T}(\mathbf{k})|, \quad
\Sigma_{0} = \sum\limits_{c_{1}}\sum\limits_{\substack{\mathbf{k}\in\, \mathcal{A}_{r}^{\circ}(c_{0},c_{1},D) \\ \|\mathbf{k}\|\leqslant 4r+1}}1, \\
\Sigma_{1} = \sum\limits_{j=1}^{r}U_{j},\quad \Sigma_{2} = \sum\limits_{j=1}^{r}V_{j},\quad U_{j} = \sum\limits_{c_{1}}\sum\limits_{\substack{\mathbf{k}\in\,\mathcal{A}_{r}^{\circ}(D,c_{0},c_{1}) \\ 4r+2\leqslant \|\mathbf{k}\|= k_{j} \leqslant 2Q}}\frac{1}{k_{j}}, \quad V_{j} = \sum\limits_{c_{1}}\sum\limits_{\substack{\mathbf{k}\in\,\mathcal{A}_{r}^{\circ}(D,c_{0},c_{1}) \\ \|\mathbf{k}\|= k_{j} > 2Q}}\frac{1}{k_{j}^{3}}.
\end{multline*}
\emph{the sign $\displaystyle \sum\limits_{c_{1}}$ means the summation over} $0\leqslant c_{1}\leqslant D-1$, $c_{1}\ne c_{0}$, $\text{GCD}(c_{0},c_{1},D)=1$, \emph{for} $\mathbf{k} = (k_{1},\ldots,k_{r})$, the symbol $\|\mathbf{k}\|$ \emph{means} $\max{\{k_{1},\ldots,k_{r}\}}$, \emph{and the constant in $O$-symbol is absolute.}
\vspace{0.3cm}

\textsc{Proof.} By the Corollary of Lemma 3,
\[
N(Q;r,D,c_{0}) = \sum\limits_{c_{1}}\sum\limits_{\mathbf{k}\in\, \mathcal{A}_{r}^{\circ}(D,c_{0},c_{1})}N_{c_{0},c_{1}}(Q\omega),\quad \omega = \mathcal{T}(\mathbf{k}) = \mathcal{T}(k_{1},\ldots,k_{r}).
\]
By Lemma 16, the non-zero contribution to the above sum comes from the tuples $\mathbf{k}$ with $\|\textbf{k}\| \leqslant 2Q$. Using the notations of Lemma 17, we get
\begin{multline*}
N(Q;r,D,c_{0}) = \sum\limits_{\substack{0\leqslant c_{1}\leqslant D-1 \\ c_{1}\ne c_{0} \\ \text{GCD}(D,c_{0},c_{1})=1}}\sum_{\substack{\mathbf{k}\in \mathcal{A}_{r}^{\circ}(D,c_{0},c_{1}) \\ \|\mathbf{k}\| = k_{j}\leqslant 2Q}}\biggl(\frac{\kappa_{D}Q^{2}}{D^{2}}|\mathcal{T}(\mathbf{k})| + O\bigl(Q(\rho_{1}+ \rho_{2}+\rho_{3})\bigr)\biggr) \\
= \frac{\kappa_{D}Q^{2}}{D^{2}}\sum\limits_{c_{1}}\sum_{\substack{\mathbf{k}\in \mathcal{A}_{r}^{\circ}(D,c_{0},c_{1}) \\ \|\mathbf{k}\| = k_{j}\leqslant 2Q}}
|\mathcal{T}(\mathbf{k})| + O\bigl(Q(R_{1}+R_{1}+R_{3})\bigr),
\end{multline*}
where $R_{s}$ means the contribution coming from the term $\rho_{s}$, $s = 1,2,3$.

Passing to the estimates of $R_{s}$ we note that $r\leqslant Q/3 < (Q-1)/2$ and hence $4r+2<2Q$. Thus we find:
\begin{multline*}
R_{1} = \sum\limits_{c_{1}}\sum_{\substack{\mathbf{k}\in \mathcal{A}_{r}^{\circ}(D,c_{0},c_{1}) \\ \|\mathbf{k}\|  = k_{j}\leqslant 2Q}}\frac{1}{k_{j}} =
\sum\limits_{c_{1}}\biggl(\;\sum_{\substack{\mathbf{k}\in \mathcal{A}_{r}^{\circ}(D,c_{0},c_{1}) \\ \|\mathbf{k}\|  = k_{j}\leqslant 4r+1}} + \sum_{\substack{\mathbf{k}\in \mathcal{A}_{r}^{\circ}(D,c_{0},c_{1}) \\ 4r+2\leqslant \|\mathbf{k}\|  = k_{j}\leqslant 2Q}} \biggr)\frac{1}{k_{j}} \leqslant \\
\leqslant \sum\limits_{c_{1}}\sum_{\substack{\mathbf{k}\in \mathcal{A}_{r}^{\circ}(D,c_{0},c_{1}) \\ \|\mathbf{k}\|  = k_{j}\leqslant 4r+1}}1 + \sum\limits_{j=1}^{r}\sum\limits_{c_{1}}
\sum_{\substack{\mathbf{k}\in \mathcal{A}_{r}^{\circ}(D,c_{0},c_{1}) \\ 4r+2\leqslant \|\mathbf{k}\|  = k_{j}\leqslant 2Q}}\frac{1}{k_{j}}\leqslant \Sigma_{0} + \Sigma_{1}.
\end{multline*}
Similarly,
\[
R_{1}\leqslant \frac{\ln{Q}}{D}\sum\limits_{c_{1}}\sum_{\substack{\mathbf{k}\in \mathcal{A}_{r}^{\circ}(D,c_{0},c_{1}) \\ \|\mathbf{k}\|  = k_{j}\leqslant 4r+1}}1 \leqslant
\frac{\ln{Q}}{D}\Sigma_{0}.
\]
Next,
\[
R_{3}\leqslant \frac{\ln{Q}}{D}\sum\limits_{j=1}^{r}\sum\limits_{c_{1}}\sum_{\substack{\mathbf{k}\in \mathcal{A}_{r}^{\circ}(D,c_{0},c_{1}) \\ 4r+2\leqslant \|\mathbf{k}\|  = k_{j}\leqslant 2Q}}\frac{1}{k_{j}}\leqslant \frac{\ln{Q}}{D}\Sigma_{1}.\quad \square
\]
Further, the contribution from the main terms has the form
\[
\frac{\kappa_{D}Q^{2}}{D^{2}}\sum\limits_{c_{1}}\sum\limits_{\substack{\mathbf{k}\in A_{r}^{\circ}(D,c_{0},c_{1}) \\ k_{j}\leqslant 2Q}}|\mathcal{T}(\mathbf{k})|.
\]
If we omit the restriction $k_{j}\leqslant 2Q$ in this sum then the error does not exceed
\[
\frac{Q^{2}}{D^{2}}\sum\limits_{j=1}^{r}\sum\limits_{c_{1}}\sum\limits_{\substack{\mathbf{k}\in A_{r}^{\circ}(D,c_{0},c_{1}) \\ \|\mathbf{k}\| = k_{j}> 2Q}}|\mathcal{T}(\mathbf{k})|.
\]
Since $\mathcal{T}(\mathbf{k}) \subseteq T^{-(j-1)}\mathcal{T}(k_{j})$ then, by Lemma 7,
\[
|\mathcal{T}(\mathbf{k})|\leqslant |T^{-(j-1)}\mathcal{T}(k_{j})| = |\mathcal{T}(k_{j})|\ll k_{j}^{-3}.
\]
Therefore, the above sum is less than
\[
\frac{Q^{2}}{D^{2}}\Sigma_{2},\quad \Sigma_{2} \ll \sum\limits_{j=1}^{r}V_{j}.
\]
Thus we find
\[
N(Q;r,D,c_{0}) = \mathfrak{c}_{r}\,\frac{\kappa_{D}}{D^{2}}\,Q^{2} + R,\quad R\ll Q\biggl(1+\frac{\ln{Q}}{D}\biggr)\bigl(\Sigma_{0}+\Sigma_{1}\bigr)+\frac{Q^{2}}{D^{2}}\Sigma_{2}.\quad \square
\]

\section{Precise formulas for the continuants and the pro\-por\-ti\-on $\boldsymbol{\nu(Q;r,D,c_{0})}$}

In order to derive uniform bounds for the remainder in the expression for $N(Q;r,D,c_{0})$ and for $\nu(Q;r,D)$ in the case $D=2,3$, we need some explicit expressions for the continuants of special type. For brevity, in this section, we use the notation $\mathbb{K}(k_{1},\ldots,k_{r})$ instead of $\mathbb{K}_{r}(k_{1},\ldots,k_{r})$.
\vspace{0.3cm}

\textsc{Lemma 19.} \emph{Let $n\geqslant 1$. Then, setting}
\[
A_{n} = \mathbb{K}\bigl((4,1)^{n}\bigr),\quad B_{n} = \mathbb{K}\bigl((4,1)^{n},4\bigr),\quad C_{n} = \mathbb{K}\bigl(1,(4,1)^{n}\bigr),\quad D_{n} = \mathbb{K}\bigl(2^{n}\bigr),
\]
\emph{we get:} $A_{n} = 2n+1$, $B_{n} = 4(n+1)$, $C_{n} = D_{n} = n+1$.
\vspace{0.3cm}

\textsc{Proof.} Using the relations (\ref{lab_07}), (\ref{lab_10}) we find
\begin{align*}
A_{n+1} & = \mathbb{K}\bigl((4,1)^{n},4,1\bigr) = \mathbb{K}\bigl((4,1)^{n},4\bigr)-\mathbb{K}\bigl((4,1)^{n}\bigr) = B_{n}-A_{n},\\
B_{n+1} & = \mathbb{K}\bigl((4,1)^{n+1},4\bigr) = \mathbb{K}\bigl((4,1)^{n},4,1,4\bigr) = \mathbb{K}(4,1)\mathbb{K}\bigl((4,1)^{n},4\bigr) - \mathbb{K}(4)\mathbb{K}\bigl((4,1)^{n}\bigr) = \\
& = 3B_{n}-4A_{n}.
\end{align*}
Further, (\ref{lab_09}) implies
\begin{align*}
& C_{n+1} = \mathbb{K}\bigl(1,(4,1)^{n},4,1\bigr) = \mathbb{K}(1,4)\mathbb{K}\bigl(1,(4,1)^{n}\bigr)-\mathbb{K}(1)\mathbb{K}\bigl(1,(4,1)^{n-1},4\bigr) = 3C_{n}-A_{n},\\
& D_{n+1} = 2\mathbb{K}(2^{n}) - \mathbb{K}(2^{n-1}) = 2D_{n}-D_{n-1}.
\end{align*}
Finally, using the relations $A_{1} = 3$, $B_{1} = 8$, $C_{1} = D_{1} = 2$, we prove the desired assertion by induction. $\square$
\vspace{0.3cm}

\textsc{Lemma 20.} \emph{Let $a,b,c$ and $d$ be any numbers, and suppose that $m,n\geqslant 0$ are integers. Then the following identities hold true:}
\begin{equation}\label{lab_52}
\begin{cases}
& \mathbb{K}(a,2^{n}) = \mathbb{K}(2^{n},a) = (a-1)n+a;\\
& \mathbb{K}(a,2^{n},b) = (a-1)(b-1)n+ab-1;\\
& \mathbb{K}(2^{n},b,c) = (bc-c-1)n+bc-1;\\
& \mathbb{K}(a,2^{n},b,c) = (a-1)(bc-c-1)n+abc-a-c;
\end{cases}
\end{equation}

\begin{equation}\label{lab_53}
\begin{cases}
& \mathbb{K}(a,(4,1)^{n}) = (2a-1)n+a;\\
& \mathbb{K}((4,1)^{n},b) = 2(b-2)n+b;\\
& \mathbb{K}(a,(4,1)^{n},b) = (2a-1)(b-2)n+ab-1;\\
& \mathbb{K}((4,1)^{n},b,c) = 2(bc-2c-1)n+bc-1; \\
& \mathbb{K}(a,(4,1)^{n},b,c) = (2a-1)(bc-2c-1)n+abc-a-c;
\end{cases}
\end{equation}

\begin{equation}\label{lab_54}
\begin{cases}
& \mathbb{K}(a,(1,4)^{n}) = 2(a-2)n+a;\\
& \mathbb{K}((1,4)^{n},b) = (2b-1)n+b;\\
& \mathbb{K}(a,(1,4)^{n},b) = (a-2)(2b-1)n+ab-1;
\end{cases}
\end{equation}

\begin{equation}\label{lab_55}
\begin{cases}
& \mathbb{K}((4,1)^{n},3,2^{m}) = 2(m+n+1)+1;\\
& \mathbb{K}(a,(4,1)^{n},3,2^{m}) = (2a-1)(m+n+1)+a;\\
& \mathbb{K}((4,1)^{n},3,2^{m},c) = 2(c-1)(m+n+1)+c+1;\\
& \mathbb{K}(a,(4,1)^{n},3,2^{m},c) = (2a-1)(c-1)(m+n+1)+ac+a-1;
\end{cases}
\end{equation}

\begin{equation}\label{lab_56}
\begin{cases}
& \mathbb{K}((1,4)^{n},1,3,2^{m}) = m+n+2;\\
& \mathbb{K}((1,4)^{n},1,3,2^{m},b) = (b-1)(m+n+1)+b;\\
& \mathbb{K}(a,(1,4)^{n},1,3) = (a-2)n+2a-3;\\
& \mathbb{K}(a,(1,4)^{n},1,3,2^{m},b) = (a-2)(b-1)(m+n+1)+ab-b-1.
\end{cases}
\end{equation}

\textsc{Proof.} All these identities follow from the previous lemma and formulas (\ref{lab_06})-(\ref{lab_11}). For example, Lemma 19 together with (\ref{lab_09}) imply
\[
\mathbb{K}(a,2^{n}) = \mathbb{K}(2^{n},a) = a\mathbb{K}(2^{n}) - \mathbb{K}(2^{n-1}) = aD_{n}-D_{n-1} = (a-1)n+a.
\]
Thus, using (\ref{lab_07}), we find
\[
\mathbb{K}(a,2^{n},b) = b\mathbb{K}(a,2^{n}) - \mathbb{K}(a,2^{n-1}) = b((a-1)n+a)-(a-1)(n-1)-a = (a-1)(b-1)n+ab-1.
\]
Further,
\[
\mathbb{K}(2^{n},b,c) = \mathbb{K}(b,c)\mathbb{K}(2^{n}) - \mathbb{K}(c)\mathbb{K}(2^{n-1}) = (bc-1)(n+1)-cn = (bc-c-1)n+bc-1
\]
and finally
\begin{multline*}
\mathbb{K}(a,2^{n},b,c) = \mathbb{K}(b,c)\mathbb{K}(a,2^{n}) - \mathbb{K}(c)\mathbb{K}(a,2^{n-1}) = (bc-1)((a-1)n+a) - c(a-1)(n-1) - ac =\\
=\ (a-1)(bc-c-1)n+abc-a-c.
\end{multline*}
This concludes the proof of (\ref{lab_52}). The rest relations are established by the similar way. $\square$
\vspace{0.3cm}

Lemmas 19, 20 allow one to estimate the remainder in the formula for $N(Q;r,D,c_{0})$. Such estimate is based on the following auxilliary assertion.
\vspace{0.3cm}

\textsc{Lemma 21.} \emph{Suppose that $r\geqslant 3$, $1\leqslant j\leqslant r$, and suppose that the tuple $\mathbf{k}_{j} = (k_{1},\ldots,k_{r})$ has the form:}
\[
k_{j}\geqslant 4r+2,\quad k_{i}=1\quad \textit{for}\quad |i-j|=1 \quad\textit{and}\quad k_{i}=2\quad \textit{for}\quad |i-j|\geqslant 2.
\]
\emph{Further, let $0\leqslant c_{0},c_{1}\leqslant D-1$, $c_{0}\ne c_{1}$ and} $\text{GCD}(c_{0},c_{1},D)=1$. \emph{If $\mathbf{k}_{j}$ belongs to the set $\mathcal{A}_{r}^{\circ}(D,c_{0},c_{1})$ then} (a) $r\leqslant 2D$ \emph{and} (b) $k_{j}\equiv h\pmod{d}$ \emph{where $d = D$ or $d = D/2$ and $0\leqslant h< d$}.
\vspace{0.3cm}

\textsc{Proof.} Suppose that $\mu,\nu\geqslant 0$, $k\geqslant 1$. Then, from Lemma 20 and from (\ref{lab_10}) it follows that
\begin{multline}\label{lab_57}
\mathbb{K}_{\mu+2}(2^{\mu},1,k) = k-\mu-1,\quad \mathbb{K}_{\mu+3}(2^{\mu},1,k,1) = k-\mu-2,\\
\mathbb{K}_{\mu+\nu+3}(2^{\mu},1,k,1,2^{\nu}) = k-(\mu+\nu+2).
\end{multline}
Next, suppose that $a, B, \kappa$ and $m$ satisfy to the following conditions:
\begin{equation}\label{lab_58}
a\not\equiv 0\pmod{B},\quad \kappa\not\equiv ja\pmod{B}\quad\text{for}\quad j = 0,1,\ldots, m-1,\quad \kappa\equiv ma\pmod{B}.
\end{equation}
Then
\begin{equation}\label{lab_59}
m\leqslant \frac{B}{\text{GCD}(a,B)}-1.
\end{equation}
Indeed, (\ref{lab_58}) implies that $ja\not\equiv 0\pmod{M}$ for $1\leqslant j\leqslant m$. Since the period of $ja\pmod{B}$ is equal to $B/\text{GCD}(a,B)$, we get (\ref{lab_59}).

Now we write the precise form of the conditions (\ref{lab_13}) for any tuple $\mathbf{k}_{j}$, $j = 1,2,\ldots, r$.
\vspace{0.3cm}

(a) \textsc{Case} $j=1$. Using (\ref{lab_57}), we find that
\[
c_{1}\mathbb{K}_{i}(k_{1},1,2^{i-2})-c_{0}\mathbb{K}_{i-1}(1,2^{i-2}) = c_{1}(k_{1}-i+1)-c_{0},\quad i=2,3,\ldots,r.
\]
Thus (\ref{lab_13}) takes the form
\[
c_{1}k_{1}-2c_{0}\not\equiv ic_{1}\pmod{D},\quad 0\leqslant i\leqslant r-2,\quad c_{1}k_{1}-2c_{0}\equiv (r-1)c_{1}\pmod{D}.
\]
Taking $a = c_{1}$, $B = D$, $\kappa = c_{1}k_{1}-2c_{0}$ and $m=r-1$, in (\ref{lab_58}), we conclude from (\ref{lab_59}) that
\[
r-1\leqslant \frac{D}{\Delta_{0}}-1,\quad \Delta_{0} = \text{GCD}(c_{1},D),\quad \text{and therefore}\quad r\leqslant \frac{D}{\Delta_{0}}.
\]
Moreover, if the last congruence is solvable, then $2c_{0}\equiv 0\pmod{\Delta_{0}}$. The condition $\text{GCD}(c_{0},c_{1},D)=1$ implies that $\Delta_{0}\,|\,2$ and hence $\Delta_{0} =1$ or $\Delta_{0} = 2$. In the first case we get $k_{1}\equiv h_{1}\pmod{D}$, and in the second one we have $k_{1}\equiv h_{2}\pmod{D/2}$ (for some $h_{1}, h_{2}$).
\vspace{0.3cm}

(b) \textsc{Case} $j=2$. In this case the system (\ref{lab_13}) has the form
\begin{equation*}
\begin{cases}
c_{1}\not\equiv 2c_{0}\pmod{D},\\
(c_{1}-c_{0})(k_{2}-1)-2c_{0}\not\equiv i(c_{1}-c_{0})\pmod{D},\quad i = 0,1,\ldots, r-3,\\
(c_{1}-c_{0})(k_{2}-1)-2c_{0}\equiv (r-2)(c_{1}-c_{0})\pmod{D}.
\end{cases}
\end{equation*}
Hence, setting $a=c_{1}-c_{0}$, $B=D$, $\kappa = (c_{1}-c_{0})(k_{2}-1)-2c_{0}$ and $m=r-2$, by (\ref{lab_59}), we get
\[
r-2\leqslant \frac{D}{\Delta_{1}}-1,\quad \Delta_{1} = \text{GCD}(c_{1}-c_{0},D),\quad \text{and therefore}\quad r\leqslant \frac{D}{\Delta_{1}}+1\leqslant D+1.
\]
The last congruence imply that $\Delta_{1}\,|\,2c_{0}$, so we have $\Delta_{1} = 1$ or $\Delta_{1} = 2$. Hence, $k$ belongs to some progression modulo $D$ or $D/2$.
\vspace{0.3cm}

(c) \textsc{Case} $3\leqslant j\leqslant r-1$. In this case, the conditions (\ref{lab_13}) corresponding  to $1\leqslant i\leqslant j-2$ yield
\[
c_{1}(i+1)-c_{0}i\not\equiv c_{0}\pmod{D},\quad\text{that is,}\quad (c_{1}-c_{0})(i+1)\not \equiv 0\pmod{D}.
\]
In view of (\ref{lab_58}), (\ref{lab_59}), this means that
\begin{equation}\label{lab_60}
j-2\leqslant \frac{D}{\Delta_{1}}-1,\quad \text{and}\quad j\leqslant \frac{D}{\Delta_{1}}+1.
\end{equation}
Next, the conditions (\ref{lab_13}) corresponding to $j\leqslant i\leqslant r$ takes the form
\begin{equation*}
\begin{cases}
i(c_{1}-c_{0})\not\equiv (k_{j}+1)(c_{1}-c_{0})-2c_{0}\pmod{D},\quad i = j,j+1,\ldots,r-1,\\
r(c_{1}-c_{0})\equiv (k_{j}+1)(c_{1}-c_{0})-2c_{0}\pmod{D}.
\end{cases}
\end{equation*}
Taking $a=c_{1}-c_{0}$, $B = Q$, $\kappa = (c_{1}-c_{0})(k_{j}-j+1)-2c_{0}$ and $m=r-j$ in (\ref{lab_58}), we get
\begin{equation}\label{lab_61}
r-j\leqslant \frac{D}{\Delta_{1}}-1.
\end{equation}
Summing (\ref{lab_60}), (\ref{lab_61}) we find
\[
r\leqslant \frac{2D}{\Delta_{1}}.
\]
If the last congruence is solvable then we conclude that $2c_{0}\equiv 0\pmod{\Delta_{1}}$ and therefore $1\leqslant \Delta_{1}\leqslant 2$. Hence, $k$ belongs to some progression modulo $D$ or $D/2$.

(d) \textsc{Case} $j=r$. Here the conditions (\ref{lab_13}) corresponding to $1\leqslant i\leqslant r-2$ yield
\[
(c_{1}-c_{0})(i+1)\not\equiv 0\pmod{D}
\]
and therefore
\[
r-2\leqslant \frac{D}{\Delta_{1}}-1,\quad\text{so we have}\quad r\leqslant \frac{D}{\Delta_{1}}+1.
\]
The congruence (\ref{lab_13}) corresponding to $i=r$ takes the form
\[
(c_{1}-c_{0})k_{r}\equiv (c_{1}-c_{0})(r-1)+2c_{0}\pmod{D},
\]
and we conclude that $k_{r}$ belongs to some progression modulo $D$ or $D/2$.

In all above cases we obviously have $r\leqslant D$. Lemma is proved. $\square$
\vspace{0.3cm}

\textsc{Corollary.} \emph{Under the assumption of Lemma 18, one has}
\[
N(Q;r,D,c_{0}) = \mathfrak{c}_{r}\,\frac{\kappa_{D}}{D^{2}}\,Q^{2} + O\biggl(Q\Bigl(1+\frac{\ln{Q}}{D}\Bigr)\Sigma_{0} + \Bigl(D+r\ln{Q}+\frac{r}{D^{2}}(\ln{Q})^{2}\Bigr)\mathbbm{1}_{\{r\leqslant 2D\}}\biggr).
\]

\textsc{Proof.} Let us fix $j$, $1\leqslant j\leqslant r$. Then Lemma 14 implies that the summation in the sums $U_{j}, V_{j}$ is taken over the tuples $\textbf{k}_{j}$ (in the notations of Lemma 21). By Lemma 21, if $r>2D$ then this sum is empty. Suppose now that $r\leqslant 2D$. Then Lemma 21 implies that $k_{j} = \|\textbf{k}_{j}\|$ belongs to some progression modulo $d$, where $d\geqslant D/2$, that is, $k_{j} = h + dm$, $m\geqslant 0$. Since $k_{j}\geqslant 4r+2$ then
\[
U_{j}\leqslant \frac{1}{4r+2}+\sum\limits_{1\leqslant m\leqslant 2Q/d}\frac{1}{md}\ll \frac{1}{r} + \frac{1}{D}\ln{\frac{2Q}{d}}\ll \frac{1}{r}+\frac{\ln{Q}}{D}.
\]
This means that
\[
\Sigma_{1}\ll rD\biggl(\frac{1}{r}+\frac{\ln{Q}}{D}\biggr)\mathbbm{1}_{\{r\leqslant 2D\}}\ll \bigl(D+r\ln{Q}\bigr)\mathbbm{1}_{\{r\leqslant 2D\}}.
\]
Similarly, if $r\leqslant 2D$ then
\[
V_{j}\leqslant \frac{1}{(2Q)^{3}} + \sum\limits_{m>2Q/d}\frac{1}{(md)^{3}}\ll \frac{1}{Q^{3}} + \frac{1}{DQ^{2}}\ll \frac{1}{DQ^{2}}
\]
and therefore
\[
\frac{Q^{2}}{D^{2}}\Sigma_{2}\ll \frac{\mathbbm{1}_{\{r\leqslant 2D\}}}{D^{3}}\ll 1.
\]
Now it remains to note that
\[
\biggl(1+\frac{\ln{Q}}{D}\biggr)\Sigma_{1}\ll \Bigl(D+r\ln{Q}+\frac{r}{D^{2}}(\ln{Q})^{2}\Bigr)\mathbbm{1}_{\{r\leqslant 2D\}}. \square
\]

To establish the desired formula for the proportion $\nu(Q;r;D,c_{0})$, we need the following lemma concerning the summation of Euler function over the progression.
\vspace{0.3cm}

\textsc{Lemma 22.} \emph{Suppose that} $0\leqslant c\leqslant D-1$, $\Delta = \text{GCD}(c,D)$. \emph{Then, uniformly in $c$ and $D\leqslant Q$, one has}
\[
N(Q;D,c) = \sum\limits_{\substack{q\equiv c\;(\mmod{D}) \\ q\leqslant Q}}\varphi(q) = \frac{\kappa_{D}}{2D}\frac{\varphi(\Delta)}{\Delta}\,Q^{2} + O(Q\ln{Q}),
\]
\emph{where the implied constant is absolute.}
\vspace{0.3cm}

For the proof, see \cite[\S 4.2, pp.~190-192]{Postnikov_1971}.
\vspace{0.3cm}

\textsc{Lemma 23.} \emph{Uniformly in $Q$, $1\leqslant r\leqslant Q/3$, $D\geqslant 2$, $D = o(\sqrt{Q\mathstrut})$, $0\leqslant c_{0}\leqslant D-1$ one has}
\begin{multline*}
\nu(Q;r,D,c_{0}) = \frac{2}{D}\frac{\Delta}{\varphi(\Delta)}\bigl(\mathfrak{c}_{r} + O(R)\bigr),\quad\textit{where}\\
\quad \Delta = \text{GCD}(c_{0},D),\quad \mathfrak{c}_{r} = \sum\limits_{\substack{0\leqslant c_{1}\leqslant D-1 \\ c_{1}\ne c_{0} \\ \text{GCD}(c_{1},\Delta)=1}}\sum\limits_{\mathbf{k}\in\,\mathcal{A}_{r}^{\circ}(D,c_{0},c_{1})}|\mathcal{T}(\mathbf{k})|, \\
R = \frac{D^{2}}{Q}\biggl(1+\frac{\ln{Q}}{D}\biggr)\Sigma_{0} + \bigl(D^{3}+Dr\ln{Q}+r(\ln{Q})^{2}\bigr)\frac{\mathbbm{1}_{\{r\leqslant 2D\}}}{Q^{2}}+\mathfrak{c}_{r}\,\frac{D\Delta}{\varphi(\Delta)}\,\frac{\ln{Q}}{Q}\\
\textit{and}\quad \Sigma_{0} = \sum\limits_{c_{1}}\sum\limits_{\substack{\mathbf{k}\in\, \mathcal{A}_{r}^{\circ}(D,c_{0},c_{1}) \\ \|\mathbf{k}\|\leqslant 4r+1}}1.
\end{multline*}
\emph{In particular, if $D$ is a fixed constant then}
\[
R\ll_{D} (\mathfrak{c}_{r}+\Sigma_{0})\,\frac{\ln{Q}}{Q}.
\]

\textsc{Proof.} Rewriting the formulas of Lemma 22 and Corollary of Lemma 21 in the form
\[
N(Q;D,c_{0}) = \frac{\kappa_{D}}{2D}\frac{\varphi(\Delta)}{\Delta}\,Q^{2}\bigl(1+O(\delta_{1})\bigr),\quad
N(Q;r,D,c_{0}) = \frac{\kappa_{D}}{D^{2}}\,Q^{2}\bigl(\mathfrak{c}_{r}+O(\delta_{2})\bigr)
\]
where
\[
\delta_{1} = \frac{D\Delta}{\varphi(\Delta)}\,\frac{\ln{Q}}{Q},\quad \delta_{2} =
\frac{D^{2}}{Q}\biggl(1+\frac{\ln{Q}}{D}\biggr)\Sigma_{0} + \bigl(D^{3}+Dr\ln{Q}+r(\ln{Q})^{2}\bigr)\frac{\mathbbm{1}_{\{r\leqslant 2D\}}}{Q^{2}}
\]
we easily obtain that
\[
\nu(Q;r,D,c_{0}) = \frac{2}{D}\frac{\Delta}{\varphi(\Delta)}\biggl(\mathfrak{c}_{r} + O(\mathfrak{c}_{r}\delta_{1}+\delta_{2})\biggr).\square
\]

\section{Theorem 1: case $\boldsymbol{1\leqslant r\leqslant 7}$}

In \S\S 6-8, we assume that $D = 3$, $c_{0}=0$. In this case, for any $c_{1}$, $1\leqslant c_{1}\leqslant 2$, the conditions (\ref{lab_21}) take the form
\begin{equation}\label{lab_62}
k_{1}\equiv 0\pmod{3}
\end{equation}
if $r=1$, and
\begin{equation}\label{lab_63}
\begin{cases}
\mathbb{K}_{i}(k_{1},\ldots,k_{i})\not\equiv 0 \pmod{3},\\
\mathbb{K}_{r}(k_{1},\ldots,k_{r})\equiv 0 \pmod{3},
\end{cases}
\end{equation}
if $r\geqslant 2$. This means, in particular, that the sets $\mathcal{A}_{r}^{\circ}(3;0,1)$ and $\mathcal{A}_{r}^{\circ}(3;0,2)$ coincide and therefore
\[
\mathfrak{c}_{r} = \sum\limits_{c_{1}=1}^{2}\sum\limits_{\mathbf{k}\in\, \mathcal{A}_{r}^{\circ}(3,0,c_{1})}|\mathcal{T}(\mathbf{k})| = 2\sum\limits_{\mathbf{k}\in\, \mathcal{A}_{r}^{\circ}}|\mathcal{T}(\mathbf{k})|
\]
where $\mathcal{A}_{r}^{\circ} = \mathcal{A}_{r}^{\circ}(3,0,1)$.

Next, if $1\leqslant r\leqslant 7$ then the sum $\Sigma_{0}$,
\[
\Sigma_{0} = \sum\limits_{c_{1}=1}^{2}\sum\limits_{\substack{\mathbf{k}\in \mathcal{A}_{r}^{\circ}(3;0,c_{1}) \\ \|\mathbf{k}\|\leqslant 4r+1}}1 = 2\sum\limits_{\substack{\mathbf{k}\in \mathcal{A}_{r}^{\circ} \\ \|\mathbf{k}\|\leqslant 4r+1}}1
\]
is bounded by some absolute constant (for example, one has $\Sigma_{0}<(4r+1)^{r}\leqslant 2\cdot 29^{7}$). Hence, the formula of Lemma 23 implies that
\[
\nu(Q;r,3,0) = \mathfrak{c}_{r} + O\biggl(\frac{\ln{Q}}{Q}\biggr).
\]
and thus
\[
\nu(r) = \nu(r,3,0) = \lim_{Q\to +\infty}\nu(Q;r;3,0) = \mathfrak{c}_{r}.
\]
The problem is to give a precise description of the sets $\mathcal{A}_{r}^{\circ}$.
\vspace{0.3cm}

\textsc{Case} $r=1$. By (\ref{lab_62}), $\mathcal{A}_{r}^{\circ}$ consists of the tuples $(k_{1})$ of unit length such that $\mathbb{K}_{1}(k_{1}) = k_{1}\equiv 0\pmod{3}$. Thus Lemma 7 implies
\[
\nu(1) = 2\sum\limits_{\substack{k_{1}\equiv 0\;(\mmod{3}) \\ k_{1} \geqslant 1}}|\mathcal{T}(k_{1})| = 8\Sigma,\quad \Sigma = \sum\limits_{k=1}^{+\infty}\frac{1}{3k(3k+1)(3k+2)}.
\]
On can easily check that
\begin{multline*}
\Sigma = \frac{1}{2}\sum\limits_{k=1}^{+\infty}\biggl(\frac{1}{3k}-\frac{1}{3k+2}\biggr) -\sum\limits_{k=1}^{+\infty}\biggl(\frac{1}{3k+1}-\frac{1}{3k+2}\biggr) = \\
= \frac{1}{6}\sum\limits_{k=1}^{+\infty}\biggl(\frac{1}{k}-\frac{1}{k+\tfrac{2}{3}}\biggr) - \frac{1}{3}\sum\limits_{k=1}^{+\infty}\biggl(\frac{1}{k+\tfrac{1}{3}}-\frac{1}{k+\tfrac{2}{3}}\biggr) = \\
= \frac{1}{6}\biggl(\psi\biggl(\frac{5}{3}\biggr) - \psi(1)\biggr) - \frac{1}{3}\biggl(\psi\biggl(\frac{2}{3}\biggr) - \psi\biggl(\frac{1}{3}\biggr)\biggr),
\end{multline*}
where the symbol $\psi(x)$ stands for digamma function. Using the identities
\[
\psi(x+1)-\psi(x) = \frac{1}{x},\quad \psi(1-x)-\psi(x) = \pi\ctg{(\pi x)}
\]
together with well-knows formulas
\[
\psi(1) = -\,\gamma,\quad \psi\biggl(\frac{2}{3}\biggr) = -\,\gamma + \frac{\pi}{2\sqrt{3}}-\frac{3}{2}\ln{3}
\]
(here $\gamma$ is Euler constant; see formulas 8.365.1, 8.365.8, 8.366.1 and 8.366.7 in \cite{Ryjik_Gradstein_1962}), we find
\[
\Sigma = \frac{3}{4}-\frac{1}{4}\biggl(\frac{\pi}{\sqrt{3}}+\ln{3}\biggr),\quad \nu(1) = 6-2\biggl(\frac{\pi}{\sqrt{3}}+\ln{3}\biggr).
\]

In the \textsc{Case} $r=2$, (\ref{lab_63}) takes the form
\[
\mathbb{K}_{1}(k_{1}) = k_{1}\not\equiv 0\pmod{3},\quad \mathbb{K}_{2}(k_{1},k_{2}) = k_{1}k_{2}-1\equiv 0\pmod{3}.
\]
Therefore, $\mathcal{A}_{2}^{\circ}$ contains only the tuples with $k_{1}\equiv k_{2}\not\equiv 0\pmod{3}$. Hence, $\nu(2) = \sigma_{1}+\sigma_{2}$, where
\[
\sigma_{j} = \sum\limits_{\substack{k_{1}\equiv k_{2}\equiv j\;(\mmod{3})\\ k_{1},k_{2} \geqslant 1}}|\mathcal{T}(k_{1},k_{2})|,\quad j = 1,2.
\]
In view of Lemma 9, the contribution from $k_{1} = 1$ to the sum $\sigma_{1}$ is equal to
\[
\sigma_{11} = \sum\limits_{\substack{k_{2}\equiv 1\;(\mmod{3}) \\ k_{2} \geqslant 4}}|\mathcal{T}(1,k_{2})| = |\mathcal{T}(1,4)| + \sum\limits_{\substack{k_{2}\equiv 1\;(\mmod{3}) \\ k_{2} \geqslant 7}}|\mathcal{T}(k_{2})|.
\]
Since $|\mathcal{T}(1,4)| = 1/35$ then Lemma 7 implies
\[
\sigma_{11} = \frac{1}{35} + 4\biggl(\sigma - \frac{1}{4\cdot 5\cdot 6}\biggr) = 4\sigma - \frac{1}{210},\quad \sigma = \sum\limits_{k =1}^{+\infty}\frac{1}{(3k+1)(3k+2)(3k+3)}.
\]
Further, by Lemma 10, the contribution coming from $k_{1} = 4$ to the sum $\sigma_{1}$ equals to
\[
\sigma_{12} = \sum\limits_{\substack{k_{2}\equiv 1\;(\mmod{3}) \\ k_{2} \geqslant 1}}|\mathcal{T}(4,k_{2})| = |\mathcal{T}(4,1)| = \frac{1}{35},
\]
and, by Lemma 12, the contribution of $k_{1}\geqslant 7$, $k_{1}\equiv 1\pmod{3}$ coincides with
\begin{multline*}
\sigma_{13} = \sum\limits_{\substack{k_{1}\equiv 1\;(\mmod{3}) \\ k_{1} \geqslant 7}}\sum\limits_{\substack {k_{2}\equiv 1\;(\mmod{3}) \\ k_{2}\geqslant 1}}|\mathcal{T}(k_{1},k_{2})| =
\sum\limits_{\substack{k_{1}\equiv 1\;(\mmod{3}) \\ k_{1} \geqslant 7}}|\mathcal{T}(k_{1},1)| = \\
=\sum\limits_{\substack{k_{1}\equiv 1\;(\mmod{3}) \\ k_{1} \geqslant 7}}|\mathcal{T}(k_{1})| = 4\biggl(\sigma - \frac{1}{4\cdot 5\cdot 6}\biggr) = 4\sigma - \frac{1}{30}.
\end{multline*}
Hence,
\[
\sigma_{1} = \sigma_{11}+\sigma_{12}+\sigma_{13} = \biggl(4\sigma-\frac{1}{210}\biggr) + \frac{1}{35} + \biggl(4\sigma - \frac{1}{30}\biggr) = 8\sigma - \frac{1}{105}.
\]
Finally, Lemma 10 implies that the sum $\sigma_{2}$ contains only one non-zero term, corresponding to $k_{1} = k_{2} = 2$, and therefore
\[
\sigma_{2} = |\mathcal{T}(2,2)| = \frac{1}{10}.
\]
Thus, the set $\mathcal{A}_{2}^{\circ}$ consists of the tuples
\begin{align*}
& (k_{1},1),\quad k_{1}\equiv 1\pmod{3},\quad k_{1}\geqslant 4, \\
& (1,k_{2}),\quad k_{2}\equiv 1\pmod{3},\quad k_{2}\geqslant 4, \\
& (2,2),
\end{align*}
and we have
\[
\nu(2) = 2\biggl(8\sigma - \frac{1}{105}+\frac{1}{10}\biggr) = 16\sigma + \frac{19}{105}.
\]
Since
\[
\sigma = \frac{1}{4}\biggl(\frac{\pi}{\sqrt{3}} - \ln{3}\biggr)-\frac{1}{6},
\]
we have
\[
\nu(2) = 4\biggl(\frac{\pi}{\sqrt{3}} - \ln{3}\biggr)-\frac{87}{35}.
\]

The calculation of proportions $\nu(r)$, $3\leqslant r\leqslant 6$, will follow the following scheme.

For any $r$, $1\leqslant r\leqslant 5$, we determine the set $\mathcal{A}_{r}^{*}$ of tuples $\mathbf{k}_{r} = (k_{1},\ldots,k_{r})$ that correspond to non-empty regions $\mathcal{T}(\mathbf{k}_{r})$ and obey an extra conditions $\mathbb{K}_{i}(k_{1},\ldots,k_{i})\not\equiv 0\pmod{3}$, $i = 1,2,\ldots,r$.

Obviously, $\mathcal{A}_{1}^{*}$ consists of the tuples $(k_{1})$ of the length one such that $k_{1}\not\equiv 0\pmod{3}$.

\textsc{Case} $r=2$. Using the above description of the set $\mathcal{A}_{2}^{\circ}$ together with Lemmas 9, 19, we conclude that the set $\mathcal{A}_{2}^{*}$ is formed by the tuples
\begin{align*}
& (k_{1},1), \quad k_{1}\equiv 2\pmod{3},\quad k_{1}\geqslant 2; \\
& (1,k_{2}), \quad k_{2}\not\equiv 1\pmod{3},\quad k_{2}\geqslant 2; \\
& (2,3),\;(2,4),\;(4,2).
\end{align*}

Given the set $\mathcal{A}_{r-1}^{*}$, we take any tuple $\mathbf{k}_{r-1} = (k_{1},\ldots,k_{r-1})\in \mathcal{A}_{r-1}^{*}$ and then check all the tuples $\mathbf{k}_{r} = (\mathbf{k}_{r-1},k) = (k_{1},\ldots,k_{r-1},k)$, $k = 1,2,3,\ldots$ whether they satisfy the above conditions. We also note that Lemma 10 allows one to shorten the search.
Since the conditions $\mathbb{K}_{i}(k_{1},\ldots,k_{i})\not\equiv 0\pmod{3}$ are satisfied automatically for any $1\leqslant i\leqslant r-1$, it is sufficient to check whether the determinant
$\mathbb{K}_{r}(k_{1},\ldots,k_{r-1},k) = \mathbb{K}_{r}$ is divisible by $3$.

If $\mathbb{K}_{r}\not\equiv 0\pmod{3}$ for $k$ under considering and $\mathcal{T}(\mathbf{k}_{r})\ne \varnothing$ then we put the tuple $\mathbf{k}_{r}$ into the set $\mathcal{A}_{r}^{*}$.
Otherwise, such a tuple satisfies to (\ref{lab_13}). In the last case, it is included into the set $\mathcal{A}_{r}^{\circ}$ and should be taken into account in the calculation of $\nu(r;3)$.

The search of the tuples and the checking of non-emptiness of the regions $\mathcal{T}(\mathbf{k}_{r})$, $\mathbf{k}_{r} = (k_{1},\ldots,k_{r})$ (that is, the calculation of its area) uses the package Wolfram Mathematica. In particular, it allow one to find the precise form of the functions
\[
f_{i}(x) = f_{i}(x;\mathbf{k}_{r}),\quad g_{i}(x) = g_{i}(x;\mathbf{k}_{r}),\quad i = 1,2,\ldots
\]
from Lemma 2 and the functions
\begin{multline*}
F_{r}(x) = F_{r}(x;\mathbf{k}_{r}) = \max{(1-x,f_{1}(x),\ldots,f_{r}(x))},\\
G_{r}(x) = G_{r}(x;\mathbf{k}_{r}) = \min{(1,g_{1}(x),\ldots,g_{r}(x))}.
\end{multline*}
The area of $\mathcal{T}(\mathbf{k}_{r})$ is calculated by the integration of the function
\[
H(x) = H(x;\mathbf{k}_{r}) = \max{(0,G_{r}(x)-F_{r}(x))}.
\]
Such search can be shorten by application of Lemma 10. In what follows, we give only the final results of such calculations.

\textsc{Case} $r = 3$. The set $\mathcal{A}_{3}^{*}$ consists of the tuples
\begin{align*}
& (2,1,k_{3}), \quad k_{3}\not\equiv 2\pmod{3},\quad k_{3}\geqslant 6; \\
& (1,k_{2},1), \quad k_{2}\equiv 0\pmod{3},\quad k_{2}\geqslant 3; \\
& (1,2,2),\;(1,2,3),\;(2,3,2),\;(2,4,1),\;(5,1,3),\;(5,1,4),\;(8,1,3),
\end{align*}
while the set $\mathcal{A}_{3}^{\circ}$ consists of the tuples
\begin{align*}
& (k_{1},1,2), \quad k_{1}\equiv 2\pmod{3},\quad k_{1}\geqslant 8; \\
& (1,k_{2},1), \quad k_{2}\equiv 2\pmod{3},\quad k_{2}\geqslant 5; \\
& (2,1,k_{3}), \quad k_{3}\equiv 2\pmod{3},\quad k_{3}\geqslant 8; \\
& (1,2,4),\;(1,3,2),\;(2,3,1),\;(4,2,1).
\end{align*}
Therefore, the formula for $\nu(3)$ involves the values
\begin{align*}
& |\mathcal{T}(1,2,4)| = |\mathcal{T}(4,2,1)| = \frac{1}{210};\quad |\mathcal{T}(1,3,2)| = |\mathcal{T}(2,3,1)| = \frac{3}{140}, \\[6pt]
& |\mathcal{T}(k,1,2)| = |\mathcal{T}(2,1,k)| =
\begin{cases}
  11/2\,340, \; k = 8, \\[6pt]
  |\mathcal{T}(k)|, \;k\equiv 2\pmod{3},\;k\geqslant 11,
\end{cases} \\
& |\mathcal{T}(1,k,1)| =  |\mathcal{T}(k)|,\quad k\equiv 2\pmod{3},\quad k\geqslant  5.
\end{align*}
Hence,
\begin{multline*}
\nu(3) = 2\biggl(\frac{2}{210} + \frac{2\cdot 3}{140} + \frac{2\cdot 11}{2\,340} + 2\sum\limits_{\substack{k\equiv 2\;(\mmod{3}) \\ k\geqslant 11}}|\mathcal{T}(k)| +
\sum\limits_{\substack{k\equiv 2\;(\mmod{3}) \\ k\geqslant 5}}|\mathcal{T}(k)|\biggr) = \\
= \frac{103}{4\,095} + 24\sum\limits_{k=1}^{+\infty}\frac{1}{(3k+2)(3k+3)(3k+4)} = 12\ln{3} - \frac{53\,132}{4\,095}.
\end{multline*}

\textsc{Case} $r=4$. The set $\mathcal{A}_{4}^{*}$ consists of the tuples
\begin{align*}
& (1,2,2,2) && (1,3,1,7) && (2,4,1,4) \\
& (1,2,2,3) && (1,6,1,3) && (5,1,3,2) \\
& (1,2,3,1) && (2,3,2,2) && (5,1,4,1) \\
& (1,3,1,6) && (2,4,1,3) && (2,1,k_{3},1), && k_{3}\equiv 1\pmod{3},\;k_{3}\geqslant 7,
\end{align*}
and the set $\mathcal{A}_{4}^{\circ}$ consists of the tuples
\begin{align*}
& (1,k_{2},1,2), && k_{2}\equiv 0\pmod{3},\;k_{2}\geqslant 6, \\
& (2,1,k_{3},1), && k_{3}\equiv 0\pmod{3},\;k_{3}\geqslant 6, \\
& (1,2,3,2) && (1,3,1,5) && (1,3,1,8) \\
& (2,3,2,1) && (5,1,3,1) && (8,1,3,1).
\end{align*}
Consequently, the expression for $\nu(4)$ involves the areas
\begin{align*}
& |\mathcal{T}(1,k,1,2)| = |\mathcal{T}(2,1,k,1)| =
\begin{cases}
  1/462, \; k = 6, \\[6pt]
  |\mathcal{T}(k)|, \;k\equiv 0\pmod{3},\;k\geqslant 9,
\end{cases} \\
& |\mathcal{T}(1,2,3,2)| = |\mathcal{T}(2,3,2,1)| = \frac{9}{1\,540},\\[6pt]
& |\mathcal{T}(1,3,1,8)| = |\mathcal{T}(8,1,3,1)| = \frac{1}{1\,170}, \\[6pt]
& |\mathcal{T}(1,3,1,5)| =  |\mathcal{T}(5,1,3,1)|= \frac{1}{195}.
\end{align*}
Therefore,
\[
\nu(4) = 2\biggl(\frac{2}{195}+\frac{2}{1\,170}+\frac{2\cdot 9}{1\,540} + \frac{2}{462} + 2\sum\limits_{\substack{k\equiv 0\;(\mmod{3}) \\ k\geqslant 9}}|\mathcal{T}(k)|\biggr) =
\frac{528\,904}{45\,045} - 4\biggl(\frac{\pi}{\sqrt{3}}+\ln{3}\biggr).
\]
\textsc{Case} $r = 5$. The set $\mathcal{A}_{5}^{*}$ consists of the tuples
\begin{align*}
& (1,2,2,2,2) && (1,3,1,7,1) && (5,1,4,1,3) \\
& (1,2,2,2,3) && (1,6,1,3,2) && (5,1,4,1,4) \\
& (1,2,2,3,1) && (2,4,1,3,2) && (2,1,7,1,3) \\
& (1,2,3,1,4) && (2,4,1,4,1) \\
& (1,2,3,1,6) && (5,1,3,2,2),
\end{align*}
and the set $\mathcal{A}_{5}^{\circ}$ consists of the tuples
\begin{align*}
& (2,1,k_{3},1,2),&& k_{3}\equiv 1\pmod{3},\;k_{3}\geqslant 7,  \\
& (1,2,2,3,2) && (2,3,2,2,1)  \\
& (1,2,3,1,5) && (5,1,3,2,1) \\
& (1,3,1,6,1) && (1,6,1,3,1).
\end{align*}
Thus the formula for $\nu(5)$ involves the areas
\begin{align*}
& |\mathcal{T}(1,2,2,3,2)| = |\mathcal{T}(2,3,2,2,1)| = \frac{1}{770}, \quad |\mathcal{T}(1,2,3,1,5)| = |\mathcal{T}(5,1,3,2,1)| = \frac{11}{2\,210} \\[6pt]
& |\mathcal{T}(1,3,1,6,1)| = |\mathcal{T}(1,6,1,3,1)| = \frac{1}{143}, \quad |\mathcal{T}(2,1,k,1,2)| =
\begin{cases}
1/616, \; k = 7, \\[6pt]
|\mathcal{T}(k)|, \;k\equiv 1\,(\mmod{3}),\;k\geqslant 10,
\end{cases}
\end{align*}
so we have
\[
\nu(5) = 2\biggl(\frac{2}{770}+\frac{2\cdot 11}{2\,210}+\frac{2}{143}+\frac{1}{616}+\sum\limits_{\substack{k\equiv 1\;(\mmod{3}) \\ k\geqslant 10}}|\mathcal{T}(k)|\biggr)
= 2\biggl(\frac{\pi}{\sqrt{3}}-\ln{3}\biggr) - \frac{4\,164\,383}{3\,063\,060}.
\]

\begin{center}
\includegraphics{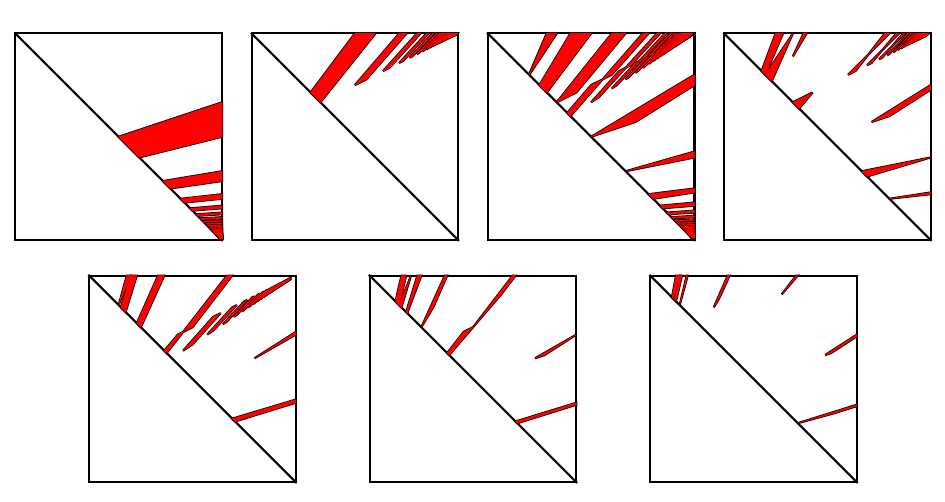}
\fontsize{10}{12pt}\selectfont
\emph{Fig.~2}. Polygons whose areas are involved in the calculation of the proportions $\nu(r)$, $1\leqslant r\leqslant 7$.
\fontsize{12}{15pt}\selectfont
\end{center}

\textsc{Case} $r = 6$. The set $\mathcal{A}_{6}^{*}$ consists of the tuples
\begin{align*}
& (1,2,2,2,2,2) && (1,2,2,3,1,6) && (2,4,1,4,1,3) && (5,1,4,1,4,1)\\
& (1,2,2,2,2,3) && (1,3,1,7,1,3) && (2,4,1,4,1,4)\\
& (1,2,2,2,3,1) && (1,6,1,3,2,2) && (5,1,3,2,2,2) \\
& (1,2,2,3,1,4) && (2,4,1,3,2,2) && (5,1,4,1,3,2)
\end{align*}
while the set $\mathcal{A}_{6}^{\circ}$ is formed by the tuples
\begin{align*}
& (1,2,2,3,1,5) && (5,1,3,2,2,1) && (1,2,3,1,6,1) && (1,6,1,3,2,1) \\
& (1,2,3,1,4,2) && (2,4,1,3,2,1) && (1,3,1,7,1,2) && (2,1,7,1,3,1).
\end{align*}
Hence,
\begin{multline*}
\nu(6) = 2\bigl(|\mathcal{T}(1,2,2,3,1,5)| + |\mathcal{T}(5,1,3,2,2,1)|  + |\mathcal{T}(1,2,3,1,4,2)|  + |\mathcal{T}(2,4,1,3,2,1)| + \\[6pt]
 + |\mathcal{T}(1,2,3,1,6,1)|  + |\mathcal{T}(1,6,1,3,2,1)|  + |\mathcal{T}(1,3,1,7,1,2)|+ |\mathcal{T}(2,1,7,1,3,1)|\bigr) = \\[6pt]
 =4\biggl(\frac{37}{12\,376} + \frac{1}{910} + \frac{23}{9\,724} + \frac{3}{1\,144}\biggr) = \frac{3\,089}{85\,085}.
\end{multline*}
Case $r=7$. The set $\mathcal{A}_{7}^{\circ}$ consists of the tuples
\begin{align*}
& (1,2,2,2,3,1,5) && (5,1,3,2,2,2,1) && (1,2,2,3,1,6,1) && (1,6,1,3,2,2,1) \\
& (1,2,2,3,1,4,2) && (2,4,1,3,2,2,1) && (1,3,1,7,1,3,1), &&
\end{align*}
and therefore
\begin{multline*}
\nu(7) = \\[6pt]
=2\bigl(|\mathcal{T}(1,2,2,2,3,1,5)| + |\mathcal{T}(5,1,3,2,2,2,1)|  + |\mathcal{T}(1,2,2,3,1,4,2)|  + |\mathcal{T}(2,4,1,3,2,2,1)| + \\[6pt]
 + |\mathcal{T}(1,2,2,3,1,6,1)|  + |\mathcal{T}(1,6,1,3,2,2,1)|  + |\mathcal{T}(1,3,1,7,1,3,1)|\bigr) = \\[6pt]
 =2\biggl(\frac{2\cdot 13}{9\,576} + \frac{2\cdot 15}{12\,376} + \frac{2}{2\,618} + \frac{1}{936}\biggr) = \frac{54\,097}{3\,879\,876}.
\end{multline*}
Thus, Theorem 1 is proved for all $1\leqslant r\leqslant 7$. $\square$

\section{Explicit form of some polygons $\boldsymbol{\mathcal{T}(\mathbf{k})}$.}

To calculate the proportions $\nu(r)$, $r\geqslant 8$, we need the precise description of regions $\mathcal{T}(\mathbf{k})$ corresponding to the tuples $\mathbf{k}$ of the forms $(2,(4,1)^{n})$, $(5,(1,4)^{n})$ and $(1,2^{n})$ ($n = 1,2,3,\ldots$). Such description is contained in the following three lemmas.
\vspace{0.3cm}

\textsc{Lemma 24.} \emph{Suppose that $n\geqslant 1$ is integer. Then the polygon $\mathcal{T}(2,(4,1)^{n})$ is defined by the inequalities}
\[
\frac{6n-2}{6n-1}<x\leqslant 1,\quad \frac{1}{3}(1+x)<y\leqslant \frac{1+4nx}{6n+1};
\]
\emph{and is a triangle with vertices}
\[
\biggl(\frac{6n-2}{6n-1},\frac{4n-1}{6n-1}\biggr),\quad \biggl(1,\frac{4n+1}{6n+1}\biggr),\quad \biggl(1,\frac{2}{3}\biggr)\quad\textit{and with area}\quad \frac{1}{6(6n-1)(6n+1)}.
\]
\textsc{Proof} (induction by $n$). Suppose that $n=1$. Then the desired region $\mathcal{T}_{3} = \mathcal{T}(2,4,1)$ is described by the inequalities (here and in what follows, the lower index is equal to the total number of components of the tuple)
\begin{equation}\label{lab_64}
\max{(1-x,f_{1}(x),f_{2}(x),f_{3}(x))}<y\leqslant \min{(1,g_{1}(x),g_{2}(x),g_{3}(x))},
\end{equation}
where
\begin{align*}
& f_{1}(x) = f_{1}(x;2) = \frac{1}{3}(1+x), && g_{1}(x) = g_{1}(x;2) = \frac{1}{2}(1+x),\\
& f_{2}(x) = f_{2}(x;2,4) = \frac{1}{9}(1+5x), && g_{2}(x) = g_{2}(x;2,4) = \frac{1}{7}(1+4x),\\
& f_{3}(x) = f_{3}(x;2,4,1) = \frac{1}{12}(1+7x), && g_{3}(x) = g_{3}(x;2,4,1) = \frac{1}{5}(1+3x).
\end{align*}
Since $f_{2}(x), f_{3}(x)\leqslant f_{1}(x)$ и $g_{2}(x)\leqslant g_{1}(x), g_{3}(x)$ for $0\leqslant x\leqslant 1$, then (\ref{lab_64}) take the form
\[
\max{\Bigl\{1-x,\frac{1}{3}(1+x)\Bigr\}}<y\leqslant \frac{1}{7}(1+4x)
\]
or, that is the same,
\[
\frac{4}{5} < x\leqslant 1,\quad \frac{1}{3}(1+x) < y\leqslant \frac{1}{7}(1+4x).
\]
One can check that such region appears to be the triangle with the vertices $\bigl(\tfrac{4}{5},\tfrac{3}{5}\bigr)$, $\bigl(1,\tfrac{5}{7}\bigr)$, $\bigl(1,\tfrac{2}{3}\bigr)$, whose area equals to $1/210 = 1/(6\cdot 5\cdot 7)$.

Assume now that the assertion is checked for all $n\leqslant m-1$. Then, passing from the polygon и $\mathcal{T}_{2m-1} = \mathcal{T}(2,(4,1)^{m-1})$ describing by the inequalities
\begin{equation}\label{lab_65}
\frac{6m-8}{6m-7}<x\leqslant 1,\quad \frac{1}{3}(1+x)<y\leqslant \frac{1+4(m-1)x}{6m-5},
\end{equation}
to the polygon $\mathcal{T}_{2m+1} = \mathcal{T}(2,(4,1)^{m})$, we should add to the system (\ref{lab_65}) the following conditions:
\begin{equation}\label{lab_66}
\begin{cases}
f_{2m}(x) = f_{2m}(x;2,(4,1)^{m-1},4) < y \leqslant g_{2m}(x) = g_{2m}(x;2,(4,1)^{m-1},4),\\
f_{2m+1}(x) = f_{2m+1}(x;2,(4,1)^{m}) < y \leqslant g_{2m+1}(x) = g_{2m+1}(x;2,(4,1)^{m}).
\end{cases}
\end{equation}
Now Lemmas 19, 20 yield
\begin{align*}
& f_{2m}(x) = \frac{1+x\mathbb{K}((4,1)^{m-1},5)}{\mathbb{K}(2,(4,1)^{m-1},5)} = \frac{1+(6m-1)x}{9m},\\
& g_{2m}(x) = \frac{1+x\mathbb{K}((4,1)^{m-1},4)}{\mathbb{K}(2,(4,1)^{m-1},4)} = \frac{1+4mx}{6m+1},\\
& f_{2m+1}(x) = \frac{1+x\mathbb{K}((4,1)^{m-1},4,2)}{\mathbb{K}(2,(4,1)^{m-1},4,2)} = \frac{1+(6m+1)x}{9m+3},\\
& g_{2m+1}(x) = \frac{1+x\mathbb{K}((4,1)^{m})}{\mathbb{K}(2,(4,1)^{m})} = \frac{1+(2m+1)x}{3m+2}.
\end{align*}
It is easy to see that $f_{2m+1}(x) \leqslant f_{2m}(x)$ and $g_{2m}(x) < g_{2m+1}(x)$ for all $x$, $0\leqslant x\leqslant 1$. Hence, the conditions (\ref{lab_66}) take the form
\begin{equation}\label{lab_67}
\frac{1+(6m-1)x}{9m}<y\leqslant \frac{1+4mx}{6m+1}.
\end{equation}
Therefore, the region $\mathcal{T}_{2m+1}$ is described by the system of inequalities (\ref{lab_65}), (\ref{lab_66}). Since
\[
\frac{1+(6m-1)x}{9m}\leqslant \frac{1}{3}(1+x),\quad \frac{1+4mx}{6m+1}\leqslant \frac{1+4(m-1)x}{6m-5}
\]
for any $x$, $0\leqslant x\leqslant 1$, then both the inequalities (\ref{lab_65}), (\ref{lab_65}) can be replaced by the system
\[
\frac{6m-2}{6m-1}<x\leqslant 1,\quad \frac{1}{3}(1+x)<y\leqslant \frac{1+4mx}{6m+1}.
\]
Thus, $\mathcal{T}_{2m+1}$  is a triangle with vertices in the points of intersection of the lines
\[
y = \frac{1}{3}(1+x)\quad\text{и}\quad y = \frac{1+4mx}{6m+1},\quad y = \frac{1+4mx}{6m+1} \quad\text{и}\quad x=1,\quad y = \frac{1}{3}(1+x) \quad\text{и}\quad x=1,
\]
that is, in the points
\[
\biggl(\frac{6m-2}{6m-1},\frac{4m-1}{6m-1}\biggr),\quad \biggl(1,\frac{4m+1}{6m+1}\biggr),\quad \biggl(1,\frac{2}{3}\biggr).
\]
Further, if $\mathcal{M}$ is a planar $k$-gon with vertices $(x_{i},y_{i})$, $i = 1,2,\ldots,k$, then its area $|\mathcal{M}|$ is given by the formula
\begin{equation}\label{lab_67}
|\mathcal{M}| = \frac{1}{2}\biggl|\,\sum\limits_{i=1}^{k}\Delta_{i}\biggr|,\quad\text{where}\quad
\Delta_{i} = \text{det\,}\begin{pmatrix} x_{i} & y_{i} \\ x_{i+1} & y_{i+1} \end{pmatrix}
\end{equation}
(under the assumption $x_{k+1} = x_{1}$, $y_{k+1} = y_{1}$).

The application of (\ref{lab_67}) to the triangle $\mathcal{T}_{2m+1}$ yields
\begin{equation*}
-\frac{1}{2}\begin{vmatrix} \displaystyle \frac{6m-2}{6m-1} & \displaystyle \frac{4m-1}{6m-1} \\[12pt]  1 & \displaystyle \frac{4m+1}{6m+1} \end{vmatrix} -
\frac{1}{2}\begin{vmatrix}1 & \displaystyle \frac{4m+1}{6m+1} \\[12pt]  1 & \displaystyle \frac{2}{3} \end{vmatrix} - \frac{1}{2}\begin{vmatrix} 1 & \displaystyle \frac{2}{3} \\[12pt]  \displaystyle \frac{6m-2}{6m-1} & \displaystyle \frac{4m-1}{6m-1} \end{vmatrix} = \frac{1}{6(6m-1)(6m+1)}.
\end{equation*}
Lemma is proved. $\square$
\vspace{0.5cm}

\textsc{Lemma 25.} \emph{Suppose that $n\geqslant 1$ is integer. Then the region $\mathcal{T}(5,(1,4)^{n})$ is defined by the inequalities}
\[
\frac{6n-1}{6n+1}<x\leqslant 1,\quad \frac{1}{6}(1+x) < y\leqslant \frac{1+(2n+1)x}{6n+5}
\]
\emph{and is the triangle with the vertices}
\[
\biggl(\frac{6n-1}{6n+1},\frac{2n}{6n+1}\biggr),\quad \biggl(1,\frac{2n+2}{6n+5}\biggr),\quad \biggl(1,\frac{1}{3}\biggr)\quad\textit{and area}\quad \frac{1}{3(6n+1)(6n+5)}.
\]

\textsc{Lemma 26.} \emph{Suppose that $n\geqslant 1$ is integer. Then the region $\mathcal{T}(1,2^{n})$ is defined by the inequalities}
\[
0<x\leqslant \frac{1}{2n+3},\quad 1-x<y\leqslant 1,\quad \frac{1}{2n+3} < x\leqslant \frac{1}{2n+1},\quad \frac{1}{2}(1+(2n+1)x) < y\leqslant 1
\]
\emph{and is the triangle with the vertices}
\[
(0,1),\quad \biggl(\frac{1}{2n+1},1\biggr),\quad \biggl(\frac{1}{2n+3},\frac{2n+2}{2n+3}\biggr)\quad\textit{and area}\quad \frac{1}{2(2n+1)(2n+3)}.
\]
The proofs of Lemmas 25, 26 are similar to the proof of Lemma 24 and are based on the formulas of Lemmas 19, 20. $\square$

Let us call the tuples of the $\mathcal{A}_{r}^{\circ}$ and $\mathcal{A}_{r}^{*}$ as degenerate and non-degenerate, cor\-res\-pon\-din\-g\-ly.
We say that the non-degenerate tuple $\mathbf{k}_{r} = (k_{1},\ldots,k_{r})\in \mathcal{A}_{r}^{*}$ \textit{generates} the tuple  $\mathbf{k}_{r+1} = (\mathbf{k}_{r},k_{r+1}) =(k_{1},\ldots,k_{r},k_{r+1})$, if both regions $\mathcal{T}(\mathbf{k}_{r})$ and $\mathcal{T}(\mathbf{k}_{r+1})$ are non-empty. In this case, we call the tuple $\mathbf{k}_{r+1}$ the \textit{descendant} of the tuple $\mathbf{k}_{r}$. For example, the tuple $(2)$ generates the tuples $(2,1)$, $(2,3)$ and $(2,4)$, the tuple $(5,1)$ generates $(5,1,3)$ and $(5,1,4)$, and so on. In what follows, we are interested only in descendants of non-degenerate tuples; so, we will not trace the descendants of the tuples $\mathbf{k}_{r}\in \mathcal{A}_{r}^{\circ}$.

We will also say that the tuple $\mathbf{k}_{r+\ell} = (\mathbf{k}_{r},k_{r+1},\ldots,k_{r+\ell})$ \textit{is generated} by the tuple $\mathbf{k}_{r} \in \mathcal{A}_{r}^{*}$ (is its descendant), if for any $i$, $1\leqslant i\leqslant \ell-1$, the tuple $\mathbf{k}_{r+i} = (\mathbf{k}_{r},k_{r+1},\ldots,k_{r+i})$ is non-degenerate and generates the tuple $\mathbf{k}_{r+i+1} = (\mathbf{k}_{r},k_{r+i},\ldots,k_{r+i+1})$. Thus, in the chain $(1)$, $(1,3)$, $(1,3,1)$, $(1,3,1,6)$, $(1,3,1,6,1)$ each subsequent tuple is the ascendant of the previous one; all of them are the descendant of the initial tuple $(1)$.

It is natural to draw the descendants of a given tuple as an oriented tree. The degenerate tuples (whose descendants are not shown in such trees) becomes the leaf nodes.
Fig. 3 shows all the descendants of the tuple $(2)$ in first five generations. To get the best view of these, we label (on Fig. 3 and the rest figures) degenerate and non-degenerate tuples by the signs <<$\circ$>> and <<$*$>> consequently.

It appears that in the case $D=3$, $c_{0}=0$, any tuple $\mathbf{k}_{r}$ from the summation formula of Lemma 23 for $\nu(r)$, $r\geqslant 8$ is the descendant of one of three tuples: $(2,4)$, $(5,1)$ or $(1,2)$. We denote the corresponding trees as $\mathbb{A}$, $\mathbb{B}$ and $\mathbb{C}$. The tuples involved into the formula for $\nu(r)$ correspond to the vertices labeled with the sign <<$\circ$>> and living on $r$th <<level>> of each of these trees. The following three lemmas describe the structure of the trees $\mathbb{A}$, $\mathbb{B}$ and $\mathbb{C}$.

\textsc{Lemma 27.} \emph{The tree $\mathbb{A}$ has the form shown on Fig.~4 (to the left). To the right, the structure of subtree $\mathbb{A}_{n}$, $n=1,2,3,\ldots$, is given.
Moreover, the following formulas for the areas corresponding to the degenerate tuples of the subtree $\mathbb{A}_{n}$ hold true}:

\begin{align*}
& |\mathcal{T}_{5n+1}^{\,\circ}| = |\mathcal{T}(2,(4,1)^{n},3,2^{3n-2},1)| = \frac{1}{2(6n-1)(6n+1)(12n+1)};\\[6pt]
& |\mathcal{T}_{5n+2}^{\,\circ}| = |\mathcal{T}(2,(4,1)^{n},3,2^{3n-1},1)| = \frac{3(4n+1)}{2(3n+1)(6n+1)(12n+1)(12n+5)};\\[6pt]
& |\mathcal{T}_{5n+3}^{\,\circ}| = |\mathcal{T}(2,(4,1)^{n},3,2^{3n},1)| = \frac{3(2n+1)}{2(3n+1)(3n+2)(12n+5)(12n+7)};\\[6pt]
& |\mathcal{T}_{5n+4}^{\,\circ}| = |\mathcal{T}(2,(4,1)^{n},3,2^{3n+1},1)| = \frac{3(4n+3)}{2(3n+2)(6n+5)(12n+7)(12n+11)};\\[6pt]
& |\mathcal{T}_{5n+5}^{\,\circ}| = |\mathcal{T}(2,(4,1)^{n},3,2^{3n+2},1)| = \frac{1}{2(6n+5)(6n+7)(12n+11)}.
\end{align*}

\pagebreak

\begin{landscape}
\begin{center}
\includegraphics{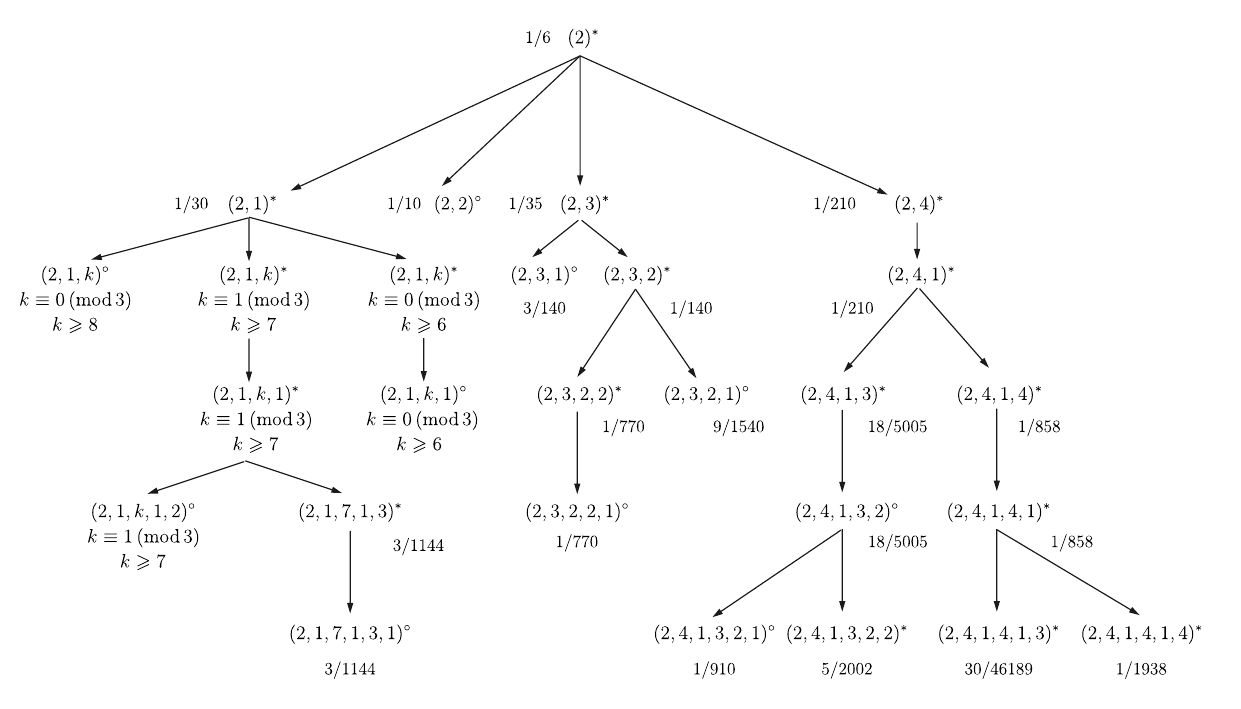}

\fontsize{10}{12pt}\selectfont
\emph{Fig.~3.} The tuple $(2)$ with the descendants till fifth generation. The fractions give the areas of the corresponding polygons.
\fontsize{12}{15pt}\selectfont
\end{center}
\end{landscape}

\pagebreak

\begin{center}
\includegraphics{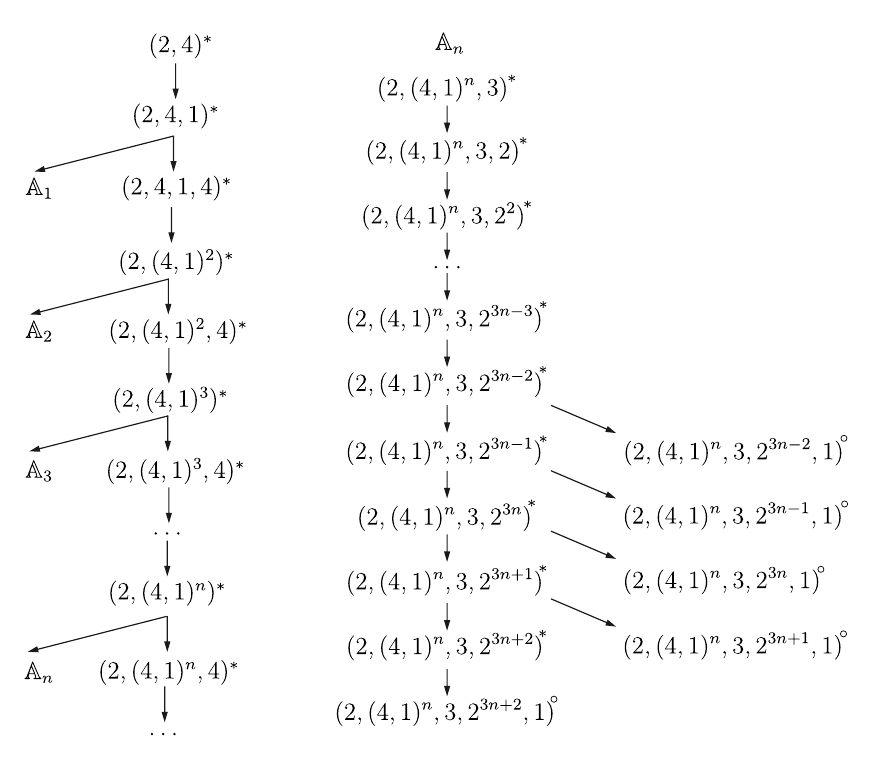}

\fontsize{10}{12pt}\selectfont
\emph{Fig.~4.} The structure of the tree $\mathbb{A}$ (to the left) and its subtree $\mathbb{A}_{n}$ (to the right).
\fontsize{12}{15pt}\selectfont
\end{center}

\textsc{Proof.} First of all, we check the degeneracy (non-degeneracy) of the tuples that label the vertices of the subtree $\mathbb{A}_{n}$. Using the formulas (\ref{lab_53}), (\ref{lab_55}), we get:
\begin{align*}
& \mathbb{K}(2,(4,1)^{n-1},4) = 6n+1\not\equiv 0\pmod{3}, \\
& \mathbb{K}(2,(4,1)^{n}) = 3n+2 \not\equiv 0 \pmod{3},\\
& \mathbb{K}(2,(4,1)^{n},3,2^{m}) = 3(m+n+1)+2\not\equiv 0\pmod{3}, \\
& \mathbb{K}(2,(4,1)^{n},3,2^{m},1) = 3 \equiv 0 \pmod{3}.
\end{align*}

In view of Lemma 24, when passing from $\mathcal{T}_{2n+1}^{*} = \mathcal{T}(2,(4,1)^{n})$ to $\mathcal{T}_{2n+2}^{*} = \mathcal{T}(2,(4,1)^{n},3)$, we need to add to the system
\[
\frac{6n-2}{6n-1}<x\leqslant 1,\quad \frac{1}{3}(1+x)<y\leqslant g_{2n}(x) = g_{2n}(x;2,(4,1)^{n-1},4) = \frac{1+4nx}{6n+1}
\]
and additional conditions
\begin{multline}\label{lab_67}
f_{2n+2}(x) = f_{2n+2}(x;2,(4,1)^{n},3) = \frac{1+(4n+4)x}{6n+7}<y\leqslant \\ \leqslant g_{2n+2}(x) = g_{2n+2}(x;2,(4,1)^{n},3) = \frac{1+(2n+3)x}{3n+5}.
\end{multline}
It is easy to check that
\begin{equation}\label{lab_68}
g_{2n}(x)\leqslant g_{2n+2}(x)
\end{equation}
for any $0\leqslant x\leqslant 1$ if $n\geqslant 2$, and for any $1/3\leqslant x\leqslant 1$ if $n=1$. Hence, the inequality (\ref{lab_68}) holds true for
$(6n-2)/(6n-1)\leqslant x \leqslant 1$. Hence, the triangle $\mathcal{T}_{2n+1}^{\,*}$ lies entirely under the segment of the line $y = g_{2n+2}(x)$, $0\leqslant x\leqslant 1$. Therefore, the upper bound for $y$ in (\ref{lab_67}) does not change the region $\mathcal{T}_{2n+1}^{\,*}$.

Further, the line $y = f_{2n+2}(x)$ intersects the line $x = 1$ at the point $\bigl(1, \tfrac{4n+5}{6n+7}\bigr)$ and intersects the line $y = (1+x)/3$ at the point $\bigl(\tfrac{6n+4}{6n+5}, \frac{4n+3}{6n+5}\bigr)$. The inequalities
\[
\frac{6n-2}{6n-1} < \frac{6n+4}{6n+5}<1,\qquad \frac{2}{3} < \frac{4n+5}{6n+7} < \frac{4n+1}{6n+1}
\]
imply that the line $y = f_{2n+2}(x)$ cuts the triangle $\mathcal{T}_{2n+1}^{\,*}$ into two parts (Fig.~5). The first one lying under this line coincides with the region $\mathcal{T}_{2n+2}^{\,*}$. It is defined by the inequalities
\begin{equation*}
\begin{cases}
\displaystyle \frac{6n-2}{6n-1} < x\leqslant \frac{6n+4}{6n+5},\quad \frac{1}{3}(1+x)<y\leqslant \frac{1+4nx}{6n+1},\\[12pt]
\displaystyle \frac{6n+4}{6n+5} < x\leqslant 1,\quad \frac{1+(4n+4)x}{6n+7}<y\leqslant \frac{1+4nx}{6n+1}
\end{cases}
\end{equation*}
and is the quadrangle with vertices
\[
\biggl(\frac{6n-2}{6n-1},\,\frac{4n-1}{6n-1}\biggr),\quad \biggl(1,\,\frac{4n+1}{6n+1}\biggr),\quad \biggl(1,\,\frac{4n+5}{6n+7}\biggr),
\quad \biggl(\frac{6n+4}{6n+5},\,\frac{4n+3}{6n+5}\biggr)
\]
Its area equals to
\[
|\mathcal{T}_{2n+2}^{\,*}| = \frac{6(2n+1)}{(6n-1)(6n+1)(6n+5)(6n+7)}.
\]

\begin{center}
\includegraphics{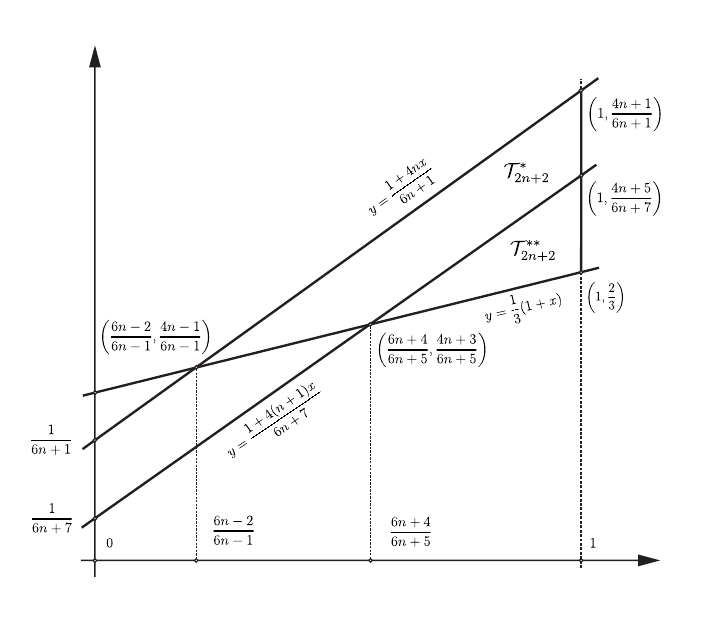}
\fontsize{10}{12pt}\selectfont

\emph{Fig.~5.} The regions $\mathcal{T}_{2n+2}^{\,*}$ and $\mathcal{T}_{2n+2}^{\,**}$.
\fontsize{12}{15pt}\selectfont
\end{center}

Similarly, in change $\mathcal{T}_{2n+1}^{\,*}\rightarrow \mathcal{T}_{2n+2}^{\,**} = \mathcal{T}_{2n+2}(2,(4,1)^{n},4)$ (here, the sign <<$**$>> is used instead of <<$*$>> because the notation $\mathcal{T}_{2n+2}^{*}$ was used earlier for the neighbour region $\mathcal{T}_{2n+2}(2,(4,1)^{n},3)$), an extra conditions appear, namely
\begin{equation}\label{lab_69}
f_{2n+2}(x;2,(4,1)^{n},4) = \frac{1+(6n+5)x}{9n+9}<y\leqslant g_{2n+2}(x;2,(4,1)^{n},4) = \frac{1+(4n+4)x}{6n+7}.
\end{equation}
Arguing as above, one can check that the lower bound for $y$ in (\ref{lab_69}) does not change the region $\mathcal{T}_{2n+1}^{\,*}$, and the upper bound coincides with the lower one given by (\ref{lab_67}). Hence, the region $\mathcal{T}_{2n+2}^{\,**}$ coincides with the part of $\mathcal{T}_{2n+1}^{\,*}$ that lies under the the line $y = g_{2n+2}(x)$ and is defined by the inequalities
\[
\frac{6n+4}{6n+5}<y\leqslant 1,\qquad \frac{1}{3}(1+x)<y\leqslant \frac{1+(4n+4)x}{6n+7}.
\]
Obviously, it is the triangle with vertices
\[
\biggl(\frac{6n+4}{6n+5},\,\frac{4n+3}{6n+5}\biggr),\quad \biggl(1,\,\frac{4n+5}{6n+7}\biggr),\quad \biggl(1,\,\frac{2}{3}\biggr)
\]
and with area
\[
|\mathcal{T}_{2n+2}^{\,**}| = \frac{1}{6(6n+5)(6n+7)}.
\]

The above arguments imply that all the regions $\mathcal{T}(2,(4,1)^{n},k)$ are empty for $k\ne 3,4$. This fact is expressed by the following equation in areas:
\begin{multline*}
|\mathcal{T}_{2n+1}^{\,*}| = \frac{1}{6(6n-1)(6n+1)} = \\[6pt]
=\frac{6(2n+1)}{(6n-1)(6n+1)(6n+5)(6n+7)} + \frac{1}{6(6n+5)(6n+7)} = |\mathcal{T}_{2n+2}^{\,*}| +  |\mathcal{T}_{2n+2}^{\,**}|.
\end{multline*}

The extra condition
\[
f_{2n+3}(x;2,(4,1)^{n+1}) = \frac{1+(6n+7)x}{9n+12}<y\leqslant g_{2n+3}(x;2,(4,1)^{n+1}) = \frac{1+(2n+3)x}{3n+5}
\]
that appears in the change $\mathcal{T}_{2n+2}^{\,**} \rightarrow \mathcal{T}_{2n+3}^{\,*} = \mathcal{T}(2,(4,1)^{n+1})$ does not change the region
$\mathcal{T}_{2n+2}^{\,**}$ due to the inequalities
\[
\frac{1+(6n+7)x}{9n+12} \leqslant \frac{1}{3}(1+x),\quad \frac{1+(4n+4)x}{6n+7}\leqslant \frac{1+(2n+3)x}{3n+5},\quad 0\leqslant x\leqslant 1.
\]

Now we show that any region $\mathcal{T}_{2n+2+m}^{\,*} = \mathcal{T}(2,(4,1)^{n},3,2^{m})$, $0\leqslant m\leqslant 3n-2$, coincides with the quadrangle $\mathcal{T}_{2n+2}^{\,*}$. For $m = 0$, it is obvious. Suppose that this assertion is checked for any $m$, $0\leqslant m\leqslant \ell-1$ (we assume that $\ell\leqslant 3n-2$), and prove it for $m = \ell$.

In the change $\mathcal{T}_{2n+2+(\ell-1)}^{\,*}\rightarrow \mathcal{T}_{2n+2+\ell}^{\,*}$ (or, that is the same in view of the induction, in the change $\mathcal{T}_{2n+2}^{\,*}\rightarrow \mathcal{T}_{2n+2+\ell}^{\,*}$), the following extra condition arises:
\begin{equation}\label{lab_70}
f(x) = f(x;2,(4,1)^{n},3,2^{\ell}) <y\leqslant g(x) = g(x;2,(4,1)^{n},3,2^{\ell}).
\end{equation}
Using Lemma 20, we find the explicit form of the inequalities (\ref{lab_70}), namely
\begin{multline*}
f(x) = \frac{1 + x\cdot\mathbb{K}((4,1)^{n},3,2^{\ell-1},3)}{\mathbb{K}(2,(4,1)^{n},3,2^{\ell-1},3)} = \frac{1 + 4(n+\ell+1)x}{6(n+\ell+1)+1},\\
g(x) = \frac{1 + x\cdot\mathbb{K}((4,1)^{n},3,2^{\ell})}{\mathbb{K}(2,(4,1)^{n},3,2^{\ell})} = \frac{1 + (2(n+\ell)+3)x}{3(n+\ell)+5}.
\end{multline*}
Note that the inequality $f(x)\leqslant (1+x)/3$ is equivalent to the inequality
\[
x\leqslant \frac{6n+4+6\ell}{6n+5+6\ell},
\]
and the last one is satisfied by all $x\leqslant (6n+4)/(6n+5)$ for any $\ell\geqslant 0$. Similarly, the inequality
 $f(x)\leqslant (1+(4n+4)x)/(6n+7)$ is equivalent to the inequality $2\ell x\leqslant 6\ell$, which holds true for all $x$, $0\leqslant x\leqslant 1$, and for any $\ell\geqslant 0$. Therefore, the condition $y>f(x)$ in (\ref{lab_70}) does not change the region $\mathcal{T}_{2n+2}^{*}$ for any $\ell\geqslant 0$.

Further, the inequality $g(x)>(1+4nx)/(6n+1)$ is reduced to $3(\ell-n)+4\leqslant (2\ell+3)x$.
If $0\leqslant \ell \leqslant n-2$, then $3(\ell - n)+4\leqslant -2<0\leqslant (2\ell+3)x$ for all $x$, $0\leqslant x\leqslant 1$. If $\ell\geqslant n-1$, then this inequality is equivalent to
\[
x\geqslant \frac{3(\ell - n)+4}{2\ell+3}.
\]
However, if $\ell\leqslant 3n-2$ then this inequality is satisfied by any $x$ such that
\[
\frac{6n-2}{6n-1}\leqslant x\leqslant 1.
\]
This follows from the obvious relations
\[
\frac{6n-2}{6n-1} - \frac{3(\ell - n)+4}{2\ell+3} = \frac{(6n+1)(3n-2-\ell)}{(6n-1)(2\ell+3)}\geqslant 0.
\]
Thus, the condition $y\leqslant g(x)$ in (\ref{lab_70}) does not split the region $\mathcal{T}_{2n+2}^{\,*}$ for any $0\leqslant\ell\leqslant 3n-2$. Hence, we have the following chain of region changes:
\begin{multline*}
\mathcal{T}_{2n+2}^{\,*}\rightarrow \mathcal{T}_{2n+3}^{\,*} = \mathcal{T}(2,(4,1)^{n},3,2)\rightarrow \mathcal{T}_{2n+4}^{\,*} = \mathcal{T}(2,(4,1)^{n},3,2^{2})\rightarrow\cdots\\
\cdots \rightarrow \mathcal{T}_{5n}^{\,*} = \mathcal{T}(2,(4,1)^{n},3,2^{3n-2}).
\end{multline*}

In the change $\mathcal{T}_{5n}^{\,*}\rightarrow \mathcal{T}_{5n+1}^{\,\circ} = \mathcal{T}(2,(4,1)^{n},3,2^{3n-2},1)$, the following extra conditions arise
\[
f(x;2,(4,1)^{n},3,2^{3n-2},1) = \frac{1+(8n+1)x}{12n+2} <y\leqslant g(x;2,(4,1)^{n},3,2^{3n-2},1) = \frac{1}{3}(1+2x).
\]
Similarly, in the change $\mathcal{T}_{5n}^{\,*}\rightarrow \mathcal{T}_{5n+1}^{\,*} = \mathcal{T}(2,(4,1)^{n},3,2^{3n-1})$, an extra conditions take the form
\[
f(x;2,(4,1)^{n},3,2^{3n-1}) = \frac{1+16nx}{24n+1} <y\leqslant g(x;2,(4,1)^{n},3,2^{3n-1}) = \frac{1+(8n+1)x}{12n+2}.
\]
It is not difficult to check that the inequalities
\[
y\leqslant \frac{1}{3}(1+x)\quad\text{и}\quad y > \frac{1+16nx}{24n+1}
\]
do not cut split the region $\mathcal{T}_{5n}^{\,*}$. Further, the line $y = (1+(8n+1)x)/(12n+2)$ passes through the vertex $\bigl(1,\tfrac{4n+1}{6n+1}\bigr)$ of the quadrangle $\mathcal{T}_{5n}^{\,*}$ and intersect the line $y = (1+x)/3$ at the point $\bigl(\tfrac{12n-1}{12n+1},\tfrac{8n}{12n+1}\bigr)$. By the inequalities
\[
\frac{6n-2}{6n-1} < \frac{12n-1}{12n+1} < \frac{6n+4}{6n+5},
\]
we conclude that the line $y = (1+(8n+1)x)/(12n+2)$ split $\mathcal{T}_{5n}^{\,*}$ into two parts, namely, to $\mathcal{T}_{5n+1}^{\,\circ}$ (which lies over this secant) and to $\mathcal{T}_{5n+1}^{\,*}$ (lying under the secant; see. Fig.~6). Hence, $\mathcal{T}_{5n+1}^{\circ}$ is the triangle with vertices
\[
\biggl(\frac{6n-2}{6n-1},\,\frac{4n-1}{6n-1}\biggr),\quad \biggl(1,\,\frac{4n+1}{6n+1}\biggr),\quad \biggl(\frac{12n-1}{12n+1},\,\frac{8n}{12n+1}\biggr)
\]
and with area
\[
|\mathcal{T}_{5n+1}^{\,\circ}| = \frac{1}{2(6n-1)(6n+1)(12n+1)}.
\]
Consequently, $\mathcal{T}_{5n+1}^{\,*}$ is defined by the system
\begin{equation*}
\begin{cases}
\displaystyle \frac{12n-1}{12n+1} < x\leqslant \frac{6n+4}{6n+5},\quad  \frac{1}{3}(1+x)<y\leqslant \frac{1+(8n+1)x}{12n+2},\\[12pt]
\displaystyle \frac{6n+4}{6n+5} < x\leqslant 1,\quad  \frac{1+(4n+4)x}{6n+7}<y\leqslant \frac{1+(8n+1)x}{12n+2}
\end{cases}
\end{equation*}
and is the quadrangle with vertices
\[
\biggl(\frac{12n-1}{12n+1},\,\frac{8n}{12n+1}\biggr),\quad \biggl(1,\,\frac{4n+1}{6n+1}\biggr),\quad \biggl(1,\,\frac{4n+5}{6n+7}\biggr),\quad \biggl(\frac{6n+4}{6n+5},\,\frac{4n+3}{6n+5}\biggr),
\]
whose area equals to
\[
|\mathcal{T}_{5n+1}^{\,*}| = \frac{42n+23}{2(6n+1)(6n+5)(6n+7)(12n+1)}.
\]
One can also check the following equation in areas:
\begin{multline*}
|\mathcal{T}_{5n+1}^{\,\circ}|+|\mathcal{T}_{5n+1}^{\,*}|  = \frac{1}{2(6n-1)(6n+1)(12n+1)} + \frac{42n+23}{2(6n+1)(6n+5)(6n+7)(12n+1)} = \\[6pt]
= \frac{6(2n+1)}{(6n-1)(6n+1)(6n+5)(6n+7)} = |\mathcal{T}_{5n}^{\,*}|  = |\mathcal{T}_{2n+2}^{\,*}|.
\end{multline*}

\begin{center}
\includegraphics{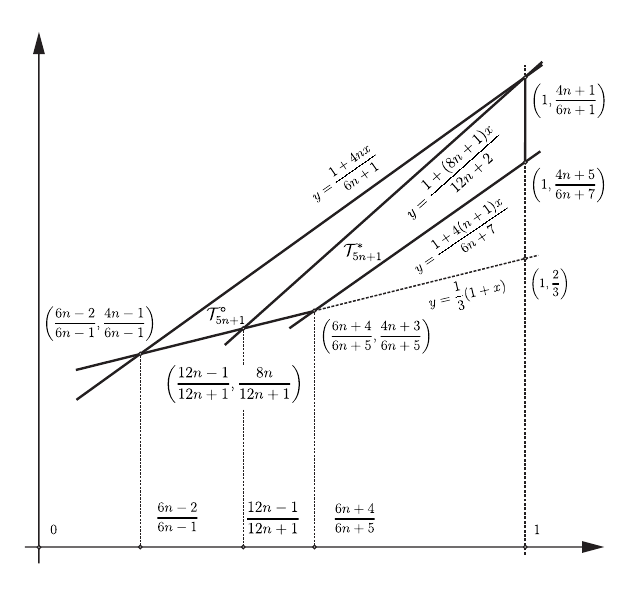}
\fontsize{10}{12pt}\selectfont

\emph{Fig.~6.} The polygons $\mathcal{T}_{5n+1}^{\,*}$ and $\mathcal{T}_{5n+1}^{\,\circ}$.
\fontsize{12}{15pt}\selectfont
\end{center}
In the change $\mathcal{T}_{5n+1}^{\,*} \rightarrow \mathcal{T}_{5n+2}^{\,\circ} = \mathcal{T}(2,(4,1)^{n},3,2^{3n-1},1)$, the following extra condition arises:
\[
f(x;2,(4,1)^{n},3,2^{3n-1},1) < y\leqslant g(x;2,(4,1)^{n},3,2^{3n-1},1).
\]
that is,
\[
\frac{1+(8n+3)x}{12n+5}<y\leqslant \frac{1}{3}(1+2x).
\]
Similarly, an extra condition that arises in the change $\mathcal{T}_{5n+1}^{\,*} \rightarrow \mathcal{T}_{5n+2}^{\,*} = \mathcal{T}(2,(4,1)^{n},3,2^{3n})$, has the form
\[
f(x;2,(4,1)^{n},3,2^{3n}) < y\leqslant g(x;2,(4,1)^{n},3,2^{3n})
\]
or, that is the same,
\[
\frac{1+(16n+4)x}{24n+7}<y\leqslant \frac{1+(8n+3)x}{12n+5}.
\]
One can check that
\[
\frac{1+(8n+1)x}{12n+2}< \frac{1}{3}(1+2x)\quad\text{for}\quad 0\leqslant x\leqslant 1,
\]
and
\[
\frac{1+(16n+4)x}{24n+7}\leqslant \frac{1}{3}(1+x)\quad\text{for}\quad x\leqslant \frac{24n+4}{24n+5}
\]
(and hence for $(12n-1)/(12n+1)\leqslant x\leqslant (6n+4)/(6n+5)$), and check that
\[
\frac{1+(16n+4)x}{24n+7}\leqslant \frac{1+(4n+4)x}{6n+7}\quad\text{for}\quad 0\leqslant x\leqslant 1
\]
(and hence for $(6n+4)/(6n+5)\leqslant x\leqslant 1$). Therefore, the inequalities
\[
y\leqslant \frac{1}{3}(1+2x),\quad y>\frac{1+(16n+4)x}{24n+7}
\]
do not split the region $\mathcal{T}_{5n+1}^{\,*}$. Further, the line $y = (1+(8n+3)x)/(12n+5)$ intersects the lined $y = (1+x)/3$ and $x = 1$ at the points
\[
\biggl(\frac{6n+1}{6n+2},\,\frac{4n+1}{6n+2}\biggr)\quad\text{and}\quad \biggl(1,\,\frac{8n+4}{12n+5}\biggr)
\]
consequently. Thus, it splits $\mathcal{T}_{5n+1}^{\,*}$ into two quadrangles: $\mathcal{T}_{5n+2}^{\,\circ}$ и $\mathcal{T}_{5n+2}^{\,*}$ (Fig.~7). The first one has its vertices at the points
\[
\biggl(\frac{12n-1}{12n+1},\frac{8n}{12n+1}\biggr),\quad \biggl(1,\frac{4n+1}{6n+1}\biggr),\quad \biggl(1,\frac{8n+4}{12n+5}\biggr),\quad \biggl(\frac{6n+1}{6n+2},\frac{4n+1}{6n+2}\biggr),
\]
and its area equals to
\[
|\mathcal{T}_{5n+2}^{\circ}| = \frac{3(4n+1)}{2(3n+1)(6n+1)(12n+1)(12n+5)}.
\]

\begin{center}
\includegraphics{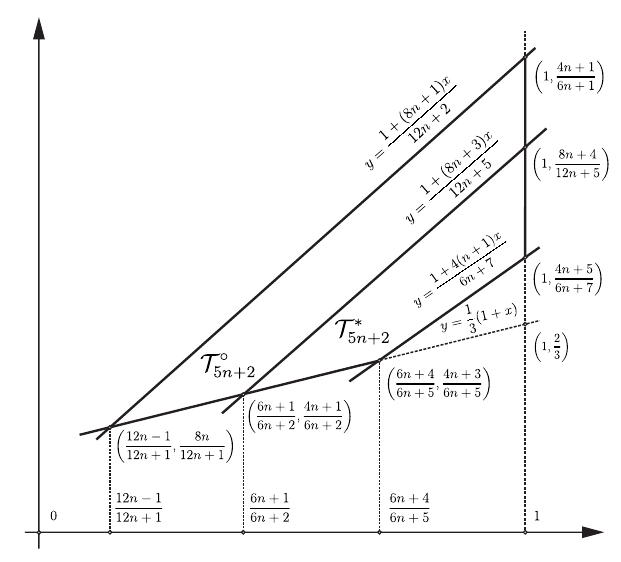}

\fontsize{10}{12pt}\selectfont
\emph{Fig.~7.} The regions $\mathcal{T}_{5n+2}^{\,*}$ и $\mathcal{T}_{5n+2}^{\,\circ}$.
\fontsize{12}{15pt}\selectfont
\end{center}

The second quadrangle, that is, $\mathcal{T}_{5n+2}^{\,*}$ is defined by the inequalities
\begin{equation*}
\begin{cases}
\displaystyle\frac{6n+1}{6n+2} < x\leqslant \frac{6n+4}{6n+5},\quad  \frac{1}{3}(1+x)<y\leqslant \frac{1+(8n+3)x}{12n+5},\\[12pt]
\displaystyle\frac{6n+4}{6n+5} < x\leqslant 1,\quad  \frac{1+(4n+4)x}{6n+7}<y\leqslant \frac{1+(8n+3)x}{12n+5}.
\end{cases}
\end{equation*}
Its vertices are the points
\[
\biggl(\frac{6n+1}{6n+2},\,\frac{4n+1}{6n+2}\biggr),\quad \biggl(1,\,\frac{8n+4}{12n+5}\biggr),\quad \biggl(1,\,\frac{4n+5}{6n+7}\biggr),\quad \biggl(\frac{6n+4}{6n+5},\,\frac{4n+3}{6n+5}\biggr),
\]
and its area is equal to
\[
|\mathcal{T}_{5n+2}^{\,*}| = \frac{5(3n+2)}{2(3n+1)(6n+5)(6n+7)(12n+5)}.
\]
One can also check that
\begin{multline*}
|\mathcal{T}_{5n+2}^{\,*}|+|\mathcal{T}_{5n+2}^{\,\circ}| = \\[6pt]
=\frac{3(4n+1)}{2(3n+1)(6n+1)(12n+1)(12n+5)} + \frac{5(3n+2)}{2(3n+1)(6n+5)(6n+7)(12n+5)} =\\[6pt]
= \frac{42n+23}{2(6n+1)(6n+5)(6n+7)(12n+1)} = |\mathcal{T}_{5n+1}^{\,*}|.
\end{multline*}

Further, in the change $\mathcal{T}_{5n+2}^{\,*}\rightarrow \mathcal{T}_{5n+3}^{\,\circ} = \mathcal{T}(2,(4,1)^{n},3,2^{3n},1)$ the following extra condition arises:
\[
f(x;2,(4,1)^{n},3,2^{3n},1)<y\leqslant g(x;2,(4,1)^{n},3,2^{3n},1),
\]
or, that is the same,
\[
\frac{1+(8n+5)x}{12n+8}<y\leqslant \frac{1}{3}(1+2x).
\]
Similarly, the extra condition that arises in the change $\mathcal{T}_{5n+2}^{\,*}\rightarrow \mathcal{T}_{5n+3}^{\,*} = \mathcal{T}(2,(4,1)^{n},3,2^{3n+1})$ has the form
\[
f(x;2,(4,1)^{n},3,2^{3n+1})<y\leqslant g(x;2,(4,1)^{n},3,2^{3n+1})
\]
or, that is the same,
\[
\frac{1+(16n+8)x}{24n+13}<y\leqslant \frac{1+(8n+5)x}{12n+8}.
\]
One can easily check that
\[
\frac{1+(8n+3)x}{12n+5}\leqslant \frac{1}{3}(1+2x)\quad\text{for}\quad 0\leqslant x\leqslant 1.
\]
Hence, the condition $y \leqslant (1+2x)/3$ does not affect the region $\mathcal{T}_{5n+2}^{\,*}$. It easy to check that
\[
\frac{1+(16n+8)x}{24n+13}\leqslant \frac{1+(4n+4)x}{6n+7}\quad\text{for any}\quad 0\leqslant x\leqslant 1.
\]
Therefore, the condition $y > (1+(16n+8)x)/(24n+13)$ does not affect the region $\mathcal{T}_{5n+2}^{\,*}$. Finally, the line $y = (1+(8n+5)x)/(12n+8)$ meets the lines $x=1$ and $y = (1+x)/3$ at the points
\[
\biggl(1,\,\frac{4n+3}{6n+4}\biggr),\quad \biggl(\frac{12n+5}{12n+7},\,\frac{8n+4}{12n+7}\biggr).
\]
Its comparison with the corresponding coordinates of vertices of the quadrangle $\mathcal{T}_{5n+2}^{\,*}$ yields:
\[
\frac{6n+1}{6n+2} < \frac{12n+5}{12n+7} < \frac{6n+4}{6n+5},\quad \quad \frac{4n+5}{6n+7} < \frac{4n+3}{6n+7} < \frac{8n+4}{12n+5}.
\]
Hence, the line $y = (1+(8n+5)x)/(12n+8)$ split the region $\mathcal{T}_{5n+2}^{\,*}$ into two quadrangles, namely $|\mathcal{T}_{5n+3}^{\,\circ}|$ and $|\mathcal{T}_{5n+3}^{\,*}|$ (Fig.~8). The first one has its vertices at the points
\[
\biggl(\frac{6n+1}{6n+2},\,\frac{4n+1}{6n+2}\biggr),\quad \biggl(1,\,\frac{8n+4}{12n+5}\biggr),\quad \biggl(1,\,\frac{4n+3}{6n+4}\biggr),\quad \biggl(\frac{12n+5}{12n+7},\,\frac{8n+4}{12n+7}\biggr),
\]
and its area is equal to
\[
|\mathcal{T}_{5n+3}^{\circ}| = \frac{3(2n+1)}{2(3n+1)(3n+2)(12n+5)(12n+7)}.
\]
The vertices of the second one are the points
\[
\biggl(\frac{12n+5}{12n+7},\,\frac{8n+4}{12n+7}\biggr),\quad \biggl(1,\,\frac{4n+3}{6n+4}\biggr),\quad \biggl(1,\,\frac{4n+5}{6n+7}\biggr),\quad \biggl(\frac{6n+4}{6n+5},\,\frac{4n+3}{6n+5}\biggr),
\]
its are equals to
\[
|\mathcal{T}_{5n+3}^{\,*}| = \frac{9n+7}{2(3n+2)(6n+5)(6n+7)(12n+7)},
\]
and it is defined by the system
\begin{equation*}
\begin{cases}
\displaystyle \frac{12n+5}{12n+7} < x\leqslant \frac{6n+4}{6n+5},\quad  \frac{1}{3}(1+x)<y\leqslant \frac{1+(8n+5)x}{12n+8},\\[12pt]
\displaystyle \frac{6n+4}{6n+5} < x\leqslant 1,\quad  \frac{1+(4n+4)x}{6n+7}<y\leqslant \frac{1+(8n+5)x}{12n+8}.
\end{cases}
\end{equation*}
One can check that
\begin{multline*}
|\mathcal{T}_{5n+3}^{\,\circ}|+|\mathcal{T}_{5n+3}^{\,*}| = \\[6pt]
= \frac{3(2n+1)}{2(3n+1)(3n+2)(12n+5)(12n+7)} + \frac{9n+7}{2(3n+2)(6n+5)(6n+7)(12n+7)} =\\[6pt]
= \frac{5(3n+2)}{2(3n+1)(6n+5)(6n+7)(12n+5)} = |\mathcal{T}_{5n+2}^{\,*}|.
\end{multline*}

\begin{center}
\includegraphics{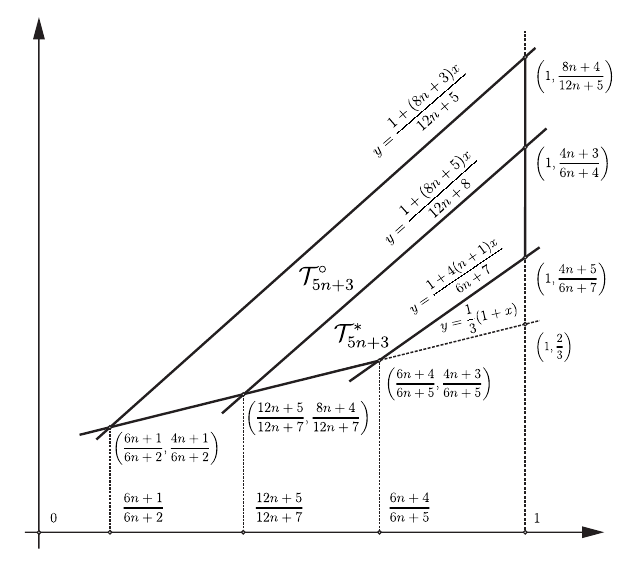}

\fontsize{10}{12pt}\selectfont
\emph{Fig.~8.} The regions $\mathcal{T}_{5n+3}^{\,*}$ and $\mathcal{T}_{5n+3}^{\,\circ}$.
\fontsize{12}{15pt}\selectfont
\end{center}

In the change $\mathcal{T}_{5n+3}^{\,*}\rightarrow \mathcal{T}_{5n+4}^{\,\circ} = \mathcal{T}(2,(4,1)^{n},3,2^{3n+1},1)$, the following extra condition arises
\[
f(x;2,(4,1)^{n},3,2^{3n+1},1)<y\leqslant g(x;2,(4,1)^{n},3,2^{3n+1},1),
\]
that is,
\[
\frac{1+(8n+7)x}{12n+11}<y\leqslant \frac{1}{3}(1+2x).
\]
The change $\mathcal{T}_{5n+3}^{\,*}\rightarrow \mathcal{T}_{5n+4}^{\,*} = \mathcal{T}(2,(4,1)^{n},3,2^{3n+2})$ yields the condition
\[
f(x;2,(4,1)^{n},3,2^{3n+2})<y\leqslant g(x;2,(4,1)^{n},3,2^{3n+2}),
\]
that is,
\[
\frac{1+(16n+12)x}{24n+19}<y\leqslant \frac{1+(8n+7)x}{12n+11}.
\]
Since
\[
\frac{1+(16n+12)x}{24n+19}\leqslant \frac{1+(4n+4)x}{6n+7}\quad\text{for}\quad 0\leqslant x\leqslant 1,
\]
we conclude that the condition $y>(1+(16n+12)x)/(24n+19)$ does not affect the region $\mathcal{T}_{5n+3}^{\,*}$. Further,
\[
\frac{1+(8n+5)x}{12n+8}\leqslant \frac{1}{3}(1+2x)\quad\text{for}\quad 0\leqslant x\leqslant 1,
\]
so the condition $y\leqslant (1+2x)/3$ does not affect the region $\mathcal{T}_{5n+3}^{*}$. Finally, the line $y = (1+(8n+7)x)/(12n+11)$ meets the lines $x=1$ and $y = (1+x)/3$ at the points
\[
\biggl(1,\,\frac{8n+8}{12n+11}\biggr)\quad\text{и}\quad \biggl(\frac{6n+4}{6n+5},\,\frac{4n+3}{6n+5}\biggr).
\]
The last one is the vertex of the quadrangle $\mathcal{T}_{5n+3}^{\,*}$. The comparison of $y$-coordinates of the points lying on the line $x=1$ yields
\[
\frac{4n+5}{6n+7} < \frac{8n+8}{12n+11} < \frac{4n+3}{6n+4}.
\]
Thus we conclude that the line $y = (1+(8n+7)x)/(12n+11)$ splits $\mathcal{T}_{5n+3}^{\,*}$ into two parts, namely $\mathcal{T}_{5n+4}^{\,\circ}$ and $\mathcal{T}_{5n+4}^{\,*}$ (Fig.~9).
The first one is the quadrangle with vertices
\[
\biggl(\frac{12n+5}{12n+7},\,\frac{8n+4}{12n+7}\biggr),\quad \biggl(1,\,\frac{4n+3}{6n+4}\biggr),\quad \biggl(1,\,\frac{8n+8}{12n+11}\biggr),\quad \biggl(\frac{6n+4}{6n+5},\,\frac{4n+3}{6n+5}\biggr);
\]
its area is equal to
\[
|\mathcal{T}_{5n+4}^{\,\circ}| = \frac{3(4n+3)}{2(3n+2)(6n+5)(12n+7)(12n+11)}.
\]
The second one is the triangle defined by the inequalities
\[
\frac{6n+4}{6n+5}<x\leqslant 1,\quad \frac{1+(4n+4)x}{6n+7}<y\leqslant \frac{1+(8n+7)x}{12n+11}.
\]
It has its vertices at the points
\[
\biggl(\frac{6n+4}{6n+5},\,\frac{4n+3}{6n+5}\biggr),\quad \biggl(1,\,\frac{8n+8}{12n+11}\biggr),\quad \biggl(1,\,\frac{4n+5}{6n+7}\biggr)
\]
and its area is equal to
\[
|\mathcal{T}_{5n+4}^{\,*}| = \frac{1}{2(6n+5)(6n+7)(12n+11)}.
\]
One can check also, that
\begin{multline*}
|\mathcal{T}_{5n+4}^{\,\circ}|+|\mathcal{T}_{5n+4}^{\,*}| = \\[6pt]
= \frac{3(4n+3)}{2(3n+2)(6n+5)(12n+7)(12n+11)} + \frac{1}{2(6n+5)(6n+7)(12n+11)} = \\[6pt]
=\frac{9n+7}{2(3n+2)(6n+5)(6n+7)(12n+7)} = |\mathcal{T}_{5n+3}^{\,*}|.
\end{multline*}

\begin{center}
\includegraphics{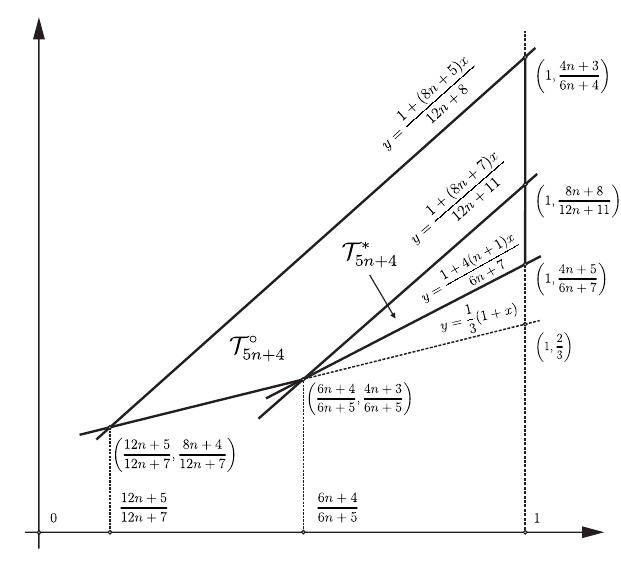}

\fontsize{10}{12pt}\selectfont
\emph{Fig.~9.} The regions $\mathcal{T}_{5n+4}^{\,*}$ and $\mathcal{T}_{5n+4}^{\,\circ}$.
\fontsize{12}{15pt}\selectfont
\end{center}

Finally, the regions $\mathcal{T}_{5n+5}^{\,\circ} = \mathcal{T}(2,(4,1)^{n},3,2^{3n+2},1)$ and $\mathcal{T}_{5n+4}^{\,*}$ are coincide. Indeed, the change $\mathcal{T}_{5n+4}^{\,*}\rightarrow \mathcal{T}_{5n+5}^{\,\circ}$ yields to the condition
\[
f(x;2,(4,1)^{n},3,2^{3n+2},1)<y\leqslant g(x;2,(4,1)^{n},3,2^{3n+2},1),
\]
that is,
\[
\frac{1+(8n+9)x}{12n+14}<y\leqslant \frac{1}{3}(1+2x).
\]
The calculation shows that
\[
\frac{1+(8n+7)x}{12n+11}\leqslant \frac{1}{3}(1+2x)\quad\text{for any}\quad 0\leqslant x\leqslant 1.
\]
This means that the condition $y\leqslant (1+2x)/3$ does not affect to the region $\mathcal{T}_{5n+4}^{\,*}$. Similarly,
\[
\frac{1+(8n+9)x}{12n+14}\leqslant \frac{1+(4n+4)x}{6n+7}\quad\text{for any}\quad 0\leqslant x\leqslant 1.
\]
Hence, the condition $y>(1+(8n+9)x)/(12n+14)$ does not affect to the region $\mathcal{T}_{5n+4}^{\,*}$ also. Therefore,
\[
|\mathcal{T}_{5n+5}^{\,\circ}| = |\mathcal{T}_{5n+4}^{*}| =\frac{1}{2(6n+5)(6n+7)(12n+11)}.
\]
Lemma 27 is proved. $\square$
\vspace{0.3cm}

\textsc{Lemma 28.} \emph{The form of the tree $\mathbb{B}$ is shown on Fig.~10 (to the left). At the same figure, the subtree $\mathbb{B}_{n}$, $n=1,2,3,\ldots$, is shown (to the right).
Moreover, the following formulas for the areas of polygons corresponding to the degenerate tuples of the subtree $\mathbb{B}_{n}$ hold true:}

\begin{align*}
& |\mathcal{T}_{5n+3}^{\,\circ}| = |\mathcal{T}(5,(1,4)^{n},1,3,2^{3n-1},1)| = \frac{1}{6(2n+1)(6n+1)(12n+5)};\\[6pt]
& |\mathcal{T}_{5n+4}^{\,\circ}| = |\mathcal{T}(5,(1,4)^{n},1,3,2^{3n},1)| = \frac{9n+5}{6(2n+1)(3n+2)(6n+5)(12n+5)};\\[6pt]
& |\mathcal{T}_{5n+5}^{\,\circ}| = |\mathcal{T}(5,(1,4)^{n},1,3,2^{3n+1},1)| = \frac{3(8n+7)}{2(3n+2)(6n+5)(12n+11)(12n+13)};\\[6pt]
& |\mathcal{T}_{5n+6}^{\,\circ}| = |\mathcal{T}(5,(1,4)^{n},1,3,2^{3n+2},1)| = \frac{3(8n+9)}{2(3n+4)(6n+7)(12n+11)(12n+13)};\\[6pt]
& |\mathcal{T}_{5n+7}^{\,\circ}| = |\mathcal{T}(5,(1,4)^{n},1,3,2^{3n+3},1)| = \frac{9n+13}{6(2n+3)(3n+4)(6n+7)(12n+19)};\\[6pt]
& |\mathcal{T}_{5n+8}^{\,\circ}| = |\mathcal{T}(5,(1,4)^{n},1,3,2^{3n+4},1)| = \frac{1}{6(2n+3)(6n+11)(12n+19)}.
\end{align*}

\begin{center}
\includegraphics{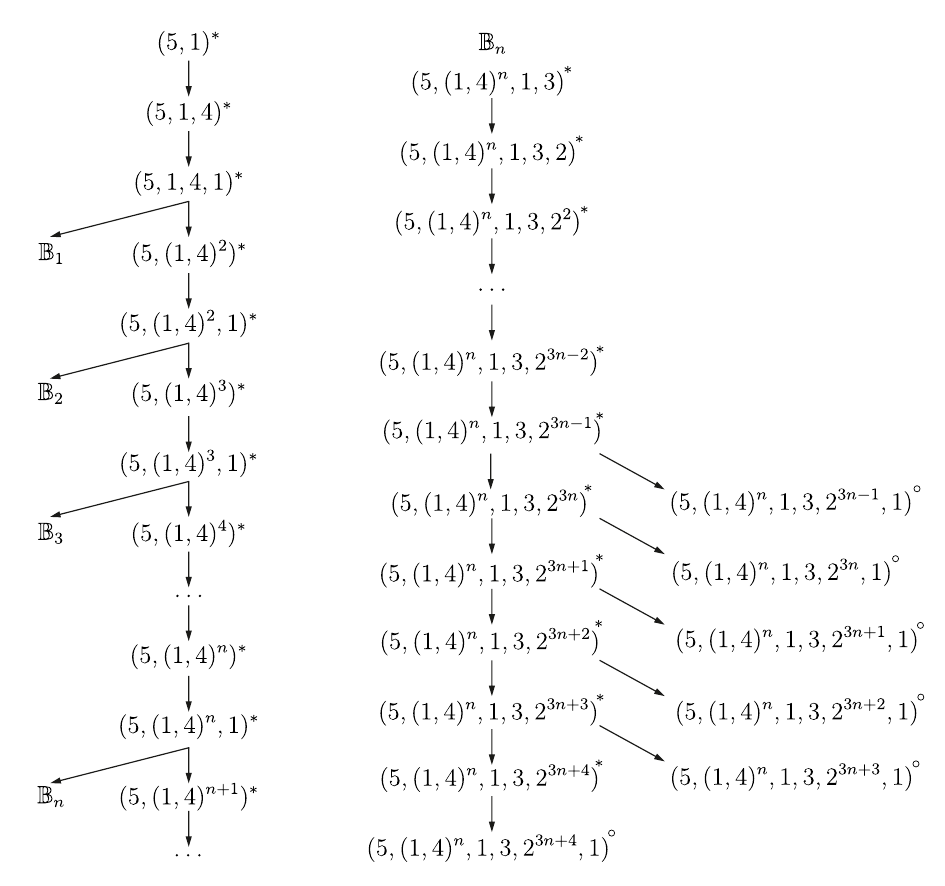}

\fontsize{10}{12pt}\selectfont
\emph{Fig.~10.} The tree $\mathbb{B}$ (to the left) and its subtree $\mathbb{B}_{n}$ (to the right).
\fontsize{12}{15pt}\selectfont
\end{center}

\textsc{Proof} of Lemma 28 is similar to the proof of the previous lemma and is based on the formulas (\ref{lab_55}) и (\ref{lab_56}) of Lemma 20. For this reason, for any leaf of the tree, we give the final system of inequalities that define the corresponding polygons, and write the list of its vertices with its area.
Let $\mathcal{T}(5,(1,4)^{n}) = \mathcal{T}_{2n+1}^{\,*}$.

The change $\mathcal{T}_{2n+1}^{\,*}\rightarrow \mathcal{T}_{2n+2}^{\,*} = \mathcal{T}(5,(1,4)^{n},1)$ does not affect the initial region. Hence, $\mathcal{T}_{2n+2}^{\,*}$ is the triangle
defined by the system
\[
\frac{6n-1}{6n+1} < x\leqslant 1,\quad \frac{1}{6}(1+x)<y\leqslant \frac{1+(2n+1)x}{6n+5}.
\]
Its vertices are the points
\[
\biggl(\frac{6n-1}{6n+1},\,\frac{2n}{6n+1}\biggr),\quad \biggl(1,\,\frac{2n+2}{6n+5}\biggr),\quad \biggl(1,\,\frac{1}{3}\biggr),
\]
and its area equals to
\[
|\mathcal{T}_{2n+2}^{*}| = \frac{1}{3(6n+1)(6n+5)}.
\]

The region $\mathcal{T}_{2n+3}^{*} = \mathcal{T}(5,(1,4)^{n},1,3)$ is defined by the system
\begin{equation*}
\begin{cases}
\displaystyle \frac{6n-1}{6n+1} < x\leqslant \frac{6n+5}{6n+7},\quad  \frac{1}{6}(1+x)<y\leqslant \frac{1+(2n+1)x}{6n+5},\\[12pt]
\displaystyle \frac{6n+5}{6n+7} < x\leqslant 1,\quad \frac{1+(2n+3)x}{6n+11}<y\leqslant \frac{1+(2n+1)x}{6n+5}.
\end{cases}
\end{equation*}
It is a quadrangle with vertices
\[
\biggl(\frac{6n-1}{6n+1},\,\frac{2n}{6n+1}\biggr),\quad \biggl(1,\,\frac{2n+2}{6n+5}\biggr),\quad \biggl(1,\,\frac{2n+4}{6n+11}\biggr),\quad \biggl(\frac{6n+5}{6n+7},\,\frac{2n+2}{6n+7}\biggr),
\]
and its area equals to
\[
|\mathcal{T}_{2n+3}^{\,*}| = \frac{24(n+1)}{(6n+1)(6n+5)(6n+7)(6n+11)}.
\]
The region $\mathcal{T}_{2n+3}^{\,**} = \mathcal{T}(5,(1,4)^{n+1})$ is defined by the system
\[
\frac{6n+5}{6n+7} < x\leqslant 1,\quad  \frac{1}{6}(1+x)<y\leqslant \frac{1+(2n+3)x}{6n+11},\\
\]
and is the quadrangle with vertices
\[
\biggl(\frac{6n+5}{6n+7},\,\frac{2n+2}{6n+7}\biggr),\quad \biggl(1,\,\frac{2n+4}{6n+11}\biggr),\quad \biggl(1,\,\frac{1}{3}\biggr).
\]
Its area is equal to
\[
|\mathcal{T}_{2n+3}^{\,**}| = \frac{1}{3(6n+7)(6n+11)}.
\]
One can check that $|\mathcal{T}_{2n+3}^{\,*}|+|\mathcal{T}_{2n+3}^{\,**}| = |\mathcal{T}_{2n+2}^{\,*}|$.

Each region $\mathcal{T}_{2n+3+m}^{\,*} = \mathcal{T}(5,(1,4)^{n},1,3,2^{m})$, $0\leqslant m \leqslant 3n-1$, coincides with $\mathcal{T}_{2n+3}^{\,*}$.

The region $\mathcal{T}_{5n+3}^{\,*} = \mathcal{T}(5,(1,4)^{n},1,3,2^{3n})$ is defined by the inequalities
\begin{equation*}
\begin{cases}
\displaystyle \frac{12n+1}{12n+5} < x\leqslant \frac{6n+5}{6n+7},\quad  \frac{1}{6}(1+x)<y\leqslant \frac{1+(4n+2)x}{12n+7},\\[12pt]
\displaystyle \frac{6n+5}{6n+7} < x\leqslant \frac{6n+2}{6n+3},\quad \frac{1+(2n+3)x}{6n+11}<y\leqslant \frac{1+(4n+2)x}{12n+7},\\[12pt]
\displaystyle \frac{6n+2}{6n+3} < x\leqslant 1,\quad \frac{1+(2n+3)x}{6n+11}<y\leqslant \frac{1+(2n+1)x}{6n+5}
\end{cases}
\end{equation*}
and is the pentagon with the vertices\footnote{It is interesting to note that the region $\mathcal{T}(5,(1,4)^{n},1,3,2^{3n})$ together with the regions $\mathcal{T}(2^{3n},3,1,(4,1)^{n},5)$, $n=1,2,3,\ldots$, that correspond to the tuples obtained from each other by the symmetry, are the single known examples of non-empty pentagons among the regions $\mathcal{T}(\mathbf{k})$. All the rest figures are triangles or quadrangles. It is not known, are there the polygons with six, seven or more vertices among the polygons $\mathcal{T}(\mathbf{k})$ .}
\begin{multline*}
\biggl(\frac{12n+1}{12n+5},\,\frac{4n+1}{12n+5}\biggr),\quad \biggl(\frac{6n+2}{6n+3},\,\frac{1}{3}\biggr),\quad \biggl(1,\,\frac{2n+2}{6n+5}\biggr),\\ \biggl(1,\,\frac{2n+4}{6n+11}\biggr),\quad \biggl(\frac{6n+5}{6n+7},\,\frac{2n+2}{6n+7}\biggr).
\end{multline*}
Its area equals to
\[
|\mathcal{T}_{5n+3}^{\,*}| = \frac{1}{6(2n+1)(6n+5)(12n+5)}+\frac{42n+43}{(6n+5)(6n+7)(6n+11)(12n+5)}.
\]

The region $\mathcal{T}_{5n+3}^{\,\circ} = \mathcal{T}(5,(1,4)^{n},1,3,2^{3n-1},1)$ is defined by the system
\begin{equation*}
\begin{cases}
\displaystyle \frac{6n-1}{6n+1} < x\leqslant \frac{12n+1}{12n+5},\quad  \frac{1}{6}(1+x)<y\leqslant \frac{1+(2n+1)x}{6n+5},\\[12pt]
\displaystyle \frac{12n+1}{12n+5} < x\leqslant \frac{6n+2}{6n+3},\quad \frac{1+(4n+2)x}{12n+7}<y\leqslant \frac{1+(2n+1)x}{6n+5}
\end{cases}
\end{equation*}
and is the triangle with the vertices
\[
\biggl(\frac{6n-1}{6n+1},\,\frac{2n}{6n+1}\biggr),\quad \biggl(\frac{6n+2}{6n+3},\,\frac{1}{3}\biggr),\quad \biggl(\frac{12n+1}{12n+5},\,\frac{4n+1}{12n+5}\biggr).
\]
Its area is equal to
\[
|\mathcal{T}_{5n+3}^{\,\circ}| = \frac{1}{6(2n+1)(6n+1)(12n+5)},
\]
so one can check that $|\mathcal{T}_{5n+3}^{\,*}|+|\mathcal{T}_{5n+3}^{\,\circ}| = |\mathcal{T}_{5n+2}^{\,*}| = |\mathcal{T}_{2n+3}^{\,*}|$.

The region $\mathcal{T}_{5n+4}^{\,*} = \mathcal{T}(5,(1,4)^{n},1,3,2^{3n+1})$ is defined by the inequalities
\begin{equation*}
\begin{cases}
\displaystyle \frac{3n+1}{3n+2} < x\leqslant \frac{6n+5}{6n+7},\quad  \frac{1}{6}(1+x)<y\leqslant \frac{1+(4n+3)x}{12n+10},\\[12pt]
\displaystyle \frac{6n+5}{6n+7} < x\leqslant 1,\quad \frac{1+(2n+3)x}{6n+11}<y\leqslant \frac{1+(4n+3)x}{12n+10}
\end{cases}
\end{equation*}
and is the quadrangle with the vertices
\[
\biggl(\frac{3n+1}{3n+2},\,\frac{2n+1}{6n+4}\biggr),\quad \biggl(1,\,\frac{2n+2}{6n+5}\biggr),\quad \biggl(1,\,\frac{2n+4}{6n+11}\biggr), \quad \biggl(\frac{6n+5}{6n+7},\,\frac{2n+2}{6n+7}\biggr)
\]
and area
\[
|\mathcal{T}_{5n+4}^{\,*}| = \frac{18n+19}{2(3n+2)(6n+5)(6n+7)(6n+11)}.
\]
The region $\mathcal{T}_{5n+4}^{\,\circ} = \mathcal{T}(5,(1,4)^{n},1,3,2^{3n},1)$ is defined by the system
\begin{equation*}
\begin{cases}
\displaystyle \frac{12n+1}{12n+5} < x\leqslant \frac{3n+1}{3n+2},\quad  \frac{1}{6}(1+x)<y\leqslant \frac{1+(4n+2)x}{12n+7},\\[12pt]
\displaystyle \frac{3n+1}{3n+2} < x\leqslant \frac{6n+2}{6n+3},\quad \frac{1+(4n+3)x}{12n+10}<y\leqslant \frac{1+(4n+2)x}{12n+7},\\[12pt]
\displaystyle \frac{6n+2}{6n+3} < x\leqslant 1,\quad \frac{1+(4n+3)x}{12n+10}<y\leqslant \frac{1+(2n+1)x}{6n+5}
\end{cases}
\end{equation*}
and is the quadrangle with the vertices
\[
\biggl(\frac{12n+1}{12n+5},\,\frac{4n+1}{12n+5}\biggr),\quad \biggl(\frac{6n+2}{6n+3},\,\frac{1}{3}\biggr),\quad \biggl(1,\,\frac{2n+2}{6n+5}\biggr),\quad \biggl(\frac{3n+1}{3n+2},\,\frac{2n+1}{6n+4}\biggr)
\]
and area
\[
|\mathcal{T}_{5n+4}^{\,\circ}| = \frac{9n+5}{6(2n+1)(3n+2)(6n+5)(12n+5)}.
\]
One can check the identity $|\mathcal{T}_{5n+4}^{\,*}| + |\mathcal{T}_{5n+4}^{\,\circ}| = |\mathcal{T}_{5n+3}^{\,*}|$

The region $\mathcal{T}_{5n+5}^{\,*} = \mathcal{T}(5,(1,4)^{n},1,3,2^{3n+2})$ is defined by the system
\begin{equation*}
\begin{cases}
\displaystyle \frac{12n+7}{12n+11} < x\leqslant \frac{6n+5}{6n+7},\quad  \frac{1}{6}(1+x)<y\leqslant \frac{1+(4n+4)x}{12n+13},\\[12pt]
\displaystyle \frac{6n+5}{6n+7} < x\leqslant 1,\quad \frac{1+(2n+3)x}{6n+11}<y\leqslant \frac{1+(4n+4)x}{12n+13}
\end{cases}
\end{equation*}
and is the quadrangle with the vertices
\[
\biggl(\frac{12n+7}{12n+11},\frac{4n+3}{12n+11}\biggr),\quad \biggl(1,\frac{4n+5}{12n+13}\biggr),\quad \biggl(1,\frac{2n+4}{6n+11}\biggr), \quad \biggl(\frac{6n+5}{6n+7},\frac{2n+2}{6n+7}\biggr)
\]
and area
\[
|\mathcal{T}_{5n+5}^{\,*}| = \frac{48n+55}{(6n+7)(6n+11)(12n+11)(12n+13)}.
\]
The region $\mathcal{T}_{5n+5}^{\,\circ} = \mathcal{T}(5,(1,4)^{n},1,3,2^{3n+1},1)$ is defined by the system
\begin{equation*}
\begin{cases}
\displaystyle \frac{3n+1}{3n+2} < x\leqslant \frac{12n+7}{12n+11},\quad  \frac{1}{6}(1+x)<y\leqslant \frac{1+(4n+3)x}{12n+10},\\[12pt]
\displaystyle \frac{12n+7}{12n+11} < x\leqslant 1,\quad \frac{1+(4n+4)x}{12n+13}<y\leqslant \frac{1+(4n+3)x}{12n+10}
\end{cases}
\end{equation*}
and is the quadrangle with the vertices
\[
\biggl(\frac{3n+1}{3n+2},\,\frac{2n+1}{6n+4}\biggr),\quad \biggl(1,\,\frac{2n+2}{6n+5}\biggr),\quad \biggl(1,\,\frac{4n+5}{12n+13}\biggr), \quad \biggl(\frac{12n+7}{12n+11},\,\frac{4n+3}{12n+11}\biggr)
\]
and area
\[
|\mathcal{T}_{5n+5}^{\,\circ}| = \frac{3(8n+7)}{2(3n+2)(6n+5)(12n+11)(12n+13)}.
\]
One can check that $|\mathcal{T}_{5n+5}^{\,*}| + |\mathcal{T}_{5n+5}^{\,\circ}| = |\mathcal{T}_{5n+4}^{\,*}|$.

The region $\mathcal{T}_{5n+6}^{\,*} = \mathcal{T}(5,(1,4)^{n},1,3,2^{3n+3})$ is defined by the system of inequalities
\begin{equation*}
\begin{cases}
\displaystyle \frac{6n+5}{6n+7} < x\leqslant 1,\quad \frac{1+(2n+3)x}{6n+11}<y\leqslant \frac{1+(4n+5)x}{12n+16}
\end{cases}
\end{equation*}
and is the triangle with the vertices
\[
\biggl(\frac{6n+5}{6n+7},\,\frac{2n+2}{6n+7}\biggr),\quad \biggl(1,\,\frac{2n+3}{6n+8}\biggr),\quad \biggl(1,\,\frac{2n+4}{6n+11}\biggr)
\]
and area
\[
|\mathcal{T}_{5n+6}^{\,*}| = \frac{1}{2(3n+4)(6n+7)(6n+11)}.
\]
The region $\mathcal{T}_{5n+6}^{\circ} = \mathcal{T}(5,(1,4)^{n},1,3,2^{3n+2},1)$ is defined by the system of inequalities
\begin{equation*}
\begin{cases}
\displaystyle \frac{12n+7}{12n+11} < x\leqslant \frac{6n+5}{6n+7},\quad  \frac{1}{6}(1+x)<y\leqslant \frac{1+(4n+4)x}{12n+13},\\[12pt]
\displaystyle \frac{6n+5}{6n+7} < x\leqslant 1,\quad \frac{1+(4n+5)x}{12n+16}<y\leqslant \frac{1+(4n+4)x}{12n+13}
\end{cases}
\end{equation*}
and is the quadrangle with the vertices
\[
\biggl(\frac{12n+7}{12n+11},\,\frac{4n+3}{12n+11}\biggr),\quad \biggl(1,\,\frac{4n+5}{12n+13}\biggr),\quad \biggl(1,\,\frac{2n+3}{6n+8}\biggr), \quad \biggl(\frac{6n+5}{6n+7},\,\frac{2n+2}{6n+7}\biggr)
\]
and area
\[
|\mathcal{T}_{5n+6}^{\,\circ}| = \frac{3(8n+9)}{2(3n+4)(6n+7)(12n+11)(12n+13)}.
\]
One can check that $|\mathcal{T}_{5n+6}^{\,*}|+|\mathcal{T}_{5n+6}^{\,\circ}| = |\mathcal{T}_{5n+5}^{\,*}|$.

the region $\mathcal{T}_{5n+7}^{\,\circ} = \mathcal{T}(5,(1,4)^{n},1,3,2^{3n+3},1)$ is defined by the system
\begin{equation*}
\begin{cases}
\displaystyle \frac{6n+5}{6n+7} < x\leqslant \frac{6n+8}{6n+9},\quad  \frac{1+(2n+3)x}{6n+11} <y\leqslant \frac{1+(4n+5)x}{12n+16},\\[12pt]
\displaystyle \frac{6n+8}{6n+9} < x\leqslant 1,\quad \frac{1+(4n+6)x}{12n+19}<y\leqslant \frac{1+(4n+5)x}{12n+16}
\end{cases}
\end{equation*}
and is the quadrangle with the vertices
\[
\biggl(\frac{6n+5}{6n+7},\,\frac{2n+2}{6n+7}\biggr),\quad \biggl(1,\,\frac{2n+3}{6n+8}\biggr),\quad \biggl(1,\,\frac{4n+7}{12n+19}\biggr), \quad \biggl(\frac{6n+8}{6n+9},\,\frac{1}{3}\biggr)
\]
and area
\[
|\mathcal{T}_{5n+7}^{\,\circ}| = \frac{9n+13}{6(2n+3)(3n+4)(6n+7)(12n+19)}.
\]
The region $\mathcal{T}_{5n+7}^{\,*} = \mathcal{T}(5,(1,4)^{n},1,3,2^{3n+4})$ is defined by the system
\[
\frac{6n+8}{6n+9} < x\leqslant 1,\quad  \frac{1+(2n+3)x}{6n+11} <y\leqslant \frac{1+(4n+6)x}{12n+19},
\]
and is the triangle with the vertices
\[
\biggl(\frac{6n+8}{6n+9},\,\frac{1}{3}\biggr),\quad \biggl(1,\,\frac{4n+7}{12n+19}\biggr),\quad \biggl(1,\,\frac{2n+4}{6n+11}\biggr)
\]
and area
\[
|\mathcal{T}_{5n+7}^{\,*}| = \frac{1}{6(2n+3)(6n+11)(12n+19)}.
\]
Finally, the change $\mathcal{T}_{5n+7}^{\,*}\rightarrow \mathcal{T}(5,(1,4)^{n},1,3,2^{3n+4},1) = \mathcal{T}_{5n+8}^{\,\circ}$ leads to the additional condition
\[
\frac{1+(4n+7)x}{12n+22} < y\leqslant \frac{1}{3}(1+x),
\]
which does not affect to the form of the region $\mathcal{T}_{5n+7}^{\,*}$. Therefore, $\mathcal{T}_{5n+8}^{\,\circ}$ coincides with $\mathcal{T}_{5n+7}^{\,*}$. In particular,
\[
|\mathcal{T}_{5n+8}^{\,\circ}| = \frac{1}{6(2n+3)(6n+11)(12n+19)}.
\]
Lemma 28 is proved. $\square$
\vspace{0.3cm}

\textsc{Lemma 29.} \emph{The tree $\mathbb{C}$ is shown in the Fig.~11 (to the left). If $n = 3m$, then the subtree $\mathbb{C}_{n}$ has the form shown in the Fig.~11 (to the right). In this case, the areas of the polygons corresponding to degenerate tuples of the subtree $\mathbb{C}_{n}$ are the following:}
\begin{align*}
& |\mathcal{T}_{5m+2}^{\,\circ}| = |\mathcal{T}(1,2^{3m},3,(1,4)^{m-1},1,5)| = \frac{9m+4}{6(2m+1)(3m+1)(6m+1)(12m+7)};\\[6pt]
& |\mathcal{T}_{5m+3}^{\,\circ}| = |\mathcal{T}(1,2^{3m},3,(1,4)^{m},2)| = \frac{3(2m+1)}{2(3m+1)(3m+2)(12m+5)(12m+7)};\\[6pt]
& |\mathcal{T}_{5m+4}^{\,\circ}| = |\mathcal{T}(1,2^{3m},3,(1,4)^{m},1,5)| = \frac{9m+5}{6(2m+1)(3m+2)(6m+5)(12m+5)}.
\end{align*}
\emph{If $n = 3m+1$, then the subset $\mathbb{C}_{n}$ has the form shown in the Fig.~12 (from above). Moreover, the areas of the polygons corresponding to degenerate tuples of the subtree $\mathbb{C}_{n}$ are the following:}
\begin{align*}
& |\mathcal{T}_{5m+3}^{\,\circ}| = |\mathcal{T}(1,2^{3m+1},3,(1,4)^{m-1},1,5)| = \frac{1}{6(2m+1)(6m+5)(12m+7)};\\[6pt]
& |\mathcal{T}_{5m+4}^{\,\circ}| = |\mathcal{T}(1,2^{3m+1},3,(1,4)^{m},2)| = \frac{3(4m+3)}{2(3m+2)(6m+5)(12m+7)(12m+11)};\\[6pt]
& |\mathcal{T}_{5m+5}^{\,\circ}| = |\mathcal{T}(1,2^{3m+1},3,(1,4)^{m},1,5)| = \frac{3(8m+7)}{2(3m+2)(6m+5)(12m+11)(12m+13)};\\[6pt]
& |\mathcal{T}_{5m+6}^{\,\circ}| = |\mathcal{T}(1,2^{3m+1},3,(1,4)^{m+1},2)| = \frac{1}{2(6m+5)(6m+7)(12m+13)}.
\end{align*}
\emph{Finally, if $n = 3m+2$ then the subtree $\mathbb{C}_{n}$ has the form shown in Fig.~12 (below). Moreover, the areas of the polygons corresponding to degenerate tuples of the subtree $\mathbb{C}_{n}$ are the following:}
\begin{align*}
& |\mathcal{T}_{5m+5}^{\,\circ}| = |\mathcal{T}(1,2^{3m+2},3,(1,4)^{m},2)| = \frac{1}{2(6m+5)(6m+7)(12m+11)};\\[6pt]
& |\mathcal{T}_{5m+6}^{\,\circ}| = |\mathcal{T}(1,2^{3m+2},3,(1,4)^{m},1,5)| = \frac{3(8m+9)}{2(3m+4)(6m+7)(12m+11)(12m+13)};\\[6pt]
& |\mathcal{T}_{5m+7}^{\,\circ}| = |\mathcal{T}(1,2^{3m+2},3,(1,4)^{m+1},2)| = \frac{3(4m+5)}{2(3m+4)(6m+7)(12m+13)(12m+17)};\\[6pt]
& |\mathcal{T}_{5m+8}^{\,\circ}| = |\mathcal{T}(1,2^{3m+2},3,(1,4)^{m+1},1,5)| = \frac{1}{6(2m+3)(6m+7)(12m+17)}.
\end{align*}

The proof of Lemma 29 is similar to the proof of Lemma 27 and is based on the formulas (\ref{lab_09}), (\ref{lab_55}) and (\ref{lab_56}).
\vspace{0.3cm}

\begin{center}
\includegraphics{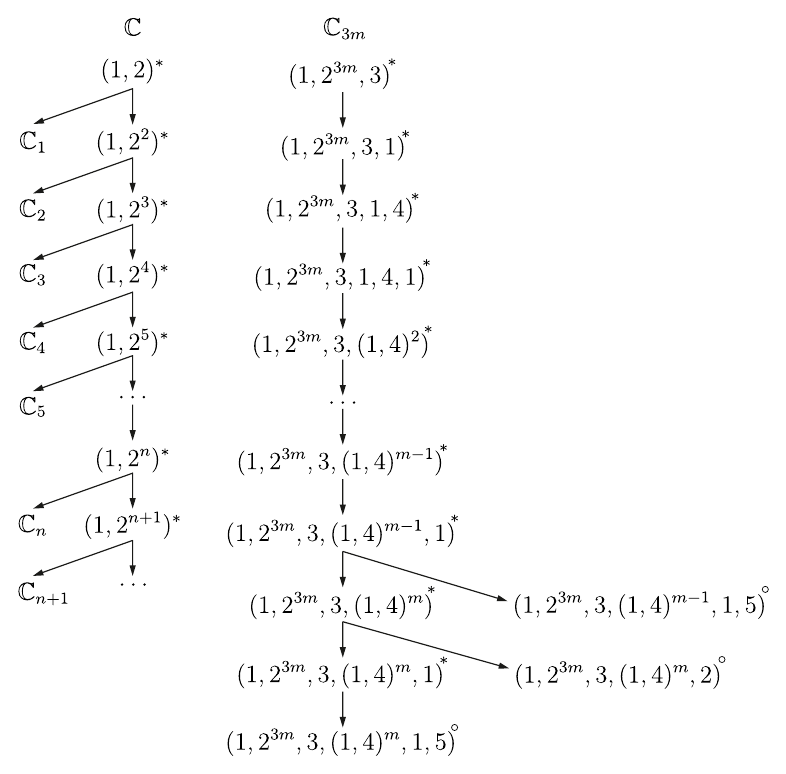}

\fontsize{10}{12pt}\selectfont
\emph{Fig.~11}. The tree $\mathbb{C}$ (to the left) and its subtree $\mathbb{C}_{n}$ for $n = 3m$ (to the right).
\fontsize{12}{15pt}\selectfont
\end{center}

\pagebreak

\begin{landscape}
\begin{center}
\includegraphics{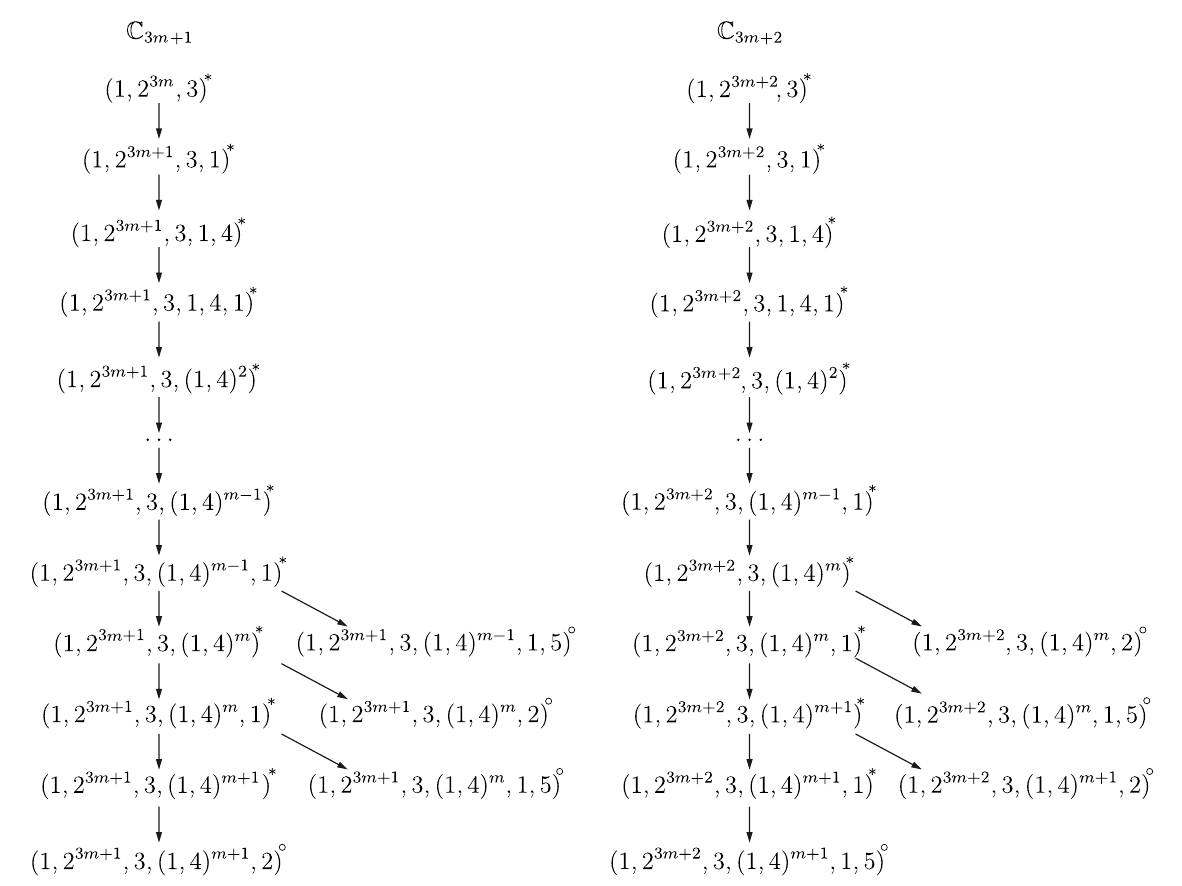}
\vspace{0.2cm}

\fontsize{10}{12pt}\selectfont
\emph{Fig.~12.} The sub-trees $\mathbb{C}_{n}$ for $n = 3m+1$ (above) and $n = 3m+2$ (below).
\fontsize{12}{15pt}\selectfont
\end{center}
\end{landscape}

\pagebreak

\section{Theorem 1: case $\boldsymbol{r\geqslant 8}$}

The above arguments imply that
\begin{equation*}
|\mathcal{A}_{r}^{\circ}| = |\mathcal{A}_{r}^{\circ}(3;0,1)| =
\begin{cases}
4, & \text{if}\quad r\geqslant 9, r\not\equiv 3\pmod{5},\\
6, & \text{if}\quad r\geqslant 8, r\equiv 3\pmod{5}.
\end{cases}
\end{equation*}
In particular, this means that the sum $\Sigma_{0}$ defined in Lemma 18 is bounded by $6$ for any $r\geqslant 8$. Therefore, Lemma 23 yields
\begin{equation}\label{lab_X}
\nu(Q;r,3,0) = \frac{2}{3}\cdot\frac{3}{2}\biggl(\mathfrak{c}_{r} + O\biggl(\frac{\ln{Q}}{Q}\biggr)\biggr) = \mathfrak{c}_{r} + O\biggl(\frac{\ln{Q}}{Q}\biggr),
\end{equation}
where the implied constant is absolute. Hence, for fixed $r = 5m+i$ we get from (\ref{lab_X}):
\[
\nu(r) = \nu(r,3,0) = \lim_{Q\to +\infty}\nu(Q;r,3,0) = \mathfrak{c}_{r} = \sum\limits_{\mathbf{\overline{k}}\in \mathcal{A}_{r}^{\circ}}|\mathcal{T}(\mathbf{\overline{k}})|.
\]

Starting from the set $\mathcal{A}_{6}^{*}$ of tuples constructed in the proof of Theorem 1, we construct by the above procedure the set $\mathcal{A}_{7}^{*}$.
It contains the tuples
\begin{equation}\label{lab_71}
\text{(a)}
\begin{cases}
& (2,(4,1)^{3}),\\
& (2,(4,1)^{2},3,2),\\
& (2,4,1,3,2^{3}),
\end{cases}
\quad
\text{(b)}
\begin{cases}
& (5,(1,4)^{3}),\\
& (5,(1,4)^{2},1,3),\\
& (5,1,4,1,3,2^{2}),\\
& (5,1,3,2^{4}),
\end{cases}
\quad
\text{(c)}
\begin{cases}
& (1,2^{6}),\\
& (1,2^{5},3),\\
& (1,2^{4},3,1),\\
& (1,2^{3},3,1,4),\\
& (1,2^{2},3,1,4,1).
\end{cases}
\end{equation}

We refer the tuples corresponding to the leaves of the trees $\mathbb{A}$, $\mathbb{B}$, and $\mathbb{C}$ as the tuples of types $A$, $B$ and $C$ correspondingly.
Then one can check that the tuples (\ref{lab_71} a) are $A$-type tuples, (\ref{lab_71} b) are $B$-type tuples, and (\ref{lab_71} c) are $C$-type tuples.
All the tuples generated by the sets (\ref{lab_71} a-c) have the corresponding types.

Therefore, to find the proportion $\nu(r)$, $r\geqslant 8$, it is sufficient to sum all the areas of the regions $\mathcal{T}(\mathbf{k})$ that correspond to the leaves of the trees $\mathbb{A}$, $\mathbb{B}$, $\mathbb{C}$ that are located on level $r$. In the trees, all the degenerate tuples are labeled by the sign <<$\circ$>>. Therefore, the calculation of the proportion $\nu(r)$ follows to the below scheme. Giver $r$, we pick out all the tuples of length $r$ of each set from Lemmas 27-29 and the sum the corresponding areas.

For example, suppose that $r = 5m$, where $m\geqslant 2$. Then there exists the single degenerate $A$-type tuple of length $r$, namely
$(2,(4,1)^{m-1},3,2^{3(m-1)+2},1)$. By Lemma 27, the area of $\mathcal{T}_{5m}^{\,\circ} = \mathcal{T}_{5(m-1)+5}^{\,\circ}$ equals to
\[
\frac{1}{2(6(m-1)+5)(6(m-1)+7)(12(m-1)+11)} = \frac{1}{2(6m-1)(6m+1)(12m-1)}.
\]
Similarly, there exists a unique $B$-type tuple of length $r$, namely
\[
\bigl(5,(1,4)^{m-1},1,3,2^{3(m-1)+1},1\bigr),
\]
and we have
\begin{multline*}
|\mathcal{T}_{5m}^{\,\circ}| = |\mathcal{T}_{5(m-1)+5}^{\,\circ}| = \frac{3(8(m-1)+7)}{2(3(m-1)+2)(6(m-1)+5)(12(m-1)+11)(12(m-1)+13)} = \\[6pt]
 = \frac{3(8m-1)}{2(3m-1)(6m-1)(12m-1)(12m+1)}.
\end{multline*}

Further, there are two degenerate $C$-type tuples, namely
\[
\bigl(1,2^{3(m-1)+1},3,(1,4)^{m-1},1,5\bigr),\quad \bigl(1,2^{3(m-1)+2},3,(1,4)^{m-1},2\bigr),
\]
and the corresponding areas are equal to
\[
\frac{3(8m-1)}{2(3m-1)(6m-1)(12m-1)(12m+1)}\quad\text{and}\quad \frac{1}{2(6m-1)(6m+1)(12m-1)}.
\]
therefore,
\begin{multline*}
\nu(5m) = \frac{6(8m-1)}{2(3m-1)(6m-1)(12m-1)(12m+1)} + \frac{2}{(6m-1)(6m+1)(12m-1)} =\\[6pt]
= P_{01}(m) + P_{02}(m).
\end{multline*}
Further, if $r = 5m+1$, $m\geqslant 2$, then we have the following degenerate tuples of length $r$
\begin{multline*}
\bigl(2,(4,1)^{m},3,2^{3m-2},1\bigr)\quad (A\text{ - type}),\quad
\bigl(5,(1,4)^{m-1},1,3,2^{3(m-1)+2},1\bigr)\quad (B\text{ - type}),\\
\bigl(1,2^{3(m-1)+1},3,(1,4)^{m},2\bigr)\quad (C\text{ - type}),\quad
\bigl(1,2^{3(m-1)+2},3,(1,4)^{m-1},1,5\bigr)\quad (C\text{ - type}).
\end{multline*}
The corresponding areas are equal to
\begin{multline*}
\frac{1}{2(6m-1)(6m+1)(12m+1)},\quad \frac{3(8m+1)}{2(3m+1)(6m+1)(12m-1)(12m+1)}, \\[6pt]
\frac{1}{2(6m-1)(6m+1)(12m+1)},\quad \frac{3(8m+1)}{2(3m+1)(6m+1)(12m-1)(12m+1)}.
\end{multline*}
Hence,
\begin{multline*}
\nu(5m+1) = \frac{6(8m+1)}{2(3m+1)(6m+1)(12m-1)(12m+1)} + \frac{2}{(6m-1)(6m+1)(12m+1)} = \\[6pt]
= P_{11}(m)+P_{12}(m).
\end{multline*}
If $r = 5m+2$, $m\geqslant 2$, then the degenerate tuples of length $r$ are:
\begin{multline*}
\bigl(2,(4,1)^{m},3,2^{3m-1},1\bigr)\quad (A\text{ - type}),\quad
\bigl(5,(1,4)^{m-1},1,3,2^{3(m-1)+3},1\bigr)\quad (B\text{ - type}),\\
\bigl(1,2^{3m},3,(1,4)^{m-1},1,5\bigr)\quad (C\text{ - type}),\quad
\bigl(1,2^{3(m-1)+2},3,(1,4)^{m},2\bigr)\quad (C\text{ - type}).
\end{multline*}
The corresponding areas are equal to
\begin{multline*}
\frac{3(4m+1)}{2(3m+1)(6m+1)(12m+1)(12m+5)},\quad \frac{9m+4}{6(2m+1)(3m+1)(6m+1)(12m+7)}, \\[6pt]
\frac{9m+4}{6(2m+1)(3m+1)(6m+1)(12m+7)},\quad \frac{3(4m+1)}{2(3m+1)(6m+1)(12m+1)(12m+5)}.
\end{multline*}
Thus,
\begin{multline*}
\nu(5m+2) = \frac{6(4m+1)}{2(3m+1)(6m+1)(12m+1)(12m+5)} + \frac{2(9m+4)}{6(2m+1)(3m+1)(6m+1)(12m+7)} = \\[6pt]
= P_{21}(m)+P_{22}(m).
\end{multline*}
In the case$r = 5m+3$, $m\geqslant 1$, the degenerate tuples of length $r$ are:
\begin{multline*}
\bigl(2,(4,1)^{m},3,2^{3m},1\bigr)\quad (A\text{ - type}),\quad
\bigl(5,(1,4)^{m},1,3,2^{3m-1},1\bigr)\quad (B\text{ - type}),\\
\bigl(5,(1,4)^{m-1},1,3,2^{3(m-1)+4},1\bigr)\quad (B\text{ - type}),\quad \bigl(1,2^{3m},3,(1,4)^{m},2\bigr)\quad (\text{C - type}),\\
\bigl(1,2^{3m+1},3,(1,4)^{m-1},1,5\bigr)\quad (C\text{ - type}),\quad \bigl(1,2^{3(m-1)+2},3,(1,4)^{m},1,5\bigr)\quad (\text{C - type}).
\end{multline*}
The corresponding areas are equal to
\begin{multline*}
\frac{3(2m+1)}{2(3m+1)(3m+2)(12m+5)(12m+7)},\quad \frac{1}{6(2m+1)(6m+1)(12m+5)}, \\[6pt]
\frac{1}{6(2m+1)(6m+5)(12m+7)},\quad \frac{3(2m+1)}{2(3m+1)(3m+2)(12m+5)(12m+7)},\\[6pt]
\frac{1}{6(2m+1)(6m+5)(12m+7)},\quad \frac{1}{6(2m+1)(6m+1)(12m+5)}.
\end{multline*}
Hence,
\begin{multline*}
\nu(5m+3) = \frac{6(2m+1)}{2(3m+1)(3m+2)(12m+5)(12m+7)} + \frac{2}{3(2m+1)(6m+1)(12m+5)} + \\[6pt]
+ \frac{2}{3(2m+1)(6m+5)(12m+7)} = P_{31}(m)+P_{32}(m)+P_{33}(m).
\end{multline*}
Finally, if $r = 5m+4$, $m\geqslant 1$, then the degenerate tuples of length $r$ are:
\begin{multline*}
\bigl(2,(4,1)^{m},3,2^{3m+1},1\bigr)\quad (A\text{ - type}),\quad
\bigl(5,(1,4)^{m},1,3,2^{3m},1\bigr)\quad (B\text{ - type}),\\
\bigl(1,2^{3m},3,(1,4)^{m},1,5\bigr)\quad (C\text{ - type}),\quad
\bigl(1,2^{3m+1},3,(1,4)^{m},2\bigr)\quad (C\text{ - type}).
\end{multline*}
The corresponding areas are equal to
\begin{multline*}
\frac{3(4m+3)}{2(3m+2)(6m+5)(12m+7)(12m+11)},\quad \frac{9m+5}{6(2m+1)(3m+2)(6m+5)(12m+5)}, \\[6pt]
\frac{9m+5}{6(2m+1)(3m+2)(6m+5)(12m+5)},\quad \frac{3(4m+3)}{2(3m+2)(6m+5)(12m+7)(12m+11)}.
\end{multline*}
Thus we conclude that
\begin{multline*}
\nu(5m+4) = \\
=\frac{6(4m+3)}{2(3m+2)(6m+5)(12m+7)(12m+11)} + \frac{2(9m+5)}{(6m+3)(6m+4)(6m+5)(12m+5)} = \\[6pt]
=P_{41}(m)+P_{42}(m).
\end{multline*}
Finally, we note that $P_{ij}(m)\asymp m^{-3}\asymp r^{-3}$ for $r = 5m+i$. Hence, the formula (\ref{lab_X}) is asymptotic if $r = o\bigl((Q/\ln{Q})^{1/3}\bigr)$.
The last condition is the main reason why we trace carefully the dependence on the parameter $r$ in the remainders in Lemmas 18, 21 and 23.
Now Theorem 1 is proved completely. $\square$

\section{The formula of three authors}

In this section, we establish the correspondence between the main result of \cite{Cobeli_Vajaitu_Zaharescu_2012} and the above formula for the limit proportion $\nu(r;D,c)$ for the case of fixed $D$ and $r$. To do this, we recall some definitions introduced in \cite{Cobeli_Vajaitu_Zaharescu_2012}.

Let $r\geqslant 1$, $Q\geqslant \max{(3,r)}$, and let $q_{0}, q_{1},\ldots, q_{r+1}$ be the tuple of the consecutive denominators of Farey fractions from $\Phi_{Q}$. Then the numbers $k_{i}$ defined by the relations
\[
k_{i} = \frac{q_{i-1}+q_{i+1}}{q_{i}} = \biggl[\frac{q_{i-1}+Q}{q_{i}}\biggr],\quad i = 1,2,\ldots, r,
\]
are integer (see \S 2). In this case, following to \cite{Cobeli_Vajaitu_Zaharescu_2012}, we say that the tuple $\mathbf{k}_{r} = (k_{1},\ldots, k_{r})$ is generated by the triple $Q, q_{0}, q_{1}$. The tuple $\mathbf{k}_{r}$ generated by some triple $Q, q_{0}, q_{1}$ such that
\begin{equation}\label{lab_74}
1\leqslant q_{0},q_{1}\leqslant Q,\quad q_{0}+q_{1}>Q,\quad \text{GCD}(q_{0},q_{1}) = 1,
\end{equation}
is called an \emph{admissible}.

Further, fix the number $D\geqslant 2$ and the tuple $\overline{\mathbf{c}} = (c_{0},\ldots,c_{r+1})$ of integers with the conditions $0\leqslant c_{i}\leqslant D-1$, $i=0,1,\ldots,r+1$.
We say that $\mathbf{k}_{r}$ belongs to the set $\mathcal{K}^{r}(\overline{\mathbf{c}},D)$ if $\mathbf{k}_{r}$ is admissible and if there exists a triple $Q, q_{0}, q_{1}$ generating $\mathbf{k}_{r}$ and such that the corresponding tuple of denominators satisfies to the following congruences:
\begin{equation}\label{lab_75}
q_{i}\equiv c_{i}\pmod{D},\quad i = 0,1,\ldots, r+1.
\end{equation}

Next, if $\ell\geqslant 3$ and $\overline{\mathbf{v}} = (v_{0},\ldots,v_{\ell-1})$ is a fixed tuple with the conditions $0\leqslant v_{i}\leqslant D-1$, $0\leqslant i\leqslant \ell-1$, and if $N(Q;\ell,D,\overline{\mathbf{v}})$ denotes the number of tuples of consecutive Farey fractions
\[
\frac{a_{0}}{q_{0}}<\frac{a_{1}}{q_{1}} < \ldots < \frac{a_{\ell-2}}{q_{\ell-2}}<\frac{a_{\ell-1}}{q_{\ell-1}},\quad \frac{a_{i}}{q_{i}}\in \Phi_{Q},\quad q_{i}\equiv v_{i}\pmod{D},
\]
then the following formula for the limit proportion
\[
\varrho(\ell,D,\overline{\mathbf{v}}) = \lim_{Q\to +\infty}\frac{N(Q;\ell,D,\mathbf{\overline{v}})}{N(Q)}
\]
is given in \cite[Corollary of Theorem 2]{Cobeli_Vajaitu_Zaharescu_2012}:
\begin{equation}\label{lab_76}
\varrho(\ell,D,\overline{\mathbf{v}}) = \frac{2}{D^{2}}\prod\limits_{p|D}\biggl(1-\frac{1}{p^{2}}\biggr)^{-1}\sum\limits_{\mathbf{\overline{k}_{\ell-2}}\in \mathcal{K}^{\ell-2}(\overline{\mathbf{v}},D)}|\mathcal{T}(\mathbf{\overline{k}_{\ell-2}})|
\end{equation}
Now we express the limit proportion
\[
\nu(r,D,c_{0}) = \lim_{Q\to +\infty}\nu(Q;r,D,c_{0})  = \lim_{Q\to +\infty}\frac{N(Q;r,D,c_{0})}{N(Q)}
\]
in terms of $\varrho(\cdot)$.

Obviously, one has
\[
N(Q;r,D,c_{0}) = \sum\limits_{\mathbf{\overline{c}}}N(Q;r+2,D,\mathbf{\overline{c}}),
\]
where the sign $\sum\limits_{\mathbf{\overline{c}}}$ means the summation over the set of integer tuples $\mathbf{\overline{c}} = (c_{0},\ldots,c_{r+1})$ obeying the conditions $0\leqslant c_{i}\leqslant D-1$, $c_{i}\not\equiv c_{0}$, $i = 1,2,\ldots, r$, $c_{r+1}= c_{0}$ and, moreover, $\text{GCD}(c_{1},\Delta)=1$ where $\Delta = \text{GCD}(c_{0},D)$. Therefore,
\[
\nu(Q;r,D,c_{0}) = \frac{N(Q)}{N(Q;D,c_{0})}\cdot \frac{N(Q;r,D,c_{0})}{N(Q)} = \frac{N(Q)}{N(Q;D,c_{0})}\sum\limits_{\mathbf{\overline{c}}}\frac{N(Q;r+2,D,\mathbf{\overline{c}})}{N(Q)}.
\]
By Lemma 22,
\[
\lim_{Q\to +\infty}\frac{N(Q)}{N(Q;D,c_{0})} = \frac{D\Delta}{\varphi(\Delta)}\prod\limits_{p|D}\biggl(1-\frac{1}{p^{2}}\biggr).
\]
Using (\ref{lab_76}), we find
\begin{multline}\label{lab_77}
\nu(r,D,c_{0}) = \frac{D\Delta}{\varphi(\Delta)}\prod\limits_{p|D}\biggl(1-\frac{1}{p^{2}}\biggr)\sum\limits_{\mathbf{\overline{c}}}\varrho(r+2,D,\mathbf{\overline{c}}) = \\
= \frac{2\Delta}{D\varphi(\Delta)}\sum\limits_{\mathbf{\overline{c}}}\sum\limits_{\mathbf{k_{\ell-2}}\in \mathcal{K}^{\ell-2}(\overline{\mathbf{c}},D)}|\mathcal{T}(\mathbf{k_{\ell-2}})|
\end{multline}
At the same time, Lemma 23 implies that
\begin{equation}\label{lab_78}
\nu(r,D,c_{0}) = \frac{2\Delta}{D\varphi(\Delta)}\sum\limits_{\substack{c_{1} = 0 \\ \text{CDG}(c_{1},\Delta)=1}}^{D-1}\sum\limits_{\mathbf{k_{r}}\in \mathcal{A}_{r}^{\circ}(D,c_{0},c_{1})}|\mathcal{T}(\mathbf{k_{r}})|.
\end{equation}
We show that the sums in (\ref{lab_77}) and (\ref{lab_78}) are equal.

First, we note that any tuple $\mathbf{k_{r}}$ from (\ref{lab_78}) is admissible. It is sufficient to prove that there exists at leas one triple $Q, q_{0}, q_{1}$ satisfying to (\ref{lab_74}) and generating $\mathbf{k_{r}}$. Suppose that $0\leqslant c_{1}\leqslant D-1$ is fixed, $\text{GCD}(c_{1},\Delta)=1$ and let $\mathbf{k_{r}}\in \mathcal{A}_{r}(D,c_{0},c_{1})$.

Since $|\mathcal{T}(\mathbf{k_{r}})|>0$ then there exists a square $\mathcal{P}$ of a small size $\varepsilon = \varepsilon(\mathbf{k_{r}})>0$ lying entirely in $\mathcal{T}(\mathbf{k_{r}})$. In what follows, we may assume that $0<\varepsilon<0.1$. Next, define
\[
Q_{0} = \biggl(\frac{\pi D}{\varepsilon}\biggr)^{2}\ln{\frac{\pi D}{\varepsilon}}
\]
and suppose that $Q>Q_{0}$. Finally, suppose that the square $Q\cdot \mathcal{P}$ has the form
\[
X<x\leqslant X+H,\quad Y<y\leqslant Y+H,\quad H = \varepsilon Q, \quad \text{and take}\quad B = \min{(X,Y)}+H.
\]

Then the number $n$ of the primitive points equals to
\begin{equation}\label{lab_79}
n = \sum\limits_{\substack{X<q_{0}\leqslant X+H \\ q_{0}\equiv c_{0}\;(\mmod D)}}\sum\limits_{\substack{Y<q_{1}\leqslant Y+H \\ q_{1}\equiv c_{1}\;(\mmod D)}}\sum\limits_{\delta|\text{GCD}(q_{0},q_{1)}}\mu(\delta)
\end{equation}
Setting $B = \min{(X,Y)}+H$, we note that $\delta\leqslant B$ for any $\delta$ in (\ref{lab_79}). Thus we find
\[
n = \sum\limits_{\delta\leqslant \Delta}\mu(\delta)\sum\limits_{\substack{X<q_{0}\leqslant X+H \\ q_{0}\equiv c_{0}\;(\mmod D) \\ q_{0}\equiv 0\;(\mmod \delta)}}1\;\sum\limits_{\substack{Y<q_{1}\leqslant Y+H \\ q_{1}\equiv c_{1}\;(\mmod D) \\ q_{1}\equiv 0\;(\mmod \delta) }}1,
\]
In view of the condition $\text{GCD}(c_{1},\Delta) = \text{GCD}(c_{0},c_{1},D)=1$, the non-zero contribution to the last sum comes from $\text{GCD}(\delta,D)$ only. Therefore,
\[
n = \sum\limits_{\substack{\delta\leqslant B \\ \text{GCD}(\delta,D)=1}}\mu(\delta)\biggl(\;\sum\limits_{\substack{X<q_{0}\leqslant X+H \\ q_{0}\equiv \sigma_{0}\pmod{D\delta}}}1\biggr)
\biggl(\;\sum\limits_{\substack{Y<q_{1}\leqslant Y+H \\ q_{1}\equiv \sigma_{1}\pmod{D\delta}}}1\biggr)
\]
where $\sigma_{j}\equiv c_{j}\delta\overline{\delta}\pmod{D\delta}$, $\delta\overline{\delta}\equiv 1\pmod{D}$. Replacing both the above sums by $H/(D\delta)+\theta_{j}$, $|\theta_{j}|\leqslant 1$, we get
\[
n = \sum\limits_{\substack{\delta\leqslant B \\ \text{GCD}(\delta,D)=1}}\biggl(\Bigl(\frac{H}{D}\Bigr)^{2}\frac{\mu(\delta)}{\delta^{2}} + \theta\Bigl(\frac{2H}{D\delta}+1\Bigr)\biggr) = \kappa_{D}\biggl(\frac{H}{D}\biggr)^{2} + R,
\]
where
\[
|R|\leqslant \biggl(\frac{H}{D}\biggr)^{2}\sum\limits_{\delta>D}\frac{1}{\delta^{2}} + \frac{2H}{D}\sum\limits_{\delta\leqslant B}\frac{1}{\delta} + \sum\limits_{\delta\leqslant B}1 <
\frac{2}{B}\biggl(\frac{H}{D}\biggr)^{2} + \frac{2H}{D}(\ln{B}+1)+ B.
\]
Since $H = Q\varepsilon < B\leqslant Q$ and $D\geqslant 2$, we obtain
\[
|R|<\frac{2H}{D^{2}}+\frac{2H}{D}(\ln{Q}+1)+Q\leqslant \varepsilon Q(\ln{Q}+1.5)+Q < Q\ln{Q}.
\]
Hence,
\[
n>\kappa_{D}\biggl(\frac{H}{D}\biggr)^{2}_Q\ln{Q} > 6\biggl(\frac{\varepsilon Q}{\pi D}\biggr)^{2}- Q\ln{Q} = Q(\ln{Q})\biggl(\frac{6Q}{\ln{Q}}\biggl(\frac{\varepsilon }{\pi D}\biggr)^{2}-1\biggr).
\]
Obviously,
\[
\frac{Q}{\ln{Q}} > \frac{Q_{0}}{\ln{Q_{0}}} = \biggl(\frac{\pi D}{\varepsilon}\biggr)^{2}\,\frac{\ln{\frac{\displaystyle \pi D}{\displaystyle \varepsilon}}}{2\ln{\frac{\displaystyle \pi D}{\displaystyle \varepsilon}}+\ln{\ln{\frac{\displaystyle \pi D}{\displaystyle \varepsilon}}}} > \frac{1}{3}\biggl(\frac{\pi D}{\varepsilon}\biggr)^{2}.
\]
Therefore,
\[
n>Q(\ln{Q})\biggl(6\biggl(\frac{\varepsilon }{\pi D}\biggr)^{2}\cdot\frac{1}{3}\biggl(\frac{\pi D}{\varepsilon}\biggr)^{2}-1 \biggr)>Q\ln{Q}>0.
\]
Thus, for any $Q>Q_{0}$, there exists at least one triple $Q, q_{0}, q_{1}$ generating $\mathbf{\overline{k}}_{r}$. This means that the tuple $\mathbf{\overline{k}}_{r}$ is admissible.

Next, if $\mathbf{\overline{k}}_{r}\in \mathcal{A}_{r}^{\circ}(D,c_{0},c_{1})$ is generated by $Q, q_{0}, q_{1}$, then the denominators $q_{2}, q_{3},\ldots, q_{r+1}$ are connected with
the components of $\mathbf{\overline{k}}_{r}$ by the relations (\ref{lab_05}) and (\ref{lab_12}). Since $q_{0}\equiv c_{0}\pmod{D}$ and $q_{1}\equiv c_{1}\pmod{D}$ then, passing from (\ref{lab_12}) to the congruences we get
\[
q_{i+1}\equiv c_{1}\mathbb{K}_{i}(k_{1},\ldots,k_{i}) - c_{0}\mathbb{K}_{i-1}(k_{2},\ldots,k_{i})\pmod{D},\quad i = 1,2,\ldots r.
\]
Hence, if we define the integers $c_{i}$ by the conditions
\[
c_{i+1}\equiv c_{1}\mathbb{K}_{i}(k_{1},\ldots,k_{i}) - c_{0}\mathbb{K}_{i-1}(k_{2},\ldots,k_{i})\pmod{D},\quad 0\leqslant c_{i+1}\leqslant D-1,\quad i = 1,2,\ldots r,
\]
we get $q_{i}\equiv c_{i}\pmod{D}$ for any $i$. Hence, $\mathbf{\overline{k}}_{r}\in \mathcal{K}^{r}(\mathbf{\overline{c}},D)$ where $\mathbf{\overline{c}}  = (c_{0},\ldots,c_{r+1})$.

Since $\text{GCD}(c_{1},\Delta)=1$, $c_{i}\not\equiv c_{0}$ for $i=1,\ldots,r$ and $c_{r+1} = c_{0}$, the tuple $\mathbf{\overline{c}}$ will appear in the sum (\ref{lab_78}).
This means that any tuple from (\ref{lab_77}) will appear in the sum (\ref{lab_78}).

Conversely, let $\mathbf{\overline{c}}$ and $\mathbf{\overline{k}}_{r}\in \mathcal{K}^{r}(\mathbf{\overline{c}},D)$ are the tuples from (\ref{lab_78}), and let $|\mathcal{T}(\mathbf{\overline{k}}_{r})|>0$. Then the tuple $\mathbf{\overline{k}}_{r}$ will appear in (\ref{lab_77}).

Indeed, let $Q, q_{0}, q_{1}$ be the triple that generates $\mathbf{\overline{k}}_{r}$. Then the corresponding sequence $q_{2}, q_{3},\ldots, q_{r+1}$ of denominators is connected with the components of $\mathbf{\overline{k}}_{r}$ by (\ref{lab_05}) and (\ref{lab_12}). If $\mathbf{\overline{c}} = (c_{0},\ldots,c_{r+1})$ then, in view of the definition of $\mathcal{K}^{r}(\mathbf{\overline{c}},D)$,
\begin{align*}
q_{i+1}\equiv c_{i+1}\equiv c_{1}\mathbb{K}_{i}(k_{1},\ldots,k_{i})-c_{0}\mathbb{K}_{i-1}(k_{2},\ldots,k_{i})\not\equiv c_{0}\pmod{D}, \quad i = 1,\ldots, r-1,\\
q_{r+1}\equiv c_{r+1}\equiv c_{1}\mathbb{K}_{r}(k_{1},\ldots,k_{r})-c_{0}\mathbb{K}_{r-1}(k_{2},\ldots,k_{r})\equiv c_{0}\pmod{D}.
\end{align*}
These relations imply that $\mathbf{\overline{k}}_{r}\in \mathcal{A}_{r}^{\circ}(D,c_{0},c_{1})$.

Finally, it is not difficult to check that the multiplicities of a given tuple $\mathbf{\overline{k}}_{r}$ in the sums (\ref{lab_77}), (\ref{lab_78}) are equal to the number of residues $c_{1}\pmod{D}$, $\text{GCD}(c_{1},\Delta)$ that satisfy the conditions
\begin{equation*}
\begin{cases}
c_{1}\mathbb{K}_{i}(k_{1},\ldots,k_{i})-c_{0}\mathbb{K}_{i-1}(k_{2},\ldots,k_{i})\not\equiv c_{0}\pmod{D}, \quad i = 1,\ldots, r-1,\\
c_{1}\mathbb{K}_{r}(k_{1},\ldots,k_{r})-c_{0}\mathbb{K}_{r-1}(k_{2},\ldots,k_{r})\equiv c_{0}\pmod{D}.
\end{cases}
\end{equation*}
Thus, the sums (\ref{lab_77}), (\ref{lab_78}) are equal.
\vspace{0.3cm}

\section{Theorem 2: the case $\boldsymbol{D = 2, c_{0}=0}$}

We give only a sketch of the proof of Theorem 2 since the derivation of the formulas for $\nu(r;2,0)$ follow to the same steps as for $\nu(r;3,0)$ for $r\geqslant 8$ (see Theorem 1). At the same time, in the case $D = 2$, the computations are not so complicated as for the case $D = 3$.

We will construct the sets $\mathcal{A}_{r}^{\circ} = \mathcal{A}_{r}^{\circ}(2,2,0)$ and $\mathcal{A}_{r}^{*} = \mathcal{A}_{r}^{*}(2,2,0)$, $r \geqslant 1$ consequently.

\textsc{Case} $r = 1$. Here $\mathcal{A}_{r}^{*}$ consists of the tuples $(k_{1})$, $k_{1}\equiv 1\pmod{2}$, and $\mathcal{A}_{r}^{\circ}$ consists of the tuples $(k_{1})$, $k_{1}\equiv 0\pmod{2}$.

\textsc{Case} $r = 2$. The set $\mathcal{A}_{r}^{*}$ consists of the tuples  $(k_{1},k_{2})$, where $k_{1}\equiv 1\pmod{2}$, $k_{2}\equiv 0\pmod{2}$. The corresponding regions $\mathcal{T}(k_{1},k_{2})$ are non-empty for the pairs $(3,2)$ and $(1,k_{2})$, $k_{2}\equiv 0\pmod{2}$, $k_{2}\geqslant 2$ only. The set $\mathcal{A}_{r}^{\circ}$ consists of the tuples $(k_{1},k_{2})$, where $k_{1}, k_{2}\equiv 1\pmod{2}$. The non-empty regions correspond to the tuples $(k_{1},1)$, $k_{1}\equiv 1\pmod{2}$, $k_{1}\geqslant 3$, and to the pairs $(1,k_{2})$, $k_{2}\equiv 1\pmod{2}$, $k_{2}\geqslant 3$.

\textsc{Case} $r = 3$. The set $\mathcal{A}_{r}^{*}$ is finite and consists of four tuples
\[
(1,2,2),\quad (1,2,4),\quad (1,4,2),\quad (3,2,2).
\]
The set $\mathcal{A}_{r}^{\circ}$ is infinite and consists of the tuples
\[
(1,k_{2},1),\quad k_{2}\equiv 0\pmod{2},\quad k_{2}\geqslant 4,\quad (1,2,3),\quad (3,2,1).
\]

\textsc{Case} $r=4$. The set $\mathcal{A}_{r}^{\circ}$ consists of
\[
(1,2,2,3),\quad (1,2,4,1), \quad (1,4,2,1),\quad (3,2,2,1).
\]
Using Corollary 1 of Lemma 10, we get
\[
\nu(2;4,0) = 4\bigl(|\mathcal{T}(1,2,2,3)| + |\mathcal{T}(1,2,4,1)|\bigr) = 4\biggl(\frac{2}{315}+\frac{1}{210}\biggr) = \frac{2}{45}.
\]
The set $\mathcal{A}_{r}^{*}$ contains only the tuples $(1,2,2,2)$ and $(3,2,2,2)$.

For any $r\geqslant 5$, the set $\mathcal{A}_{r}^{*}$ consists of the tuples $(1,2^{r-1})$, $(3,2^{r-1})$, while the set $\mathcal{A}_{r}^{\circ}$ consists of the tuples $(1,2^{r-2},3)$ and $(3,2^{r-2},1)$. Correspondingly, the region $\mathcal{T}(1,2^{r-2},3)$ is defined by the inequalities
\begin{equation*}
\begin{cases}
\displaystyle \frac{1}{2r+1} < x\leqslant \frac{1}{2r-1},\quad 1-x<y\leqslant \frac{1}{2}(1+(2r-1)x),\\[12pt]
\displaystyle \frac{1}{2r-1} < x\leqslant \frac{1}{2r-3},\quad \frac{1}{2}(1+(2r-3)x) < y \leqslant 1,
\end{cases}
\end{equation*}
and hence is the quadrangle with the vertices
\[
\biggl(\frac{1}{2r+1},\,\frac{2r}{2r+1}\biggr),\quad \biggl(\frac{1}{2r-1},\,1\biggr),\quad \biggl(\frac{1}{2r-3},\,1\biggr),
\quad \biggl(\frac{1}{2r-1},\,\frac{2r-2}{2r-1}\biggr)
\]
and area
\[
\frac{2}{(2r-3)(2r-1)(2r+1)}.
\]
Similarly, the region $\mathcal{T}(3,2^{r-2},1)$ is defined by the inequalities
\begin{equation*}
\begin{cases}
\displaystyle \frac{2r-5}{2r-3} < x\leqslant \frac{2r-3}{2r-1},\quad \frac{1}{4}(1+x)<y\leqslant \frac{1+(r-1)x}{2r-1},\\[12pt]
\displaystyle \frac{2r-3}{2r-1} < x\leqslant 1,\quad \frac{1+rx}{2r+1}<y\leqslant \frac{1+(r-1)x}{2r-1}
\end{cases}
\end{equation*}
and hence is the quadrangle with the vertices
\[
\biggl(\frac{2r-5}{2r-3},\,\frac{r-2}{2r-3}\biggr),\quad \biggl(1,\,\frac{r}{2r-1}\biggr),\quad \biggl(1,\,\frac{r+1}{2r+1}\biggr),
\quad \biggl(\frac{2r-3}{2r-1},\,\frac{r-1}{2r-1}\biggr)
\]
and area
\[
\frac{2}{(2r-3)(2r-1)(2r+1)}.
\]
Thus the assertion follows. $\square$

\section{Final remarks}

It is natural to ask whether is it possible to derive an explicit formulas for the proportions $\nu(r;D,0)$ for $D>3$? However, even in the case of prime $D=5$, the computations
increase significantly. The corresponding trees of tuples have a much more complicated structure than the above trees $\mathbb{A}$, $\mathbb{B}$ and $\mathbb{C}$.
Thus, in some cases the <<moment of stopping>> (that is, the moment when the tuple in the tree produces the degenerate descendant only) comes very late, and one has to deal with the continuants of more than one hundred variables. It is interesting to find the shorter way to obtain the answer in this problem.

At the same time, it is easy to prove that, for any $D$ and $c_{0}$, all the proportions $\nu(r;D,c_{0})$ are rational for $r>2D$.

Indeed, all linear functions $f_{i}(x;\mathbf{k}_{i})$, $g_{i}(x;\mathbf{k}_{i})$ have rational coefficients. Therefore, all vertices of any non-empty polygon $\mathcal{T}(\mathbf{k})$ have rational coordinates. Hence, the areas of such polygons are also rational.

Suppose that the set $\mathcal{A}_{r}^{\circ}(D,c_{0},c_{1})$ is infinite for some $r$ and $c_{1}$ with the conditions $0\leqslant c_{1}\leqslant D-1$, $c_{1}\ne c_{0}\leqslant D-1$, $\text{GCD}(c_{0},c_{1},D)=1$. Then there exists an index $j$, $1\leqslant j\leqslant r$ such that the set $\mathcal{A}_{r}^{\circ}(D,c_{0},c_{1})$ contains the tuples with an arbitrary large components $k_{j}$. In particular, we can choose the tuple with $k_{j}\geqslant 4r+2$.

By Lemma 14, for such tuple, one has
\[
k_{i}=1\quad \text{for}\quad |i-j|=1,\quad k_{i} = 2\quad \text{for}\quad |i-j|\geqslant 2.
\]
Moreover, this tuple satisfy to (\ref{lab_13}).

In the proof of Lemma 21, we find that in this case the length $r$ obeys to the following restrictions:
\begin{align*}
& r\leqslant \frac{D}{\text{GCD}(c_{1},D)}\quad \text{for}\quad j=1,\\
& r\leqslant \frac{D}{\text{GCD}(c_{1}-c_{0},D)}\quad \text{for}\quad j=2\quad \text{and}\quad j=r,\\
& r\leqslant \frac{2D}{\text{GCD}(c_{1},D)}\quad \text{for}\quad 3\leqslant j\leqslant r-1.
\end{align*}
In any of these cases, we have $r\leqslant 2D$. Hence, if $r>2D$ then all the sets $\mathcal{A}_{r}^{\circ}(D,c_{0},c_{1})$ are finite. From the expression for $\mathfrak{c}_{r}$ (see Lemma 23), it follows that the proportion $\nu(r,D,c_{0})$ is rational.

\pagebreak

\Large{\textbf{Appendix I. Precise form of the polygons $\boldsymbol{\mathcal{T}(\mathbf{k}_{r})}$}}
\fontsize{12}{15pt}\selectfont
\vspace{0.5cm}

For the convenience, we put in this section the list of all tuples forming the set $\mathcal{A}_{r}^{\circ}$ for any $r\geqslant 8$ and give here the description of the corresponding polygons $\mathcal{T}(\mathbf{k}_{r})$.
\vspace{0.5cm}

\begin{center}
\textbf{1.} \textsc{Case} $r = 5n+1$, $n\geqslant 2$.
\end{center}
\vspace{0.3cm}

\textbf{1.1.} $A$-type: $\mathbf{k}_{r} = \bigl(2, (4,1)^{n},3,2^{3n-2},1\bigr)$
\begin{equation*}
\begin{cases}
\displaystyle \frac{6n-2}{6n-1} < x\leqslant \frac{12n-1}{12n+1},\quad \frac{1}{3}(1+x)<y\leqslant \frac{1+4nx}{6n+1},\\[12pt]
\displaystyle \frac{12n-1}{12n+1} < x\leqslant 1,\quad \frac{1+(8n+1)x}{12n+2}<y\leqslant \frac{1+4nx}{6n+1}.
\end{cases}
\end{equation*}
Here $\mathcal{T}(\mathbf{k}_{r})$ is the triangle with the vertices
\[
A=\biggl(\frac{6n-2}{6n-1},\frac{4n-1}{6n-1}\biggr),\quad B=\biggl(1,\frac{4n+1}{6n+1}\biggr),\quad C=\biggl(\frac{12n-1}{12n+1},\frac{8n}{12n+1}\biggr).
\]

\textbf{1.2.} $B$-type: $\mathbf{k}_{r} = \bigl(5, (1,4)^{n-1},1,3,2^{3n-1},1\bigr)$
\begin{equation*}
\begin{cases}
\displaystyle \frac{12n-5}{12n-1} < x\leqslant \frac{6n-1}{6n+1},\quad \frac{1}{6}(1+x)<y\leqslant \frac{1+4nx}{12n+1},\\[12pt]
\displaystyle \frac{6n-1}{6n+1} < x\leqslant 1,\quad \frac{1+(4n+1)x}{12n+4}<y\leqslant \frac{1+4nx}{12n+1}.
\end{cases}
\end{equation*}
Here $\mathcal{T}(\mathbf{k}_{r})$ is the quadrangle with the vertices
\begin{multline*}
A=\biggl(\frac{12n-5}{12n-1},\frac{4n-1}{12n-1}\biggr),\quad B=\biggl(1,\frac{4n+1}{12n+1}\biggr),\quad C=\biggl(1,\frac{2n+1}{6n+2}\biggr),\\
D=\biggl(\frac{6n-1}{6n+1} ,\frac{2n}{6n+1}\biggr).
\end{multline*}

\textbf{1.3.} $C$-type: $\mathbf{k}_{r} = \bigl(1,2^{3n-2},3,(1,4)^{n},2\bigr)$
\begin{equation*}
\begin{cases}
\displaystyle \frac{1}{6n+1} < x\leqslant \frac{2}{12n+1},\quad 1-x<y\leqslant \frac{1}{2}(1+(6n-1)x),\\[12pt]
\displaystyle \frac{2}{12n+1} < x\leqslant \frac{1}{6n-1},\quad \frac{1}{3}(1+(12n-2)x)<y\leqslant \frac{1}{2}(1+(6n-1)x).
\end{cases}
\end{equation*}
Here $\mathcal{T}(\mathbf{k}_{r})$ is the triangle with the vertices
\[
A=\biggl(\frac{1}{6n+1},\frac{6n}{6n+1}\biggr),\quad B=\biggl(\frac{1}{6n-1},1\biggr),\quad C=\biggl(\frac{2}{12n+1},\frac{12n-1}{12n+1}\biggr).
\]

\textbf{1.4.} $C$-type: $\mathbf{k}_{r} = \bigl(1,2^{3n-1},3,(1,4)^{n-1},1,5\bigr)$
\begin{equation*}
\begin{cases}
\displaystyle \frac{1}{6n+2} < x\leqslant \frac{1}{6n+1},\quad 1-x<y\leqslant \frac{1}{3}(1+(12n+1)x),\\[12pt]
\displaystyle \frac{1}{6n+1} < x\leqslant \frac{2}{12n+1},\quad \frac{1}{3}(1+(12n-1)x)<y\leqslant \frac{1}{3}(1+(12n+1)x),\\[12pt]
\displaystyle \frac{2}{12n+1} < x\leqslant \frac{2}{12n-1},\quad \frac{1}{3}(1+(12n-1)x)<y\leqslant 1.
\end{cases}
\end{equation*}
Here $\mathcal{T}(\mathbf{k}_{r})$ is the quadrangle with the vertices
\begin{multline*}
A=\biggl(\frac{1}{6n+2},\frac{6n+1}{6n+2}\biggr),\quad B=\biggl(\frac{2}{12n+1},1\biggr),\quad C=\biggl(\frac{2}{12n-1},1\biggr),\\
D=\biggl(\frac{1}{6n+1} ,\frac{6n}{6n+1}\biggr).
\end{multline*}

\begin{center}
\textbf{2.} \textsc{Case} $r = 5n+2$, $n\geqslant 2$.
\end{center}
\vspace{0.3cm}

\textbf{2.1.} $A$-type: $\mathbf{k}_{r} = \bigl(2,(4,1)^{n},3,2^{3n-1},1\bigr)$
\begin{equation*}
\begin{cases}
\displaystyle \frac{12n-1}{12n+1} < x\leqslant \frac{6n+1}{6n+2},\quad \frac{1}{3}(1+x) <y\leqslant \frac{1+(8n+1)x}{12n+2},\\[12pt]
\displaystyle \frac{6n+1}{6n+2} < x\leqslant 1,\quad \frac{1+(8n+3)x}{12n+5}<y\leqslant \frac{1+(8n+1)x}{12n+2}.
\end{cases}
\end{equation*}
Here $\mathcal{T}(\mathbf{k}_{r})$ is the quadrangle with the vertices
\begin{multline*}
A=\biggl(\frac{12n-1}{12n+1},\frac{8n}{12n+1}\biggr),\quad B=\biggl(1,\frac{4n+1}{6n+1}\biggr),\quad C=\biggl(1,\frac{8n+4}{12n+5}\biggr),\\
D=\biggl(\frac{6n+1}{6n+2},\frac{4n+1}{6n+2}\biggr).
\end{multline*}

\textbf{2.2.} $B$-type: $\mathbf{k}_{r} = \bigl(5,(1,4)^{n-1},1,3,2^{3n},1\bigr)$
\begin{equation*}
\begin{cases}
\displaystyle \frac{6n-1}{6n+1} < x\leqslant \frac{6n+2}{6n+3},\quad \frac{1+(2n+1)x}{6n+5} <y\leqslant \frac{1+(4n+1)x}{12n+4},\\[12pt]
\displaystyle \frac{6n+2}{6n+3} < x\leqslant 1,\quad \frac{1+(4n+2)x}{12n+7}<y\leqslant \frac{1+(4n+1)x}{12n+4}.
\end{cases}
\end{equation*}
Here $\mathcal{T}(\mathbf{k}_{r})$ is the quadrangle with the vertices
\begin{multline*}
A=\biggl(\frac{6n-1}{6n+1},\frac{2n}{6n+1}\biggr),\quad B=\biggl(1,\frac{2n+1}{6n+2}\biggr),\quad C=\biggl(1,\frac{4n+3}{12n+7}\biggr),\\
D=\biggl(\frac{6n+2}{6n+3},\frac{1}{3}\biggr).
\end{multline*}

\textbf{2.3.} $C$-type: $\mathbf{k}_{r} = \bigl(1,2^{3n},3,(1,4)^{n-1},1,5\bigr)$
\begin{equation*}
\begin{cases}
\displaystyle \frac{2}{12n+7} < x\leqslant \frac{1}{6n+3},\quad 1-x <y\leqslant \frac{1}{3}(1+(12n+4)x),\\[12pt]
\displaystyle \frac{1}{6n+3} < x\leqslant \frac{1}{6n+2},\quad \frac{1}{2}(1+(6n+1)x)<y\leqslant \frac{1}{3}(1+(12n+4)x),\\[12pt]
\displaystyle \frac{1}{6n+2} < x\leqslant \frac{1}{6n+1},\quad \frac{1}{2}(1+(6n+1)x)<y\leqslant 1.
\end{cases}
\end{equation*}
Here $\mathcal{T}(\mathbf{k}_{r})$ is the quadrangle with the vertices
\begin{multline*}
A=\biggl(\frac{2}{12n+7},\frac{12n+5}{12n+7}\biggr),\quad B=\biggl(\frac{1}{6n+2},1\biggr),\quad C=\biggl(\frac{1}{6n+1},1\biggr),\\
D=\biggl(\frac{1}{6n+3},\frac{6n+2}{6n+3}\biggr).
\end{multline*}

\textbf{2.4.} $C$-type: $\mathbf{k}_{r} = \bigl(1,2^{3n-1},3,(1,4)^{n},2\bigr)$
\begin{equation*}
\begin{cases}
\displaystyle \frac{2}{12n+5} < x\leqslant \frac{1}{6n+2},\quad 1-x <y\leqslant \frac{1}{3}(1+(12n+2)x),\\[12pt]
\displaystyle \frac{1}{6n+2} < x\leqslant \frac{1}{6n+1},\quad \frac{1}{3}(1+(12n+1)x)<y\leqslant \frac{1}{3}(1+(12n+2)x),\\[12pt]
\displaystyle \frac{1}{6n+1} < x\leqslant \frac{2}{12n+1},\quad \frac{1}{3}(1+(12n+1)x)<y\leqslant 1.
\end{cases}
\end{equation*}
Here $\mathcal{T}(\mathbf{k}_{r})$ is the quadrangle with the vertices
\begin{multline*}
A=\biggl(\frac{2}{12n+5},\frac{12n+3}{12n+5}\biggr),\quad B=\biggl(\frac{1}{6n+1},1\biggr),\quad C=\biggl(\frac{2}{12n+1},1\biggr),\\
D=\biggl(\frac{1}{6n+2},\frac{6n+1}{6n+2}\biggr).
\end{multline*}

\begin{center}
\textbf{3.} \textsc{Case} $r = 5n+3$, $n\geqslant 1$.
\end{center}
\vspace{0.3cm}

\textbf{3.1.} $A$-type: $\mathbf{k}_{r} = \bigl(2,(4,1)^{n},3,2^{3n},1\bigr)$
\begin{equation*}
\begin{cases}
\displaystyle \frac{6n+1}{6n+2} < x\leqslant \frac{12n+5}{12n+7},\quad \frac{1}{3}(1+x) <y\leqslant \frac{1+(8n+3)x}{12n+5},\\[12pt]
\displaystyle \frac{12n+5}{12n+7} < x\leqslant 1,\quad \frac{1+(8n+5)x}{12n+8}<y\leqslant \frac{1+(8n+3)x}{12n+5}.
\end{cases}
\end{equation*}
Here $\mathcal{T}(\mathbf{k}_{r})$ is the quadrangle with the vertices
\begin{multline*}
A=\biggl(\frac{6n+1}{6n+2},\frac{4n+1}{6n+2}\biggr),\quad B=\biggl(1,\frac{8n+4}{12n+5}\biggr),\quad C=\biggl(1,\frac{4n+3}{6n+4}\biggr),\\
D=\biggl(\frac{12n+5}{12n+7},\frac{8n+4}{12n+7}\biggr).
\end{multline*}

\textbf{3.2.} $B$-type: $\mathbf{k}_{r} = \bigl(5,(1,4)^{n},1,3,2^{3n-1},1\bigr)$
\begin{equation*}
\begin{cases}
\displaystyle \frac{6n-1}{6n+1} < x\leqslant \frac{12n+1}{12n+5},\quad \frac{1}{6}(1+x) <y\leqslant \frac{1+(2n+1)x}{6n+5},\\[12pt]
\displaystyle \frac{12n+1}{12n+5} < x\leqslant \frac{6n+2}{6n+3},\quad \frac{1+(4n+2)x}{12n+7}<y\leqslant \frac{1+(2n+1)x}{6n+5}.
\end{cases}
\end{equation*}
Here $\mathcal{T}(\mathbf{k}_{r})$ is the triangle with the vertices
\[
A=\biggl(\frac{6n-1}{6n+1},\frac{2n}{6n+1}\biggr),\quad B=\biggl(\frac{6n+2}{6n+3},\frac{1}{3}\biggr),\quad C=\biggl(\frac{12n+1}{12n+5},\frac{4n+1}{12n+5}\biggr).
\]

\textbf{3.3.} $B$-type: $\mathbf{k}_{r} = \bigl(5,(1,4)^{n-1},1,3,2^{3n+1},1\bigr)$
\vspace{0.3cm}
\[
\frac{6n+2}{6n+3}<x\leqslant 1,\quad \frac{1+(2n+1)x}{6n+5}<y\leqslant \frac{1+(4n+2)x}{12n+7}.
\]
Here $\mathcal{T}(\mathbf{k}_{r})$ is the triangle with the vertices
\[
A=\biggl(\frac{6n+2}{6n+3},\frac{1}{3}\biggr),\quad B=\biggl(1,\frac{4n+3}{12n+7}\biggr),\quad C=\biggl(1,\frac{2n+2}{6n+5}\biggr).
\]

\textbf{3.4.} $C$-type: $\mathbf{k}_{r} = \bigl(1,2^{3n},3,(1,4)^{n},2\bigr)$
\begin{equation*}
\begin{cases}
\displaystyle \frac{1}{6n+4} < x\leqslant \frac{2}{12n+7},\quad 1-x <y\leqslant \frac{1}{3}(1+(12n+5)x),\\[12pt]
\displaystyle \frac{2}{12n+7} < x\leqslant \frac{2}{12n+5},\quad \frac{1}{3}(1+(12n+4)x)<y\leqslant \frac{1}{3}(1+(12n+5)x),\\[12pt]
\displaystyle \frac{2}{12n+5} < x\leqslant \frac{1}{6n+2},\quad \frac{1}{3}(1+(12n+4)x)<y\leqslant 1.
\end{cases}
\end{equation*}
Here $\mathcal{T}(\mathbf{k}_{r})$ is the quadrangle with the vertices
\begin{multline*}
A=\biggl(\frac{1}{6n+4},\frac{6n+3}{6n+4}\biggr),\quad B=\biggl(\frac{2}{12n+5},1\biggr),\quad C=\biggl(\frac{1}{6n+2},1\biggr),\\
D=\biggl(\frac{2}{12n+7},\frac{12n+5}{12n+7}\biggr).
\end{multline*}

\textbf{3.5.} $C$-type: $\mathbf{k}_{r} = \bigl(1,2^{3n+1},3,(1,4)^{n-1},1,5\bigr)$
\begin{equation*}
\begin{cases}
\displaystyle \frac{1}{6n+5} < x\leqslant \frac{2}{12n+7},\quad \frac{1}{2}(1+(6n+3)x) <y\leqslant \frac{1}{3}(1+(12n+7)x),\\[12pt]
\displaystyle \frac{2}{12n+7} < x\leqslant \frac{1}{6n+3},\quad \frac{1}{2}(1+(6n+3)x)<y\leqslant 1.
\end{cases}
\end{equation*}
Here $\mathcal{T}(\mathbf{k}_{r})$ is the triangle with the vertices
\[
A=\biggl(\frac{1}{6n+5},\frac{6n+4}{6n+5}\biggr),\quad B=\biggl(\frac{2}{12n+7},1\biggr),\quad C=\biggl(\frac{1}{6n+3},1\biggr).
\]

\textbf{3.6.} $C$-type: $\mathbf{k}_{r} = \bigl(1,2^{3n-1},3,(1,4)^{n},1,5\bigr)$
\begin{equation*}
\begin{cases}
\displaystyle \frac{1}{6n+3} < x\leqslant \frac{2}{12n+5},\quad 1-x <y\leqslant \frac{1}{2}(1+(6n+1)x),\\[12pt]
\displaystyle \frac{2}{12n+5} < x\leqslant \frac{1}{6n+1},\quad \frac{1}{3}(1+(12n+2)x)<y\leqslant \frac{1}{2}(1+(6n+1)x).
\end{cases}
\end{equation*}
Here $\mathcal{T}(\mathbf{k}_{r})$ is the triangle with the vertices
\[
A=\biggl(\frac{1}{6n+3},\frac{6n+2}{6n+3}\biggr),\quad B=\biggl(\frac{1}{6n+1},1\biggr),\quad C=\biggl(\frac{2}{12n+5},\frac{12n+3}{12n+5}\biggr).
\]
\vspace{0.3cm}

\begin{center}
\textbf{4.} \textsc{Case} $r = 5n+4$, $n\geqslant 1$.
\end{center}
\vspace{0.3cm}

\textbf{4.1.} $A$-type: $\mathbf{k}_{r} = \bigl(2,(4,1)^{n},3,2^{3n+1},1\bigr)$
\begin{equation*}
\begin{cases}
\displaystyle \frac{12n+5}{12n+7} < x\leqslant \frac{6n+4}{6n+5},\quad \frac{1}{3}(1+x) <y\leqslant \frac{1+(8n+5)x}{12n+8},\\[12pt]
\displaystyle \frac{6n+4}{6n+5} < x\leqslant 1,\quad \frac{1+(8n+7)x}{12n+11}<y\leqslant \frac{1+(8n+5)x}{12n+8}.
\end{cases}
\end{equation*}
Here $\mathcal{T}(\mathbf{k}_{r})$ is the quadrangle with the vertices
\begin{multline*}
A=\biggl(\frac{12n+5}{12n+7},\frac{8n+4}{12n+7}\biggr),\quad B=\biggl(1,\frac{4n+3}{6n+4}\biggr),\quad C=\biggl(1,\frac{8n+8}{12n+11}\biggr),\\
D=\biggl(\frac{6n+4}{6n+5},\frac{4n+3}{6n+5}\biggr).
\end{multline*}

\textbf{4.2.} $B$-type: $\mathbf{k}_{r} = \bigl(5,(1,4)^{n},1,3,2^{3n},1\bigr)$
\begin{equation*}
\begin{cases}
\displaystyle \frac{12n+1}{12n+5} < x\leqslant \frac{3n+1}{3n+2},\quad \frac{1}{6}(1+x) <y\leqslant \frac{1+(4n+2)x}{12n+7},\\[12pt]
\displaystyle \frac{3n+1}{3n+2} < x\leqslant \frac{6n+2}{6n+3},\quad \frac{1+(4n+3)x}{12n+10} <y\leqslant \frac{1+(4n+2)x}{12n+7},\\[12pt]
\displaystyle \frac{6n+2}{6n+3} < x\leqslant 1,\quad \frac{1+(4n+3)x}{12n+10} <y\leqslant \frac{1+(2n+1)x}{6n+5}.
\end{cases}
\end{equation*}
Here $\mathcal{T}(\mathbf{k}_{r})$ is the quadrangle with the vertices
\begin{multline*}
A=\biggl(\frac{12n+1}{12n+5},\frac{4n+1}{12n+5}\biggr),\quad B=\biggl(\frac{6n+2}{6n+3},\frac{1}{3}\biggr),\quad C=\biggl(1,\frac{2n+2}{6n+5}\biggr),\\
D=\biggl(\frac{3n+1}{3n+2},\frac{2n+1}{6n+4}\biggr).
\end{multline*}

\textbf{4.3.} $C$-type: $\mathbf{k}_{r} = \bigl(1,2^{3n},3,(1,4)^{n},1,5\bigr)$
\begin{equation*}
\begin{cases}
\displaystyle \frac{1}{6n+5} < x\leqslant \frac{1}{6n+4},\quad 1-x <y\leqslant \frac{1}{2}(1+(6n+3)x),\\[12pt]
\displaystyle \frac{1}{6n+4} < x\leqslant \frac{1}{6n+3},\quad \frac{1}{3}(1+(12n+5)x)<y\leqslant \frac{1}{2}(1+(6n+3)x),\\[12pt]
\displaystyle \frac{1}{6n+3} < x\leqslant \frac{2}{12n+5},\quad \frac{1}{3}(1+(12n+5)x)<y\leqslant 1.
\end{cases}
\end{equation*}
Here $\mathcal{T}(\mathbf{k}_{r})$ is the quadrangle with the vertices
\begin{multline*}
A=\biggl(\frac{1}{6n+5},\frac{6n+4}{6n+5}\biggr),\quad B=\biggl(\frac{1}{6n+3},1\biggr),\quad C=\biggl(\frac{2}{12n+5},1\biggr),\\
D=\biggl(\frac{1}{6n+4},\frac{6n+3}{6n+4}\biggr).
\end{multline*}

\textbf{4.4.} $C$-type: $\mathbf{k}_{r} = \bigl(1,2^{3n+1},3,(1,4)^{n},2\bigr)$
\begin{equation*}
\begin{cases}
\displaystyle \frac{2}{12n+11} < x\leqslant \frac{1}{6n+5},\quad 1-x <y\leqslant \frac{1}{3}(1+(12n+8)x),\\[12pt]
\displaystyle \frac{1}{6n+5} < x\leqslant \frac{1}{6n+4},\quad \frac{1}{3}(1+(12n+7)x)<y\leqslant \frac{1}{3}(1+(12n+8)x),\\[12pt]
\displaystyle \frac{1}{6n+4} < x\leqslant \frac{2}{12n+7},\quad \frac{1}{3}(1+(12n+7)x)<y\leqslant 1.
\end{cases}
\end{equation*}
Here $\mathcal{T}(\mathbf{k}_{r})$ is the quadrangle with the vertices
\begin{multline*}
A=\biggl(\frac{2}{12n+11},\frac{12n+9}{12n+11}\biggr),\quad B=\biggl(\frac{1}{6n+4},1\biggr),\quad C=\biggl(\frac{2}{12n+7},1\biggr),\\
D=\biggl(\frac{1}{6n+5},\frac{6n+4}{6n+5}\biggr).
\end{multline*}

\begin{center}
\textbf{5.} \textsc{Case} $r = 5n+5$, $n\geqslant 1$.
\end{center}
\vspace{0.3cm}

\textbf{5.1.} $A$-type: $\mathbf{k}_{r} = \bigl(2, (4,1)^{n},3,2^{3n+2},1\bigr)$
\vspace{0.3cm}
\[
\frac{6n+4}{6n+5}<x\leqslant 1,\quad \frac{1+(4n+4)x}{6n+7}<y\leqslant \frac{1+(8n+7)x}{12n+11}.
\]
Here $\mathcal{T}(\mathbf{k}_{r})$ is the triangle with the vertices
\[
A=\biggl(\frac{6n+4}{6n+5},\frac{4n+3}{6n+5}\biggr),\quad B=\biggl(1,\frac{8n+8}{12n+11}\biggr),\quad C=\biggl(1,\frac{4n+5}{6n+7}\biggr).
\]

\textbf{5.2.} $B$-type: $\mathbf{k}_{r} = \bigl(5, (1,4)^{n},1,3,2^{3n+1},1\bigr)$
\vspace{0.3cm}

\begin{equation*}
\begin{cases}
\displaystyle \frac{3n+1}{3n+2} < x\leqslant \frac{12n+7}{12n+11},\quad \frac{1}{6}(1+x)<y\leqslant \frac{1+(4n+3)x}{12n+10},\\[12pt]
\displaystyle \frac{12n+7}{12n+11} < x\leqslant 1,\quad \frac{1+(4n+4)x}{12n+13}<y\leqslant \frac{1+(4n+3)x}{12n+10}.
\end{cases}
\end{equation*}
Here $\mathcal{T}(\mathbf{k}_{r})$ is the quadrangle with the vertices
\begin{multline*}
A=\biggl(\frac{3n+1}{3n+2},\frac{2n+1}{6n+4}\biggr),\quad B=\biggl(1,\frac{2n+2}{6n+5}\biggr),\quad C=\biggl(1,\frac{4n+5}{12n+13}\biggr),\\
D=\biggl(\frac{12n+7}{12n+11} ,\frac{4n+3}{12n+11}\biggr).
\end{multline*}

\textbf{5.3.} $C$-type: $\mathbf{k}_{r} = \bigl(1,2^{3n+1},3,(1,4)^{n},1,5\bigr)$
\vspace{0.3cm}
\begin{equation*}
\begin{cases}
\displaystyle \frac{2}{12n+13} < x\leqslant \frac{2}{12n+11},\quad 1-x<y\leqslant \frac{1}{3}(1+(12n+10)x),\\[12pt]
\displaystyle \frac{2}{12n+11} < x\leqslant \frac{1}{6n+5},\quad \frac{1}{3}(1+(12n+8)x)<y\leqslant \frac{1}{3}(1+(12n+10)x),\\[12pt]
\displaystyle \frac{1}{6n+5} < x\leqslant \frac{1}{6n+4},\quad \frac{1}{3}(1+(12n+8)x)<y\leqslant 1.
\end{cases}
\end{equation*}
Here $\mathcal{T}(\mathbf{k}_{r})$ is the quadrangle with the vertices
\begin{multline*}
A=\biggl(\frac{2}{12n+13},\frac{12n+11}{12n+13}\biggr),\quad B=\biggl(\frac{1}{6n+5},1\biggr),\quad C=\biggl(\frac{1}{6n+4},1\biggr),\\
D=\biggl(\frac{2}{12n+11} ,\frac{12n+9}{12n+11}\biggr).
\end{multline*}

\textbf{5.4.} $C$-type: $\mathbf{k}_{r} = \bigl(1,2^{3n+2},3,(1,4)^{n},2\bigr)$
\vspace{0.3cm}
\begin{equation*}
\begin{cases}
\displaystyle \frac{1}{6n+7} < x\leqslant \frac{2}{12n+11},\quad \frac{1}{2}(1+(6n+5)x)<y\leqslant \frac{1}{3}(1+(12n+11)x),\\[12pt]
\displaystyle \frac{2}{12n+11} < x\leqslant \frac{1}{6n+5},\quad \frac{1}{2}(1+(6n+5)x)<y\leqslant 1.
\end{cases}
\end{equation*}
Here $\mathcal{T}(\mathbf{k}_{r})$ is the triangle with the vertices
\[
A=\biggl(\frac{1}{6n+7},\frac{6n+6}{6n+7}\biggr),\quad B=\biggl(\frac{2}{12n+11},1\biggr),\quad C=\biggl(\frac{1}{6n+5},1\biggr).
\]

\begin{center}
\includegraphics{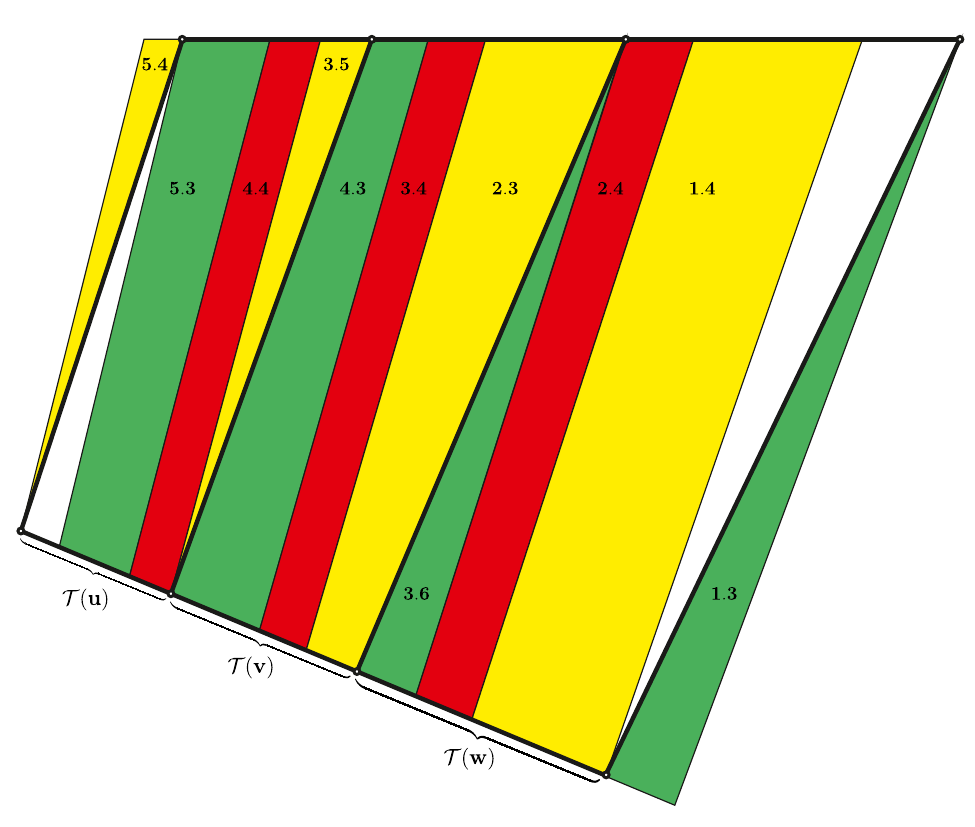}
\vspace{0.2cm}

\fontsize{10}{12pt}\selectfont
\emph{Fig.~13.} The regions corresponding to the tuples of C-type with the same $n$ from the Appendix I. \\
The bold lines bound the quadrangles $\mathcal{T}(\mathbf{u})$, $\mathcal{T}(\mathbf{v})$ and $\mathcal{T}(\mathbf{w})$ with \\ $\mathbf{u} = (1,2^{3n+1},3,(1,4)^{n-1},1)$, $\mathbf{v} = (1,2^{3n},3,(1,4)^{n-1},1)$ \\
and $\mathbf{w} = (1,2^{3n-1},2,(1,4)^{n-1},1)$.
\fontsize{12}{15pt}\selectfont
\end{center}

\pagebreak

\Large{\textbf{Appendix II. Table of the proportions $\boldsymbol{\nu(r;3,0), r\geqslant 8}$}}
\vspace{0.3cm}
\fontsize{12}{15pt}\selectfont

\begin{longtable}{|>{\fontsize{10}{9pt}\selectfont}c|>{\fontsize{10}{9pt}\selectfont}c|>{\fontsize{10}{9pt}\selectfont}c|}
\hline $r$ & $\nu(r;3,0)$       & Approximate value  \\
\hline 8 & 12\,797 / 2\,238\,390 & 0.00571705556225680\\
 \hline 9 & 18\,662 / 3\,677\,355 & 0.00507484319572084\\
 \hline 10 & 284 / 82\,225 & 0.00345393736698085\\
 \hline 11 & 1\,444 / 575\,575 & 0.00250879555227381\\
 \hline 12 & 11\,402 / 6\,135\,675 & 0.00185831224763372\\
 \hline 13 & 100\,349 / 83\,445\,180 & 0.00120257395334278\\
 \hline 14 & 3\,931 / 3\,209\,430 & 0.00122482808473779\\
 \hline 15 & 1\,607 / 1\,673\,140 & 0.000960469536320930\\
 \hline 16 & 65 / 83\,657 & 0.000776982201130808\\
 \hline 17 & 82\,244 / 130\,135\,845 & 0.000631985753041370\\
 \hline 18 & 387\,197 / 889\,847\,805 & 0.000435127218187609\\
 \hline 19 & 206\,834 / 440\,240\,493 & 0.000469820480598090\\
 \hline 20 & 5\,728 / 14\,566\,475 & 0.000393231718723988\\
 \hline 21 & 5\,776 / 17\,214\,925 & 0.000335522809422638\\
 \hline 22 & 71\,846 / 250\,675\,425 & 0.000286609666663575\\
 \hline 23 & 211\,993 / 1\,038\,512\,475 & 0.000204131394762494\\
 \hline 24 & 214\,612 / 942\,649\,785 & 0.000227668857952373\\
 \hline 25 & 4\,481 / 22\,648\,507 & 0.000197849686074230\\
 \hline 26 & 4\,511 / 25\,884\,008 & 0.000174277492110186\\
 \hline 27 & 166\,789 / 1\,087\,060\,260 & 0.000153431236645520\\
 \hline 28 & 473\,497 / 4\,241\,317\,080 & 0.000111639142056316\\
 \hline 29 & 154\,342 / 1\,214\,248\,035 & 0.000127109120666603\\
 \hline 30 & 12\,916 / 114\,103\,745 & 0.000113195232987313\\
 \hline 31 & 12\,988 / 127\,527\,715 & 0.000101844528461911\\
 \hline 32 & 1\,114\,186 / 12\,174\,765\,603 & 0.0000915160124089331\\
 \hline 33 & 4\,620\,797 / 68\,378\,820\,510 & 0.0000675764361762314\\
 \hline 34 & 315\,053 / 4\,036\,587\,555 & 0.0000780493413575913\\
 \hline 35 & 1\,759 / 24\,875\,930 & 0.0000707109241744932\\
 \hline 36 & 8\,837 / 136\,817\,615 & 0.0000645896363564005\\
 \hline 37 & 172\,624 / 2\,930\,582\,655 & 0.0000589043273375956\\
 \hline 38 & 8\,193\,149 / 186\,350\,579\,415 & 0.0000439663188905572\\
 \hline 39 & 76\,826 / 1\,497\,108\,249 & 0.0000513162625690669\\
 \hline 40 & 22\,984 / 488\,109\,335 & 0.0000470878107668234\\
 \hline 41 & 4\,616 / 106\,110\,725 & 0.0000435017289722599\\
 \hline 42 & 2\,529\,238 / 63\,042\,960\,225 & 0.0000401192772511501\\
 \hline 43 & 13\,519\,997 / 447\,799\,995\,825 & 0.0000301920436044033\\
 \hline 44 & 1\,389\,472 / 39\,113\,958\,819 & 0.0000355236862223481\\
 \hline 45 & 14\,549 / 441\,969\,385 & 0.0000329185696878077\\
 \hline 46 & 14\,603 / 475\,967\,030 & 0.0000306806965179920\\
 \hline 47 & 35\,491 / 1\,243\,364\,199 & 0.0000285443316033583\\
 \hline 48 & 4\,219\,801 / 195\,173\,958\,210 & 0.0000216207174292159\\
 \hline 49 & 3\,861\,014 / 150\,816\,240\,435 & 0.0000256007840326987\\
 \hline 50 & 35\,932 / 1\,502\,838\,029 & 0.0000239094295636832\\
 \hline 51 & 36\,052 / 1\,606\,482\,031 & 0.0000224415831016538\\
 \hline 52 & 4\,811\,746 / 228\,839\,601\,375 & 0.0000210267190254146\\
 \hline 53 & 6\,298\,009 / 393\,377\,166\,000 & 0.0000160101031385233 \\
 \hline 54 & 259\,639 / 13\,625\,703\,000 & 0.0000190550902217669 \\
 \hline 55 & 21\,743 / 1\,214\,034\,640 & 0.0000179097031366420 \\
 \hline 56 & 21\,809 / 1\,289\,911\,805 & 0.0000169073574762733 \\
 \hline 57 & 3\,171\,548 / 199\,048\,918\,929 & 0.0000159335103002055 \\
 \hline 58 & 45\,315\,197 / 3\,719\,071\,906\,305 & 0.0000121845444620677 \\
 \hline 59 & 6\,800\,162 / 466\,924\,743\,285 & 0.0000145637216656333 \\
 \hline 60 & 10\,352 / 752\,286\,535 & 0.0000137607141938278\\
 \hline 61 & 51\,904 / 3\,976\,371\,685 & 0.0000130531057234404\\
 \hline 62 & 1\,633\,814 / 132\,174\,312\,825 & 0.0000123610553751332\\
 \hline 63 & 63\,258\,749 / 6\,667\,966\,195\,275 & 0.00000948696306301401\\
 \hline 64 & 870\,908 / 76\,529\,786\,025 & 0.0000113799873909944\\
 \hline 65 & 30\,377 / 2\,812\,566\,295 & 0.0000108004565275500\\
 \hline 66 & 6\,091 / 592\,119\,220 & 0.0000102867797468219\\
 \hline 67 & 2\,578\,897 / 263\,648\,462\,490 & 0.00000978157420545479\\
 \hline 68 & 86\,067\,197 / 11\,429\,244\,813\,420 & 0.00000753043603536618\\
 \hline 69 & 10\,945\,454 / 1\,208\,027\,774\,583 &  0.00000906059796826962\\
\hline 70 & 70\,468 / 8\,163\,643\,865 & 0.00000863192970753141\\
 \hline 71 & 70\,636 / 8\,561\,870\,395 & 0.00000825006648561865\\
 \hline 72 & 2\,561\,714 / 325\,392\,460\,575 & 0.00000787269009083125\\
 \hline 73 & 22\,909\,849 / 3\,769\,931\,584\,650 & 0.00000607699330494003\\
 \hline 74 & 3\,383\,801 / 461\,571\,650\,925 & 0.00000733104165565365\\
 \hline 75 & 40\,451 / 5\,772\,789\,022 & 0.00000700718488859405\\
 \hline 76 & 40\,541 / 6\,035\,188\,523 & 0.00000671743721766088\\
 \hline 77 & 7\,836\,968 / 1\,218\,832\,670\,055 & 0.00000642989656623362\\
 \hline 78 & 29\,915\,161 / 6\,013\,356\,764\,415 & 0.00000497478565998740\\
 \hline 79 & 660\,170 / 109\,751\,661\,921 & 0.00000601512531514279\\
 \hline 80 & 92\,056 / 15\,965\,549\,615 & 0.00000576591487420585\\
\hline
\end{longtable}

\vspace{0.5cm}

\Large{\textbf{Appendix III. The sets $\boldsymbol{\mathcal{A}_{r}^{\circ}(3;0,1)}$}}
\vspace{0.3cm}
\fontsize{12}{15pt}\selectfont

\begin{longtable}{|>{\fontsize{10}{9pt}\selectfont}c|>{\fontsize{10}{9pt}\selectfont}c|>{\fontsize{10}{9pt}\selectfont}c|}
\hline $r$ & $\mathbf{\overline{k}}$       & $\mathbb{K}_{r}(\mathbf{\overline{k}})$  \\
\hline 1   & $(k_{1}), k_{1}\equiv 0\pmod{3}$, $k_{1}\geqslant 3$ & $k_{1}$\\
\hline 2   & $(k_{1},1), k_{1}\equiv 1\pmod{3}$, $k_{1}\geqslant 4$ & $k_{1}-1$\\
\hline 2   & $(1,k_{2}), k_{2}\equiv 1\pmod{3}$, $k_{2}\geqslant 4$ & $k_{2}-1$\\
\hline 2   & $(2,2)$ & $3$\\
\hline 3   & $(k_{1},1,2), k_{1}\equiv 2\pmod{3}$, $k_{1}\geqslant 8$ & $k_{1}-2$\\
\hline 3   & $(1,k_{2},1), k_{2}\equiv 2\pmod{3}$, $k_{2}\geqslant 5$ & $k_{2}-2$\\
\hline 3   & $(2,1,k_{3}), k_{3}\equiv 2\pmod{3}$, $k_{3}\geqslant 8$ & $k_{3}-2$\\
\hline 3   & $(1,2,4), (1,3,2), (2,3,1), (4,2,1)$ & $3$\\
\hline 4   & $(1,k_{2},1,2), k_{2}\equiv 0\pmod{3}$, $k_{2}\geqslant 6$ & $k_{2}-3$\\
\hline 4   & $(2,1,k_{3},1), k_{3}\equiv 0\pmod{3}$, $k_{3}\geqslant 6$ & $k_{3}-3$\\
\hline 4   & $(1,2,3,2), (1,3,1,5), (2,3,2,1), (5,1,3,1)$ & $3$\\
\hline 4   & $(1,3,1,8), (8,1,3,1)$ & $6$\\
\hline 5   & $(2,1,k_{3},1,2), k_{3}\equiv 1\pmod{3}$, $k_{3}\geqslant 7$ & $k_{3}-4$\\
\hline 5   & $(1,2,2,3,2), (2,3,2,2,1)$ & $3$\\
\hline 5   & $(1,2,3,1,5), (5,1,3,2,1)$ & $3$\\
\hline 5   & $(1,3,1,6,1), (1,6,1,3,1)$ & $3$\\
\hline 6   & $(1,2,2,3,1,5), (5,1,3,2,2,1)$ & $3$\\
\hline 6   & $(1,2,3,1,6,1), (1,6,1,3,2,1)$ & $3$\\
\hline 6   & $(1,2,3,1,4,2), (2,4,1,3,2,1)$ & $3$\\
\hline 6   & $(1,3,1,7,1,2), (2,1,7,1,3,1)$ & $3$\\
\hline 7   & $(1,2,2,2,3,1,5), (5,1,3,2,2,2,1)$ & $3$\\
\hline 7   & $(1,2,2,3,1,6,1), (1,6,1,3,2,2,1)$ & $3$\\
\hline 7   & $(1,2,2,3,1,4,2), (2,4,1,3,2,2,1)$ & $3$\\
\hline 7   & $(1,3,1,7,1,3,1)$ & $3$\\
\hline $5m,\;m\geqslant 2$   & $\bigl(2,(4,1)^{m-1},3,2^{3m-1},1\bigr), \bigl(1,2^{3m-1},3,(1,4)^{m-1},2\bigr)$ & $3$\\
\hline $5m, \;m\geqslant 2$   & $\bigl(5,(1,4)^{m-1},1,3,2^{3m-2},1\bigr), \bigl(1,2^{3m-2},3,(1,4)^{m-1},1,5\bigr)$ & $3$\\
\hline $5m+1, \;m\geqslant 2$   & $\bigl(2,(4,1)^{m},3,2^{3m-2},1\bigr), \bigl(1,2^{3m-2},3,(1,4)^{m},2\bigr)$ & $3$\\
\hline $5m+1, \;m\geqslant 2$   & $\bigl(5,(1,4)^{m-1},1,3,2^{3m-1},1\bigr), \bigl(1,2^{3m-1},3,(1,4)^{m-1},1,5\bigr)$ & $3$\\
\hline $5m+2, \;m\geqslant 2$   & $\bigl(2,(4,1)^{m},3,2^{3m-1},1\bigr), \bigl(1,2^{3m-1},3,(1,4)^{m},2\bigr)$ & $3$\\
\hline $5m+2, \;m\geqslant 2$   & $\bigl(5,(1,4)^{m-1},1,3,2^{3m},1\bigr), \bigl(1,2^{3m},3,(1,4)^{m-1},1,5\bigr)$ & $3$\\
\hline $5m+3, \;m\geqslant 1$   & $\bigl(2,(4,1)^{m},3,2^{3m},1\bigr), \bigl(1,2^{3m},3,(1,4)^{m},2\bigr)$ & $3$\\
\hline $5m+3, \;m\geqslant 1$   & $\bigl(5,(1,4)^{m},1,3,2^{3m-1},1\bigr), \bigl(1,2^{3m-1},3,(1,4)^{m},1,5\bigr)$ & $3$\\
\hline $5m+3, \;m\geqslant 1$   & $\bigl(5,(1,4)^{m-1},1,3,2^{3m+1},1\bigr), \bigl(1,2^{3m+1},3,(1,4)^{m-1},1,5\bigr)$ & $3$\\
\hline $5m+4, \;m\geqslant 1$   & $\bigl(2,(4,1)^{m},3,2^{3m+1},1\bigr), \bigl(1,2^{3m+1},3,(1,4)^{m},2\bigr)$ & $3$\\
\hline $5m+4, \;m\geqslant 1$   & $\bigl(5,(1,4)^{m},1,3,2^{3m},1\bigr), \bigl(1,2^{3m},3,(4,1)^{m},1,5\bigr)$ & $3$\\
\hline
\end{longtable}

\renewcommand{\refname}{\Large{References}}


\begin{thebibliography}{99}
\fontsize{11}{13pt}\selectfont
	
\bibitem{Cobeli_Zaharescu_2003}
\textsc{C.~Cobeli}, \textsc{A.~Zaharescu}, The Haros-Farey sequence at two hundred years. A survey. \emph{Acta Univ. Apulensis Math. Inform.}, \textbf{5} (2003), 1--38.

\bibitem{Haynes_2003}
\textsc{A.~Haynes}, A note on Farey fractions with odd denominators. \emph{J. Number Theory}, \textbf{98} (2003), №~2, 89--104.

\bibitem{Boca_Cobeli_Zaharescu_2003}
\textsc{F.P.~Boca}, \textsc{C.~Cobeli}, \textsc{A.~Zaharescu}, On the distribution of the Farey sequence with odd denominators. \emph{Michigan Math. J.}, \textbf{51} (2003), 557--573.

\bibitem{Cobeli_Zaharescu_2005}
\textsc{C.~Cobeli}, \textsc{A.~Zaharescu}, The distribution of rationals in residue classes. \texttt{arXiv:math/0511356v1 [math.NT]} (14.11.2005).

\bibitem{Cobeli_Zaharescu_2006}
\textsc{C.~Cobeli}, \textsc{A.~Zaharescu}, On the Farey Fractions  with Denominators in Arithmetic Progression. \emph{J. Integer Sequences}, \textbf{9} (2006), Article 06.3.4.

\bibitem{Alkan_Ledoan_Vajaitu_Zaharescu_2006}
\textsc{E.~Alkan}, \textsc{A.H.~Ledoan}, \textsc{M.~V\^{a}j\^{a}itu}, \textsc{A.~Zaharescu}, Discrepancy of Fractions with Divisibility Constraints. \emph{Monatsh. Math.}, \textbf{149} (2006), 179--192.

\bibitem{Haynes_2009}
\textsc{A.~Haynes}, The distribution of special subsets of the Farey sequence. \emph{J. Number Theory}, \textbf{107} (2004) 95--104. See also: \texttt{arXiv:math/0907.2171v1 [math.NT]} (13.07.2009).

\bibitem{Cobeli_Vajaitu_Zaharescu_2012}
\textsc{C.~Cobeli}, \textsc{M.~V\^{a}j\^{a}itu}, \textsc{A.~Zaharescu}, The distribution of rationals in residue classes. \emph{Math. Reports.}, \textbf{14(64)} (2012), №~1, 1--19.

\bibitem{Badziahin_Haynes_2011}
\textsc{D.A.~Badziahin}, \textsc{A.K.~Haynes}, A note on Farey fractions with denominators in arithmetic progressions. \emph{Acta Arith.}, \textbf{146} (2011), №~3, 205--215.

\bibitem{Boca_Heersink_Spiegelhalter_2013}
\textsc{F.P.~Boca}, \textsc{B.~Heersink}, \textsc{P.~Spiegelhalter}, Gap distribution of Farey fractions under some divisibility constraints. \emph{Integers}, \textbf{13} (2013), 634--648. See also: \texttt{arXiv:1301.0277v2 [math.NT]} (11.04.2013)

\bibitem{Chahal_Chaubey_Goel_2024}
\textsc{B.~Chahal}, \textsc{S.~Chaubey}, \textsc{S.~Goel}, On the distribution of index of Farey Sequences. \emph{Res. Number Theory}, \textbf{10} (2024), Article 27, 34 pp. See also: \texttt{arXiv:2203.16215v1 [math.NT]} (30.03.2022).

\bibitem{Boca_Siskaki_2022}
\textsc{F.P.~Boca}, \textsc{M.~Siskaki}, A note on the pair correlation of Farey fractions. \emph{Acta Arith.}, \textbf{205} (2022), 121--135. See also: \texttt{arXiv:2109.12744v2 [math.NT]} (19.09.2022)

\bibitem{Boca_Cobeli_Zaharescu_2001}
\textsc{F.P.~Boca}, \textsc{C.~Cobeli}, \textsc{A.~Zaharescu}, A conjecture of R. R. Hall on Farey points. \emph{J. reine angew. Math.}, \textbf{535} (2001), 207--236.

\bibitem{Baker_1975}
\textsc{A.~Baker}, Transcendental number theory. Cambr. Univ. Press, 1975.

\bibitem{Lindeman_1882}
\textsc{F.~Lindeman}, Ueber die Zahl $\pi$. \emph{Math. Ann.}, \textbf{20} (1882), 213--225.

\bibitem{Korolev_2023}
\textsc{M.A.~Korolev}, A distribution related to Farey series. \emph{Chebyshevski\v{i} Sb.}, \textbf{24} (2023), №4, 137--190.

\bibitem{Ustinov_2009}
\textsc{A.V.~Ustinov}, The solution of Arnold's problem on the weak asymptotics of Frobenius numbers with three arguments. \emph{Sb. Math.}, \textbf{200} (2009), №~4, 597--627.

\bibitem{Graham_Knuth_Patashnik_1998}
\textsc{R.L.~Graham}, \textsc{D.E.~Knuth}, \textsc{O.~Patashnik}, Concrete Mathematics. 2nd ed., Massachusetts, Addison-Wesley, 1994.

\bibitem{Smirnov_2022}
\textsc{E.Yu.~Smirnov},  Frieze patterns and continued fractions (in Russian), M., MCCME,  2022.

\bibitem{Boca_Gologan_Zaharescu_2002}
\textsc{F.P.~Boca}, \textsc{R.~Gologan}, \textsc{A.~Zaharescu}, On the index of Farey sequences. \emph{Quart. J. Math.}, \textbf{53} (2002), №~4, 377--391.

\bibitem{Postnikov_1971}
\textsc{A.G.~Postnikov}, Introduction to analytic number theory. Transl. Math. Monographs. Vol. 68. AMS, Providence, Rhode Island. 1988.

\bibitem{Ryjik_Gradstein_1962}
\textsc{I.S.~Gradshteyn}, \textsc{I.M.~Ryzhik}, \textsc{Yu.V.~Geronimus}, \textsc{M.Yu.~Tseytlin}, Table of Integrals, Series, and Products. 8th ed. Academic Press, 2015.

\end{thebibliography}
\end{document}